\documentclass{article}
\usepackage{graphicx} 
\usepackage[letterpaper,top=2cm,bottom=2cm,left=3cm,right=3cm,marginparwidth=1.75cm]{geometry}
\usepackage{authblk}

\usepackage[algo2e,ruled,vlined]{algorithm2e}
\usepackage[normalem]{ulem}
\usepackage[noadjust]{cite}
\usepackage{algorithm}
\usepackage{booktabs}
\usepackage{verbatim}
\usepackage{algpseudocode}
\usepackage{amsmath, amsfonts, amsthm,amssymb}
\usepackage{graphicx,soul}
\usepackage{makecell}
\usepackage[export]{adjustbox}
\usepackage{subcaption}
\usepackage[colorlinks=true, allcolors=blue]{hyperref}
\usepackage{appendix}
\newtheorem{theorem}{Theorem}[section]

\newtheorem{proposition}{Proposition}[section]
\usepackage{amsfonts}
\usepackage{color,xcolor}
\theoremstyle{definition}

\usepackage{tabulary} 
\usepackage{bm}
\usepackage{tabularray}
\usepackage{booktabs}
\usepackage{diagbox}
\usepackage{multirow}
\usepackage{array}
\usepackage{amsmath,amssymb}
\usepackage{caption}
\usepackage{soul}
\theoremstyle{remark}

\newtheorem{remark}[theorem]{Remark}

\newcommand{\bc}{\mathbf{c}}
\newcommand{\be}{\mathbf{e}}

\newcommand{\cC}{\mathcal{C}}

\newcommand{\cD}{\mathcal{D}}

\newcommand{\cO}{\mathcal{O}}

\newcommand{\cF}{\mathcal{F}}

\newcommand{\by}{\mathbf{y}}

\newcommand{\PP}{\mathbb P}

\newcommand{\bF}{\mathbf F}

\usepackage{hyperref}

\DeclareMathOperator*{\argmax}{arg\,max}
\DeclareMathOperator*{\argmin}{arg\,min}

\title{Phase-IDENT: Identification of Two-phase PDEs with Uncertainty Quantification}
\author{Edward L. Yang}
\author{Roy Y. He\thanks{The research of Roy Y. He is partially supported by NSFC grant 12501594, PROCORE-France/Hong Kong Joint Research Scheme by the RGC of Hong Kong and the Consulate General of France in Hong Kong (F-CityU101/24), StUp - CityU 7200779 from City University of Hong Kong, and the Hong Kong Research Grant
Council ECS grant 21309625.}}
\affil{Department of Mathematics, City University of Hong Kong, Kowloon Tong, Hong Kong}
\date{}
\begin{document}

\maketitle
\begin{abstract}
    We propose a novel method, Phase-IDENT, for identifying partial differential equations (PDEs) from noisy observations of dynamical systems that exhibit phase transitions. Such phenomena are prevalent in fluid dynamics and materials science, where they can be modeled mathematically as functions satisfying different PDEs within distinct regions separated by phase boundaries. Our approach simultaneously identifies the underlying PDEs in each regime and accurately reconstructs the phase boundaries. Furthermore, by incorporating change point detection techniques, we provide uncertainty quantification for the detected boundaries, enhancing the interpretability and robustness of our method. We conduct numerical experiments on a variety of two-phase PDE systems under different noise levels, and the results demonstrate the effectiveness of the proposed approach.
\end{abstract}

\section{Introduction}

Data-driven partial differential equation (PDE) discovery is a rapidly growing research area that has attracted considerable attention across academic and industrial domains. It aims to discover the symbolic expression of a PDE that governs observed experimental dynamics, thereby facilitating automated system modeling~\cite{christofides2002nonlinear,jansson2005computational}, principled model selection~\cite{thanasutives2024adaptive,gerardos2025principled}, and novel scientific discovery~\cite{chen2022symbolic,berg2019data}. Among the various techniques developed for this purpose~\cite{karniadakis2021physics,lucia2004reduced,koza1994genetic}, sparsity-based PDE identification methods~\cite{rudy2017data,kang2021ident,he2022robust,he2024much,messenger2021weak,tang2023weakident,tang2023fourier} have shown particular promise.

\begin{figure}
    \centering
    \begin{subfigure}[t]{0.33\textwidth}
        \centering
        \includegraphics[width=1\linewidth]{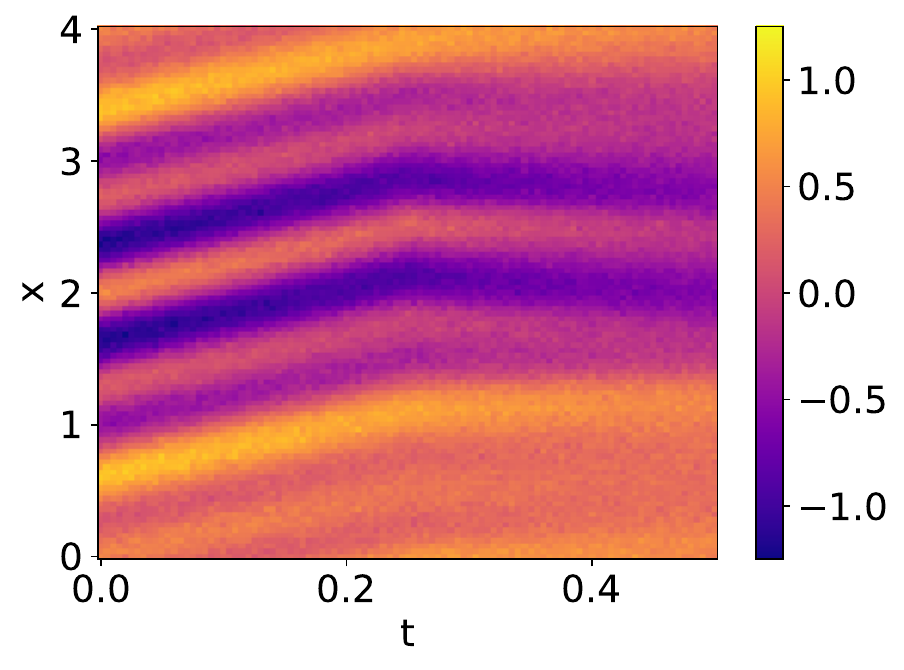}
        \caption{}
    \end{subfigure}%
    \begin{subfigure}[t]{0.33\textwidth}
        \centering
        \includegraphics[width=1\linewidth]{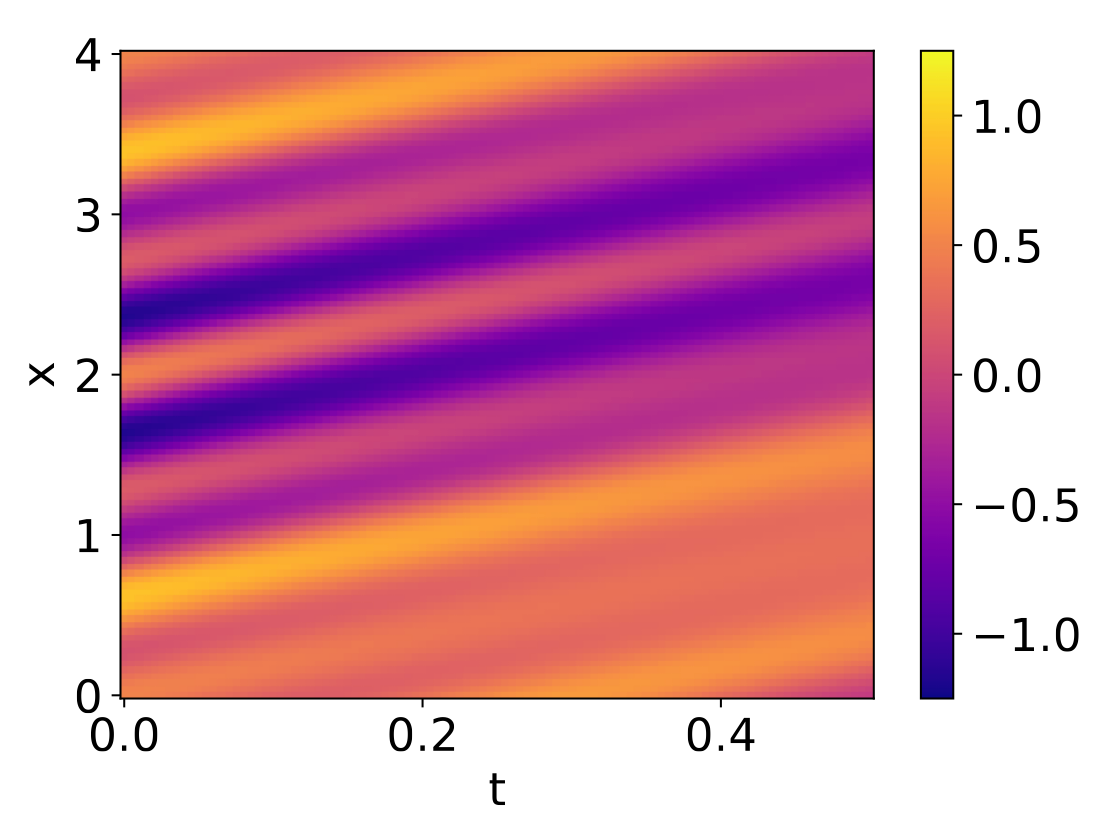}
        \caption{}
    \end{subfigure}%
    \begin{subfigure}[t]{0.33\textwidth}
        \centering
        \includegraphics[width=1\linewidth]{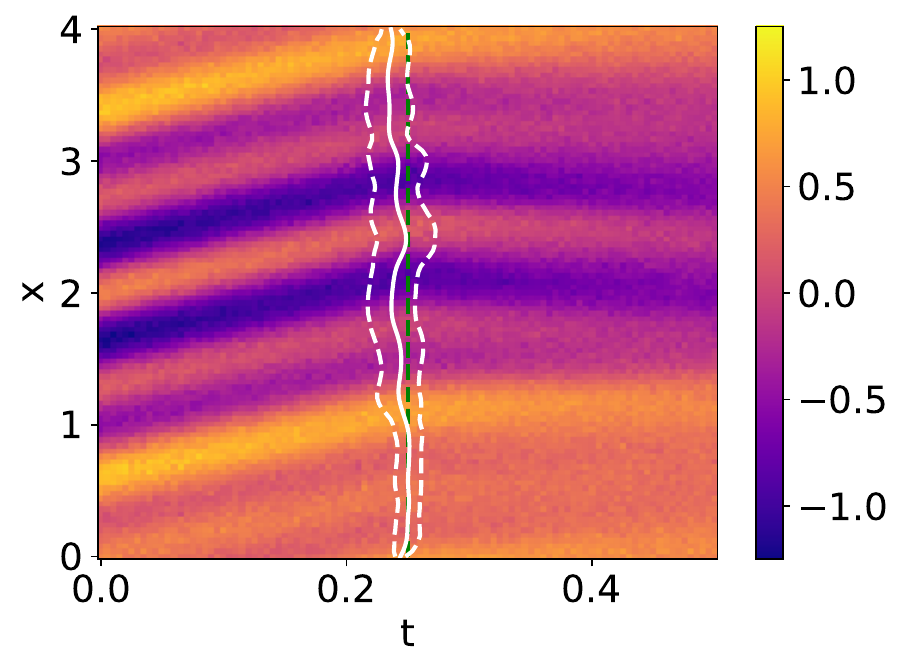}
        \caption{}
    \end{subfigure}
    \caption{Noisy observations of a solution governed by two distinct PDEs in separate domains. (a) A simulated solution governed by the transport equation followed by the viscous Burgers equation, with $10\%$ additive noise. (b) The numerical solution of a globally identified PDE using the data from (a). (c) The green dashed line indicates the ground truth location of the phase boundary. The white solid line shows the phase boundary estimated by Phase-IDENT, while the white dashed lines represent the corresponding confidence intervals at an $80\%$ confidence level (Section~\ref{boundary estimation}).}
    \label{kdvbur planform}
\end{figure}

Pioneered by Brunton et al.~\cite{brunton2016discovering}, the Sparse Identification of Nonlinear Dynamics (SINDy) framework demonstrates the feasibility of identifying dynamical systems via sparse regression. Unlike classical sparse signal recovery~\cite{candes2006stable}, SINDy constructs its dictionary using numerical differentiation and nonlinear transformations of the given data. The regression algorithm in SINDy, Sequential Threshold Ridge Regression (STRidge), was later adapted by Rudy et al.~\cite{rudy2017data} to PDE-FIND for identifying PDEs, including the viscous Burgers, reaction-diffusion, and Navier–Stokes equations. Approaching the problem from a numerical PDE perspective, Kang et al.~\cite{kang2021ident} proposed IDENT, which combines $\ell_1$-regularization with model selection based on time evolution error (TEE). To mitigate noise amplification from numerical differentiation, Robust-IDENT~\cite{he2022robust} was introduced, where the authors also found $\ell_0$-constrained regression, implemented via the Subspace Pursuit (SP) algorithm~\cite{dai2009subspace}, to be more efficient. The PDE-FIND framework continues to inspire variants such as Weak-SINDy~\cite{messenger2021weak}, SINDy-PI~\cite{kaheman2020sindy}, and Ensemble-SINDy~\cite{fasel2022ensemble}. Meanwhile, Robust-IDENT has evolved further, enabling the identification of PDEs with varying coefficients (GP-IDENT)~\cite{he2023group}, local patch-based identification (CaSLR)~\cite{he2024much}, adaptations to weak formulations (Weak-IDENT and WG-IDENT)~\cite{tang2023weakident,tang2025wg}, and Fourier-based representations (Fourier-IDENT)~\cite{tang2023fourier}, thereby addressing diverse physical constraints and data modalities. More recently, Stoch-IDENT~\cite{cui2025stoch} was proposed for identifying stochastic PDEs. For a recent survey, we refer readers to~\cite{he2025ident}. We also note that another stream of development is based on deep learning~\cite{raissi2017machine, raissi2019physics, long2019pde, kim2020integration, du2024discover, chen2024invariance}.

Notably, the aforementioned methods primarily focus on identifying single PDE models from observed trajectories. However, many physical and biological systems exhibit more complex dynamical behaviors, such as phase transitions~\cite{domb2000phase} and regime shifts~\cite{andersen2009ecological}, that require modeling beyond single-equation approaches. For instance, fluid properties can change significantly with variations in the Reynolds number~\cite{schlichting2016boundary}. Moreover, different PDEs may coexist within a common domain. In inviscid transonic flow around an airfoil, the presence of a sonic line and shock waves leads to adjacent regions governed by different PDEs~\cite{kholodar2003parametric}. Heat conduction in composite materials is governed by heat equations with distinct coefficients across regions of differing material properties~\cite{bonnetier2000elliptic}. Similarly, material deformation may be governed by different PDEs depending on whether a yield criterion is met~\cite{simo2006computational}.

Although existing methods such as PDE-FIND~\cite{rudy2017data} and Robust-IDENT~\cite{he2022robust} can be applied locally, as demonstrated in~\cite{he2024much}, they are unable to simultaneously identify multiple governing PDEs and their corresponding phase boundaries. Furthermore, for many practical applications, it is essential not only to detect these boundaries but also to accompany their identification with rigorous uncertainty quantification—an aspect that remains largely unaddressed by existing approaches.

In this paper, we propose a novel method, called \textbf{Phase-IDENT}, for identifying PDEs from noisy data that may be governed by distinct dynamics in disjoint regions separated by unknown phase boundaries. Figure~\ref{kdvbur planform}(a) shows an example trajectory of this type. Our objective is to simultaneously locate phase boundaries with uncertainty quantification and recover the governing PDE expressions within each connected region. When a phase boundary exists, a naive application of PDE identification to the entire dataset can produce models whose numerical solutions are incompatible with the observations. As shown in Figure~\ref{kdvbur planform}(b), a globally identified PDE may fit the earlier portion of the data but completely fail to match the latter part.

In this work, we show that PDEs locally identified using data near a phase boundary consistently yield larger fitting errors than those identified from data away from the boundary. We leverage this observation to efficiently construct an open cover of the underlying phase boundary, and subsequently identify PDEs separately within each connected component of its complement. To precisely locate the boundary, we numerically evolve the PDEs identified outside the open cover while monitoring simulation errors within it. Under general conditions, we demonstrate that these simulation errors exhibit significant changes when crossing the true phase boundary. We employ a statistical change point detection method~\cite{robbins2011mean} to estimate boundary points along with uncertainty quantification. Figure~\ref{kdvbur planform}(c) illustrates the detected boundary together with an associated $80\%$ confidence region.

We conduct numerical experiments to validate and analyze our method using various combinations of PDEs and noise levels. The results demonstrate that Phase-IDENT can effectively recover the underlying PDEs and accurately locate phase boundaries when they exist.

To summarize, our main contributions are:
\begin{enumerate}
\item A novel framework, \textbf{Phase-IDENT}, for the identification of two-phase PDEs and their phase boundaries from noisy data. 
\item A robust phase boundary localization method with uncertainty quantification, accompanied by the theoretical analysis of the order of  numerical evolution error across the phase boundary. 
\item Comprehensive numerical experiments demonstrating Phase-IDENT's effectiveness across various dynamical regimes.
\end{enumerate}

This paper is organized as follows. In Section~\ref{sec2}, we introduce the main steps of Phase-IDENT. In Section~\ref{sec4}, we discuss implementation details of  Phase-IDENT. In  Section~\ref{sec5}, we present  numerical results for method validation and analysis. In Section~\ref{sec6} , we conclude the paper with some discussion.

\section{Proposed Method: Phase-IDENT}\label{sec2}

\begin{figure}
    \centering
    \includegraphics[width=1\linewidth]{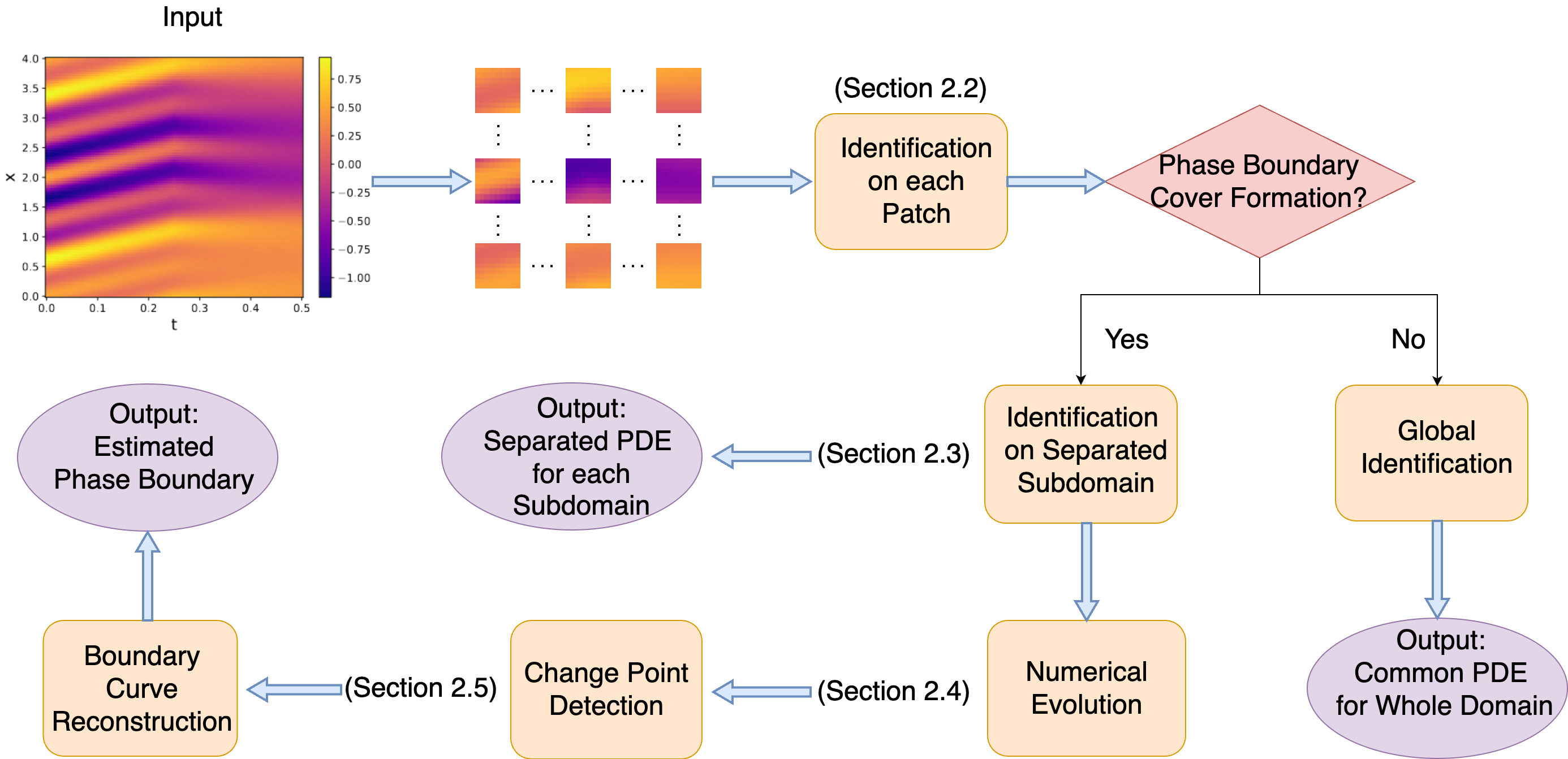}
    \caption{Schematic diagram of the proposed Phase-IDENT workflow. Phase-IDENT begins by covering the time–space domain $\Omega$ with patches. For each patch, a PDE is identified and its fit is assessed via the cross‑validation estimation error (CEE) defined in~\eqref{CEE formu}. If the distribution of high‑CEE patches is scattered (Section~\ref{patch detection}), only a single phase is inferred, and a unique PDE is identified using the entire dataset. Otherwise, the high‑CEE patches are used to construct an open set $\widetilde{\mathcal{C}}$. PDEs are then identified separately within each connected component of $\Omega \setminus \widetilde{\mathcal{C}}$ (Section~\ref{sec_PDE_ident}), and their numerical evolution errors (Section~\ref{evolution section}) inside $\widetilde{\mathcal{C}}$ are employed to locate the phase boundary via change point detection (Section~\ref{boundary estimation}).}
    \label{workflow}
\end{figure}

In this section, we propose a novel method  \textbf{Phase-IDENT} for simultaneously identifying two-phase PDEs and the phase boundaries. We describe the mathematical formulation of the problem and  discuss main ideas here.  The algorithmic details of Phase-IDENT are delegated to Section~\ref{sec4}.

\subsection{Overview of Phase-IDENT}\label{problemsetup}
We focus on demonstrating our method in the following setting. Let  
\[
\Omega = \mathcal{D} \times [0, T] \subset \mathbb{R}^2,
\]
where $\mathcal{D} = [x_{\min}, x_{\max}]$ with $x_{\min} < x_{\max}$, and $T > 0$. This defines a closed rectangular region in the spatial coordinate $x$ and the time coordinate $t$.  
Suppose we have a collection of observations  
\[
\mathcal{U} := \{(x_n, t_n, U_n) \mid n = 1, \dots, N\},
\]
where each tuple $(x_n, t_n, U_n)$ represents a noisy measurement of a physical quantity $u : \Omega \to \mathbb{R}$ at a spatial point $x_n \in \mathcal{D}$ and time $t_n \in [0, T]$.

We assume that $u$ satisfies a \textit{two-phase PDE}. More specifically, suppose $\Omega_1$ and $\Omega_2$ are two disjoint, open, connected subsets of $\Omega$, called \textit{phase domains}, which satisfy
\[
\mathcal{D} \times \{0\} \subset \overline{\Omega}_1,\quad \mathcal{D} \times \{T\} \subset \overline{\Omega}_2,\quad \text{and} \quad \overline{\Omega}_1 \cup \overline{\Omega}_2 = \Omega.
\]
Here $\overline{\Omega}_1$ and $\overline{\Omega}_2$ denote the closures of $\Omega_1$ and $\Omega_2$, respectively. We denote the \textit{phase boundary} separating $\Omega_1$ and $\Omega_2$ by
\[
\Gamma := \partial\Omega_1 \cap \partial\Omega_2.
\]
Moreover, for some $\varepsilon > 0$, we assume that $\Gamma$ can be parameterized by a continuous function $\gamma : \mathcal{D} \to (\varepsilon, T - \varepsilon)$ as follows:
\begin{equation}\label{eq_true_phase_boundary}
\begin{aligned} 
\Gamma : \mathcal{D} & \to \Omega, \\
x & \mapsto (x, \gamma(x)).
\end{aligned}  
\end{equation}
Thus, the image of $\Gamma$ is a simple curve whose endpoints intersect the boundary intervals $\{x_{\min}\} \times (\varepsilon, T - \varepsilon)$ and $\{x_{\max}\} \times (\varepsilon, T - \varepsilon)$. 
In other words, the observation domain $\Omega$ is sufficiently large to capture the phase transition for each spatial point.

Let $u : \Omega \to \mathbb{R}$ be a continuous function such that the restrictions $u|_{\Omega_1}$ and $u|_{\Omega_2}$ are at least $P$-times differentiable with respect to the spatial variable for some integer $P \geq 0$. Assume that $u$ satisfies the following PDE system:
\begin{equation}\label{eq_main_equation}
\begin{cases}
\displaystyle \partial_t u = \sum_{k=1}^K a_k f_k, & \text{in } \Omega_1, \\[10pt]
\displaystyle \partial_t u = \sum_{k=1}^K b_k f_k, & \text{in } \Omega_2,
\end{cases}
\end{equation}
with appropriate boundary conditions\footnote{We focus on identifying PDEs in the interior of $\Omega \setminus \Gamma$. Hence, the boundary conditions can be arbitrary as long as~\eqref{eq_main_equation} is well-posed.}. For $k = 1, \dots, K$, each $f_k$ is called a \textit{feature operator} that maps $u$ to monomials of partial derivatives of $u$ up to order $P$; for example, $u_x$ or $u u_{xx}$ if $P \geq 2$. The scalars $a_k, b_k \in \mathbb{R}$ are called the \textit{feature coefficients} associated with the $k$-th feature in domains $\Omega_1$ and $\Omega_2$, respectively. Features with nonzero coefficients are called \textit{active}, while those with zero coefficients are called \textit{inactive}. For better interpretability~\cite{brunton2016discovering,kang2021ident}, we assume that only a few features are active in both $\Omega_1$ and $\Omega_2$.
 
Given a collection of data $\mathcal{U} \subset \mathcal{D} \times [0,T] \times \mathbb{R}$, we propose \textbf{Phase-IDENT} to achieve the following objectives:

\begin{itemize}
\item \textbf{Locate the phase boundary.} We aim to find a continuous curve
\begin{equation}\label{widehatgamma}
\begin{aligned} 
\widehat{\Gamma}: \mathcal{D} & \to \Omega, \\
x & \mapsto \bigl(x, \widehat{\gamma}(x)\bigr),
\end{aligned}  
\end{equation}
that approximates the underlying phase boundary~\eqref{eq_true_phase_boundary}. This leads to two regions $\widehat{\Omega}_1$ and $\widehat{\Omega}_2$, which approximate the phase domains $\Omega_1$ and $\Omega_2$, respectively. In addition, for any confidence level $p \in [0,1]$ and any $x \in \mathcal{D}$, we find a value $\varepsilon(x,p) > 0$ satisfying 
\begin{equation}\label{eq_boundary_confidence_level}
\PP\Bigl( \bigl|\widehat{\gamma}(x) - \gamma(x)\bigr| \leq \varepsilon(x,p)\Bigr) \geq p ,
\end{equation}
for some probability measure $\PP$ derived from a change point detection model. This pointwise confidence bound provides informative uncertainty quantification for the estimated phase boundary.

\item \textbf{Identify the PDE for each phase domain.} We aim to determine the active features in both equations of~\eqref{eq_main_equation} and estimate the associated coefficients. Denoting the estimated coefficients by $(\widehat{a}_k)_{k=1}^K$ and $(\widehat{b}_k)_{k=1}^K$, the identified model can be expressed as
\begin{equation}\label{eq_estimated_equation}
\begin{cases}
\displaystyle \partial_t u = \sum_{k=1}^K \widehat{a}_k f_k, & \text{in } \widehat{\Omega}_1, \\[10pt]
\displaystyle \partial_t u = \sum_{k=1}^K \widehat{b}_k f_k, & \text{in } \widehat{\Omega}_2.
\end{cases}
\end{equation}
\end{itemize}

Figure~\ref{workflow} presents an overview of the proposed Phase-IDENT workflow. The method recovers the two-phase PDEs and the phase boundary through the following steps:

\begin{itemize}
\item \textbf{Initial estimation of the phase boundary} (Section~\ref{patch detection}). We begin by constructing an open covering of $\Omega$ using small patches. The union of patches where locally identified PDEs fail to fit the data serves as an initial open cover $\widetilde{\mathcal{C}} \subset \Omega$ for the underlying phase boundary. 

\item \textbf{PDE identification within connected components} (Section~\ref{sec_PDE_ident}). Inside each connected component of $\Omega \setminus \widetilde{\mathcal{C}}$, we identify a domain-specific PDE using Robust IDENT~\cite{he2022robust} refined by the trimming technique from~\cite{tang2023weakident}.

\item \textbf{Local numerical evolution} (Section~\ref{evolution section}). Within $\widetilde{\mathcal{C}}$, we numerically evolve the identified PDEs from both phase domains. The compatibility between the data and the PDEs is quantified via an error sequence, whose magnitude reflects the quality of fit.

\item \textbf{Phase boundary localization} (Section~\ref{boundary estimation}). We detect significant changes in the short-time numerical evolution errors of the identified PDEs. The estimated phase boundary $\widehat{\Gamma}$ is obtained by fitting these detected change points with a parametric curve. Uncertainty in change point detection naturally yields confidence bounds for the phase boundary location.
\end{itemize}

Further details on each of these major steps are provided below.

\begin{remark} We note that our proposed framework can be generalized to identify PDEs with \textit{more than two phases}. In higher-dimensional spaces, the extension follows a similar approach; however, special care is required when handling intersections of phase boundaries. We leave these extensions to future work.
\end{remark}

\subsection{Phase detection and initial cover construction}\label{patch detection}
\begin{figure}
    \centering
    \begin{tabular}{c@{\vspace{2pt}}c@{\vspace{2pt}}c}
    \includegraphics[width=0.32\linewidth]{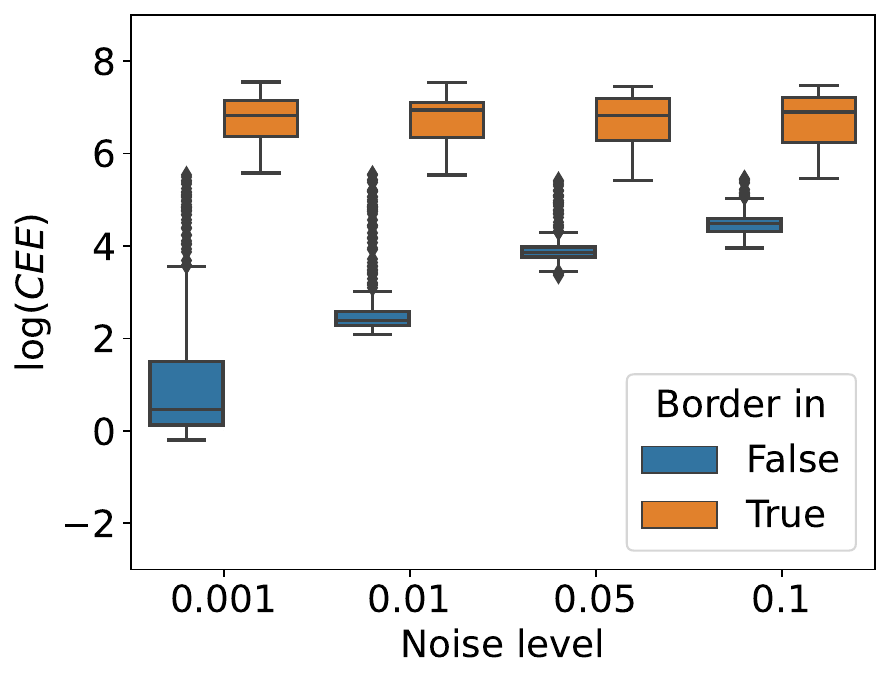}&
    \includegraphics[width=0.32\linewidth]{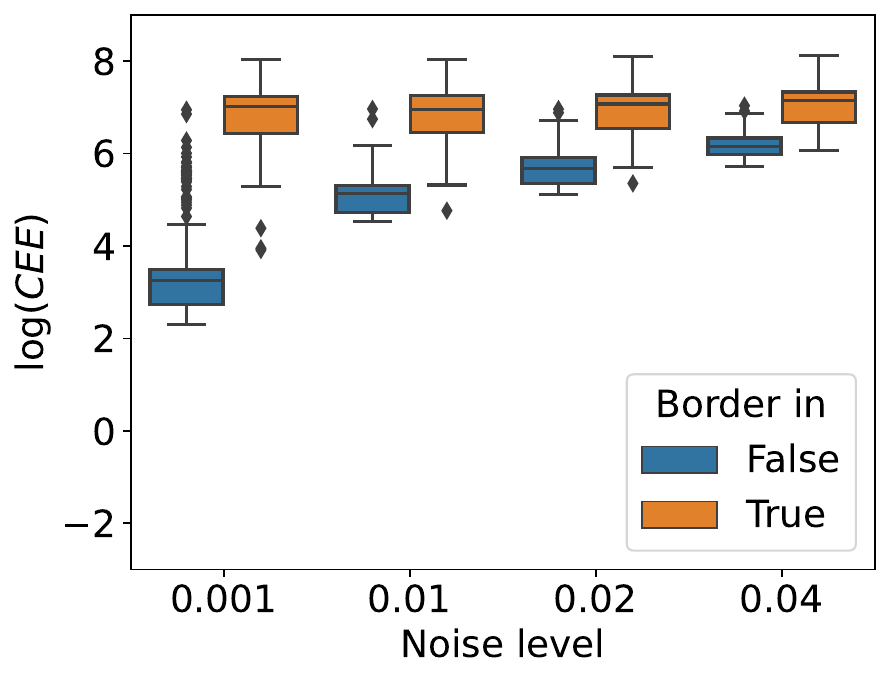}&
    \includegraphics[width=0.32\linewidth]{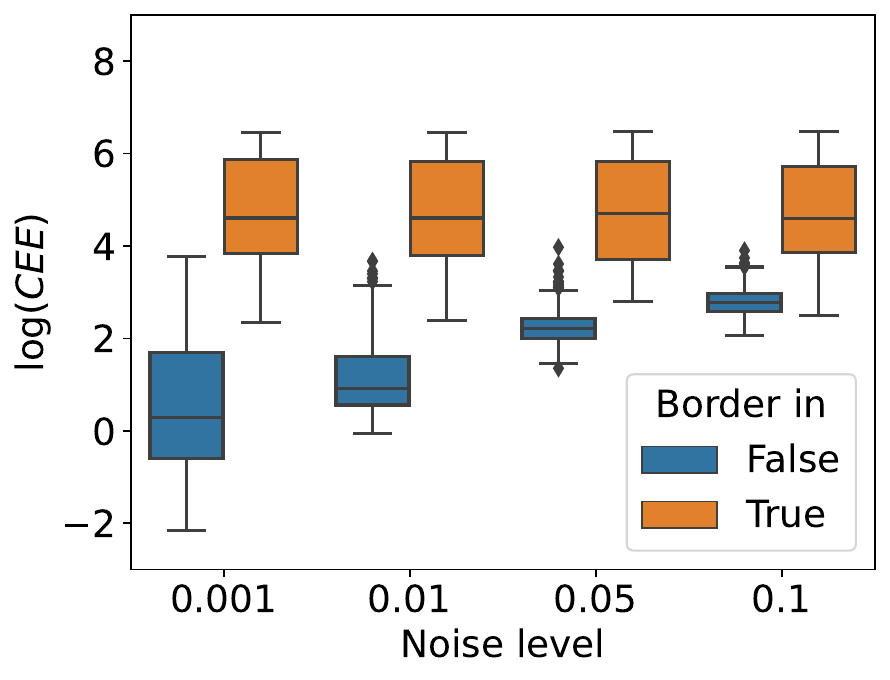}\\
    (a)&(b)&(c)
    \end{tabular}
    \caption{Comparison of the distributions of CEE values for PDEs identified from patches that contain (orange boxes) and do not contain (blue boxes) phase boundaries. The underlying two-phase PDEs are: (a) Transport $\to$ Viscous Burgers (T$\to$VB); (b) KdV $\to$ Burgers (KdV$\to$B); and (c) Burgers $\to$ Transport (B$\to$T). Across a range of noise-to-signal ratios (NSRs) in the observed data, patches with high CEE values serve as consistent indicators of phase boundary presence.
    }
    \label{CEE distri}
\end{figure}

A significant challenge in identifying two-phase PDEs is that the available data cannot be used directly for a global identification, because no single PDE accurately describes the dynamics across both phase domains (see Figure~\ref{kdvbur planform}(b)). Moreover, we lack prior knowledge of how to partition the domain into regions governed by distinct PDEs, or even whether two separate phases truly exist.

Our approach addresses this challenge by first testing for the presence of phase separation. If no phase separation is detected, we identify a single PDE using the entire dataset. If phase separation is confirmed, we construct a cover $\widetilde{\mathcal{C}}$ for the underlying phase boundary. By excluding $\widetilde{\mathcal{C}}$, we can then apply standard single-phase PDE identification algorithms separately to each remaining region.
    
We begin with a collection of closed subsets $\Pi := \{\mathcal{P}_j \subseteq \Omega : j = 1,\dots,J\}$ of the time–space domain $\Omega$, where each $\mathcal{P}_j$ is simply connected and $\bigcup_{j=1}^J \mathcal{P}_j = \Omega$. Suppose $\{f_1,\dots,f_K\}$ is the set of candidate features that may appear in either phase domain. Within each patch $\mathcal{P}_j$, we identify a PDE using Robust-IDENT~\cite{he2022robust} refined with trimming~\cite{tang2023weakident}; see Section~\ref{sec_PDE_ident} for details. Denote the PDE identified from the $j$-th patch as
\begin{equation}
\partial_t u = \sum_{k=1}^K \widehat{c}_{j,k} f_k,
\end{equation}
where $\widehat{c}_{j,k} \in \mathbb{R}$ is the estimated coefficient for the $k$-th feature. We then assess the fit of this PDE to the observed data by evaluating the cross-validation estimation error (CEE)~\cite{he2022robust}, defined by
\begin{equation}\label{CEE formu}
\begin{aligned}
\operatorname{CEE}(\mathcal{P}_j) &:= \sum_{(x,t) \in \mathcal{I}_2^j} \Bigl(\partial_t u(x,t) - \sum_{k \in \widehat{\mathcal{S}}_j} \widetilde{c}_k(\mathcal{I}_1^j) f_k(x,t)\Bigr)^2 \\
&\qquad + \sum_{(x,t) \in \mathcal{I}_1^j} \Bigl(\partial_t u(x,t) - \sum_{k \in \widehat{\mathcal{S}}_j} \widetilde{c}_k(\mathcal{I}_2^j) f_k(x,t)\Bigr)^2.
\end{aligned}
\end{equation}
Here $\mathcal{I}_1^j$ and $\mathcal{I}_2^j$ are disjoint subsets of sample points in $\mathcal{P}_j$; $\widetilde{c}_k(\mathcal{I}_1^j)$ and $\widetilde{c}_k(\mathcal{I}_2^j)$ are sparse feature coefficients estimated by least-squares regression using data indexed by $\mathcal{I}_1^j$ and $\mathcal{I}_2^j$, respectively; and $\widehat{\mathcal{S}}_j \subseteq \{1,2,\dots,K\}$ denotes the index set corresponding to the support of $(\widehat{c}_{j,1},\dots,\widehat{c}_{j,K})$.

By Theorem~3.1 of~\cite{he2022robust}, $\operatorname{CEE}(\mathcal{P}_j)$ estimates the prediction error of the identified PDE. If the dynamics exhibited by data in $\mathcal{I}_1^j$ are inconsistent with those in $\mathcal{I}_2^j$, or if one of these two datasets represents a mixture of two dynamics, the quantity $\operatorname{CEE}(\mathcal{P}_j)$ can become significantly large.  

 \begin{figure}
    \centering
    \begin{tabular}{c@{\vspace{2pt}}c@{\vspace{2pt}}c}
    \includegraphics[width=0.34\linewidth]{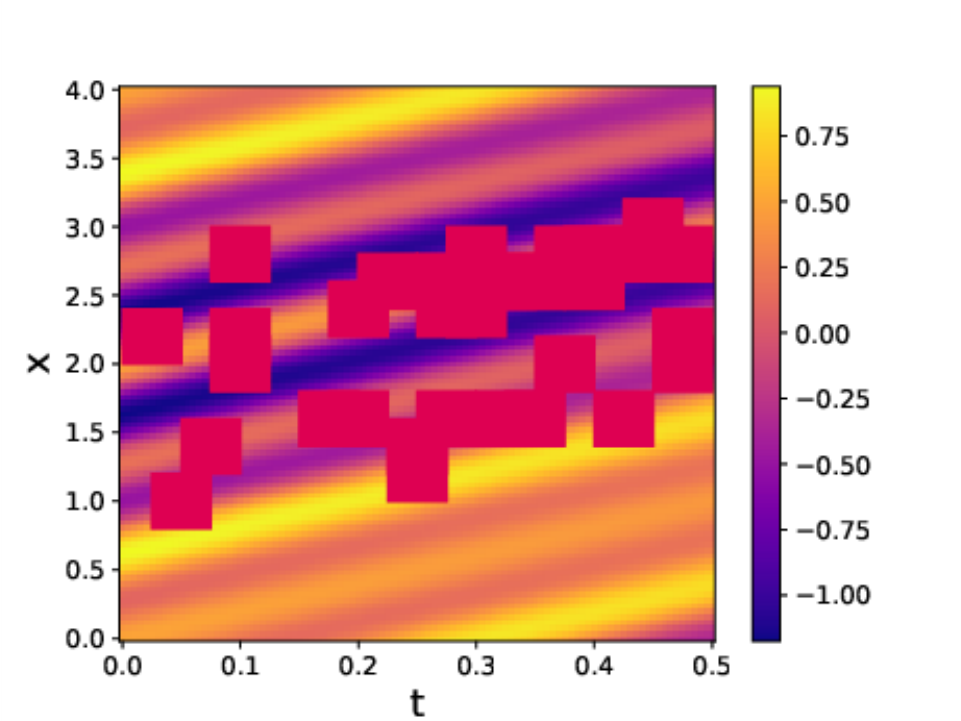}&
    \includegraphics[width=0.34\linewidth]{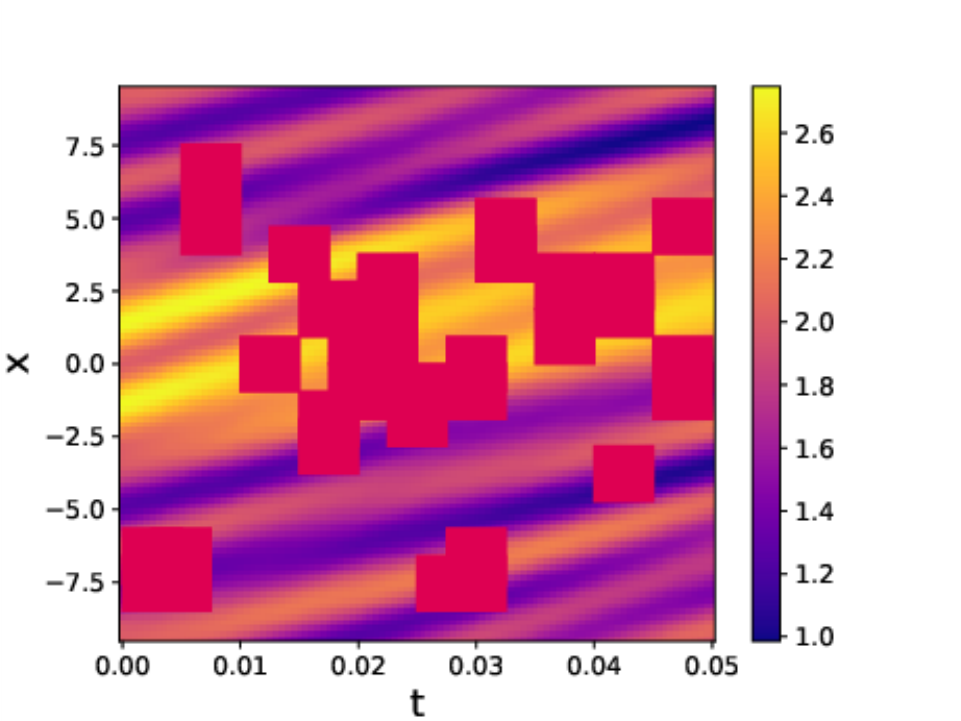}&
    \includegraphics[width=0.34\linewidth]{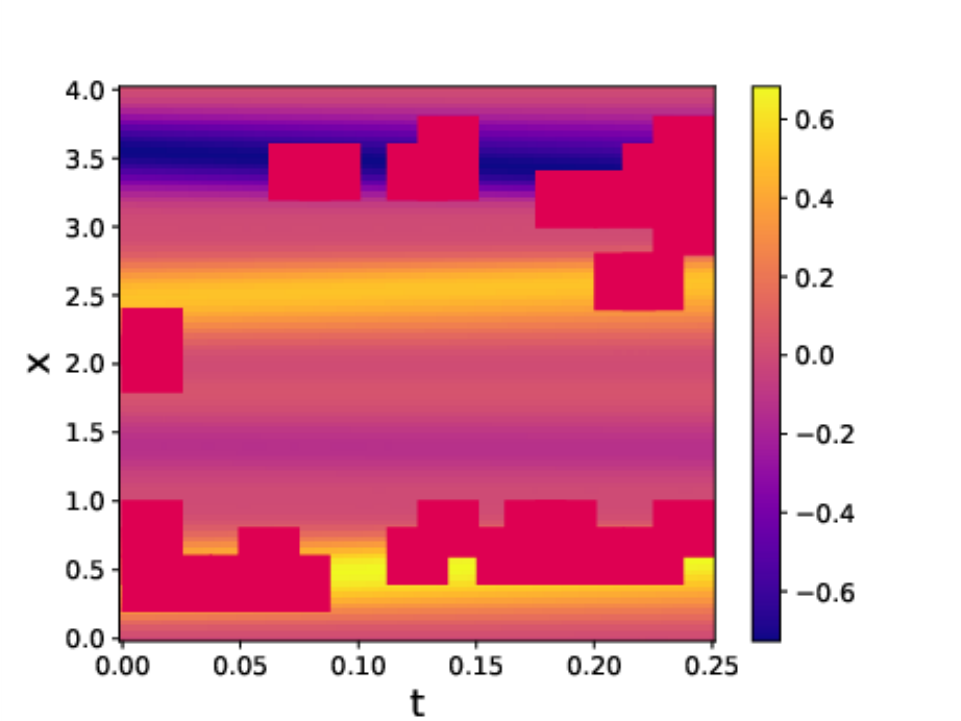}\\
    (a)&(b)&(c)\\
    \includegraphics[width=0.34\linewidth]{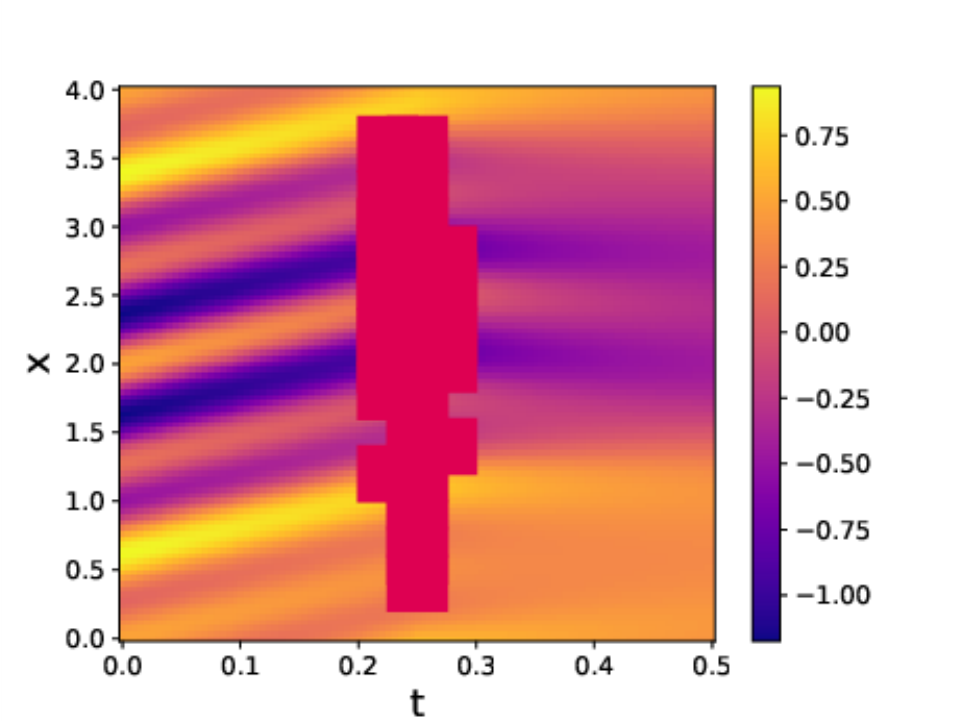}&
    \includegraphics[width=0.34\linewidth]{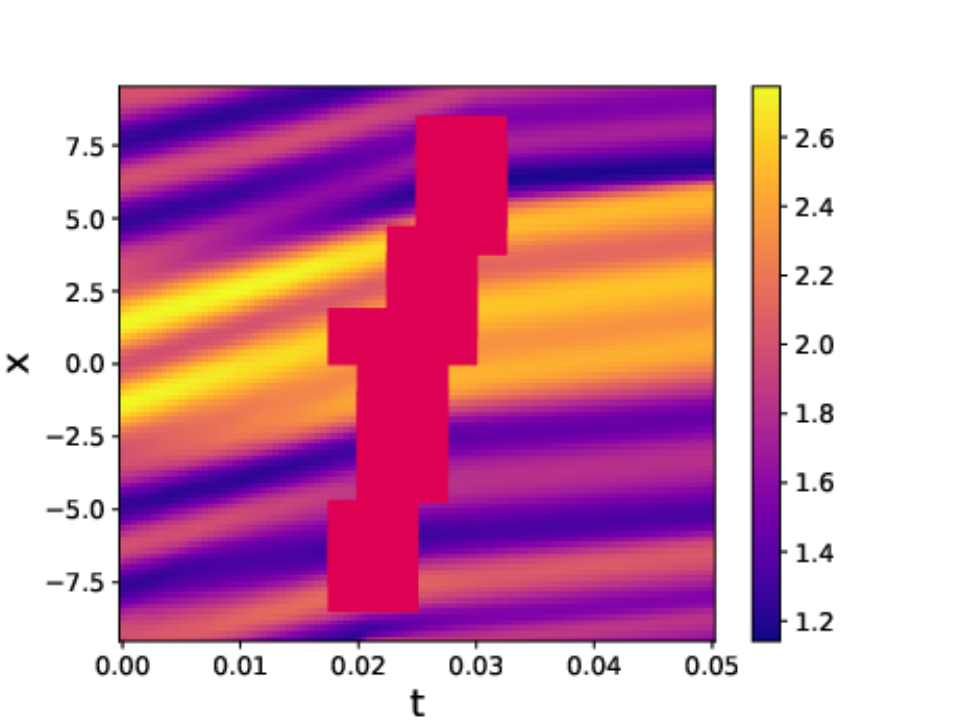}&
    \includegraphics[width=0.34\linewidth]{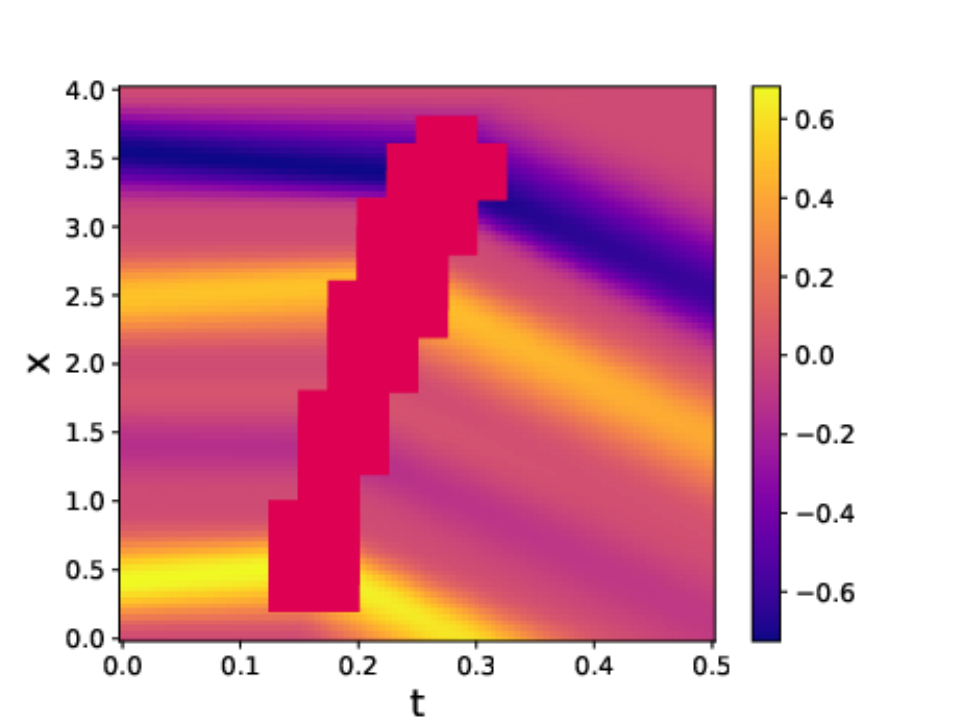}\\
    (d)&(e)&(f)
    \end{tabular}
    \caption{Comparison of the distributions of patches with high CEE values~\eqref{CEE formu} when locally identifying PDEs from data governed by single equations: (a) Transport; (b) KdV; and (c) Burgers; and from data governed by two-phase equations: (d) Transport $\to$ Viscous Burgers (T$\to$VB); (e) KdV $\to$ Burgers (KdV$\to$B); and (f) Burgers $\to$ Transport (B$\to$T). Patches with high CEE are scattered in (a)–(c), whereas they concentrate around the respective phase boundaries in (d)–(f). Using the proposed indicator $r_{\beta}$~\eqref{ra}, the values for (a), (b), and (c) are 0.28, 0.65, and 1.00, respectively; while for (d), (e), and (f), they are 38.02, 16.50, and 8.92, respectively. Throughout this experiment, $\beta = 10$.
    }
    \label{maskpure}
\end{figure}
 
In Figure~\ref{CEE distri}, we compare the CEE values for PDEs identified from patches that contain versus do not contain the phase boundary across three scenarios: Transport $\to$ Viscous Burgers (T$\to$VB), KdV $\to$ Burgers (KdV$\to$B), and Burgers $\to$ Transport (B$\to$T), where the equation before and after ``$\to$'' governs the data in $\Omega_1$ and $\Omega_2$, respectively. We observe that the CEE values from patches containing the phase boundaries are generally higher than those without. This discrepancy remains significant even when the noise level increases.  Hence, we conclude that patch-based CEE is an effective indicator for the presence of phase boundaries.  Furthermore, Figure~\ref{maskpure} shows that if there is no phase boundary, patches with relatively high CEE are scattered, whereas in two-phase PDEs, high-CEE patches consistently aggregate near the phase boundaries.

These observations imply that both the positions and spatial distribution of high-CEE patches are informative for detecting underlying phase boundaries. We also note that inaccuracies in approximating higher-order spatial derivatives near the spatial boundary can produce large CEE values. Hence, we restrict attention to patches that do not intersect the spatial boundary:
\begin{equation}
\Pi^* = \{\mathcal{P} \in \Pi \mid \mathcal{P} \cap \partial\mathcal{D} \times [0,T] = \varnothing\}.
\end{equation}
For any $\beta \in [0,100]$, we define the set of top-$\beta\%$ CEE patches as
\begin{equation}\label{apsi}
    \Pi^*_{\beta} := \{\mathcal{P} \in \Pi^* \mid \operatorname{CEE}(\mathcal{P}) \geq Q^*(\beta)\},
\end{equation}
where $Q^*(\beta)$ is the $(100-\beta)$-th percentile of the CEE values of patches in $\Pi^*$. To quantify the spatial concentration of high-CEE patches, we introduce
\begin{equation}\label{ra}
    r_{\beta} := \frac{\operatorname{Var}(\{x_{\mathcal{P}} \mid \mathcal{P} \in \Pi^*_{\beta}\}) \Delta t^2}
                      {\operatorname{Var}(\{t_{\mathcal{P}} \mid \mathcal{P} \in \Pi^*_{\beta}\}) \Delta x^2},
\end{equation}
where $(x_{\mathcal{P}}, t_{\mathcal{P}})$ denotes the center of mass of patch $\mathcal{P}$, $\operatorname{Var}(\cdot)$ denotes the sample variance, and $\Delta x$ and $\Delta t$ denote the spatial and temporal grid spacings of the data, respectively.

A large value of $r_{\beta}$ suggests that patches in $\Pi^*_{\beta}$ are aligned in a way consistent with a phase boundary, indicating a likely phase transition in the data. Conversely, a small $r_{\beta}$ indicates that high-CEE patches are scattered, and we conclude that no phase transition is present. This property allows us to decide whether the data correspond to a single phase or multiple phases. If the data are considered single-phase, a global PDE identification can be applied; otherwise, we proceed to the next stage.

The choice of $\beta$ does not strongly influence $r_\beta$ unless $\beta$ is close to zero; in this work we set $\beta = 10$. Figure~\ref{maskpure} demonstrates the effectiveness of the proposed indicator~\eqref{ra}: for PDEs without phase transitions, the values of $r_{\beta}$ are consistently smaller than those for two-phase PDEs.
 
\subsection{PDE identification in phase domains}\label{sec_PDE_ident}

When the indicator~\eqref{ra} yields a large value, the majority of high-CEE patches concentrate around the underlying phase boundary. The complement of their union then contains patches whose data are governed by a single PDE. In this study, we focus on PDE identification with two phases. To obtain a rough estimate of the phase boundary location, we consider the largest connected component $\mathcal{C} \subseteq \Omega$ of the union of high-CEE patches from
\begin{equation}\label{psi}
\Psi := \{\mathcal{P} \in \Pi \mid \mathcal{P} \cap \mathcal{D} \times \partial[0,T] = \varnothing\}.
\end{equation}
Patches intersecting the \textit{temporal} boundary are excluded so that $\Omega \setminus \mathcal{C}$ consists of two connected components, enabling separate PDE identification in each. The construction details of $\mathcal{C}$ are discussed in Section~\ref{initial cover}.

To facilitate the numerical evolution in Section~\ref{evolution section} for further phase boundary localization, we enlarge $\mathcal{C}$ to an open set $\widetilde{\mathcal{C}} \subseteq \Omega$ such that the upper and lower bounds of time $t$ within $\widetilde{\mathcal{C}}$ can be given by two functions of $x$. In particular, $\widetilde{\mathcal{C}}$ contains $\mathcal{C}$ and can be expressed as
\begin{equation}\label{eq_Ca_hat}
\widetilde{\mathcal{C}} := \{(x,t) \in \Omega \mid \gamma^{\ell}(x) < t < \gamma^{r}(x)\}
\end{equation}
for some functions $\gamma^\ell,\gamma^r: \mathcal{D} \to (0,T)$ satisfying $\gamma^\ell(x) < \gamma^r(x)$ for every $x \in \mathcal{D}$. The detailed construction of $\gamma^\ell$ and $\gamma^r$ is given in Section~\ref{Sec_convexification}. 

Denote by $R^{\ell}$ and $R^{r}$ the two connected components of $\Omega \setminus \widetilde{\mathcal{C}}$, with
\[
R^{\ell} = \{(x,t) \in \Omega \mid t \leq \gamma^\ell(x)\}, \qquad 
R^{r} = \{(x,t) \in \Omega \mid t \geq \gamma^r(x)\}.
\]
Within each connected component, we identify a single PDE governing the regional dynamics using a modified version of Robust-IDENT~\cite{he2022robust}.

Specifically, for $\delta \in \{\ell, r\}$, let $\{(x^\delta_1,t^\delta_1),\dots, (x^\delta_{N_\delta},t^\delta_{N_\delta})\} \subset R^\delta$ be the collection of sample points, where $N_\delta \geq 1$ is the number of regional data points. For a dictionary of candidate features $\{f_1,\dots,f_K\}$ as introduced in~\eqref{eq_main_equation}, we define the \textbf{regional feature matrix} and the \textbf{feature response} as
\begin{equation}\label{eq_feature_system}
\bF_\delta := \begin{bmatrix} \mathbf{f}_{\delta,1} & \mathbf{f}_{\delta,2} & \cdots & \mathbf{f}_{\delta,K} \end{bmatrix} \in \mathbb{R}^{N_\delta \times K}, \qquad 
\by_\delta := \begin{bmatrix} \partial_t u(x_{1}^\delta, t_{1}^\delta) & \cdots & \partial_t u(x_{N_\delta}^\delta, t_{N_\delta}^\delta) \end{bmatrix}^\top \in \mathbb{R}^{N_{\delta}},
\end{equation}
where $\mathbf{f}_{\delta,k} := \begin{bmatrix} f_k(x_{1}^\delta, t_{1}^\delta) & \cdots & f_k(x_{N_\delta}^\delta, t_{N_\delta}^\delta) \end{bmatrix}^\top$ for $k = 1,\dots,K$.

We then generate candidate models by solving the $\ell_0$-constrained problem
\begin{equation}\label{eq_k_sparse}
\begin{aligned}
\widehat{\bc}^\delta_k \in \arg\min_{\bc\in\mathbb{R}^K} &\ \|\bF_\delta \bc - \by_\delta\|_2^2 \\
\text{s.t.} &\ \|\bc\|_0 = k,
\end{aligned}
\end{equation}
for $k = 1,\dots,K$. Problem~\eqref{eq_k_sparse} is NP-hard~\cite{nguyen2019np}. In~\cite{he2022robust}, the greedy subspace pursuit (SP) algorithm~\cite{dai2009subspace} is shown to be effective in recovering the correct features. For completeness, we provide pseudo-code for SP in Appendix~\ref{alg_SP}. 

Once the candidate feature coefficients $\{\widehat{\bc}^\delta_1,\dots,\widehat{\bc}^\delta_{K}\}$ for the dynamics in $R^\delta$ are obtained, we select the optimal model as the one with the smallest CEE~\eqref{CEE formu}. Note that only observations of the solution $u$ are available; therefore, we approximate the entries of $\bF_\delta$ using a 5-point ENO scheme~\cite{harten1997uniformly} and those of $\by_\delta$ via forward differences. When the data contain noise, we apply successive denoised differentiation (SDD)~\cite{he2022robust}. Additionally, we adopt the trimming technique introduced in~\cite{tang2023weakident}, which effectively improves feature selection by eliminating unimportant terms. After generating the candidates, we therefore apply trimming to obtain refined models. Further details can be found in~\cite{tang2023weakident}.

To successfully identify the PDEs using regional data, it is natural to require that both $R^\ell$ and $R^r$ contain sufficiently many data. By examining the null space, an obvious necessary condition for finding the underlying PDE via~\eqref{eq_k_sparse} is $N_\delta>2k^*$ where $k^*$ is the true number of features. In~\cite{he2022asymptotic}, the  exact recovery of the PDE features as $N_\delta\to+\infty$ is guaranteed for identification frameworks using $\ell_1$-regularization; and it can be generalized to~\eqref{eq_k_sparse} under structural assumptions about the regional feature matrices~\cite{candes2006stable}. In~\cite{he2024much}, it was  shown that the regional data needs to contain enough numbers of Fourier modes for a unique identification, which was also found to be crucial in accurate coefficient reconstruction~\cite{tang2023weakident}.

We also note that a similar local PDE identification framework was employed in the patch-based method CaSLR~\cite{he2024much}, which additionally enforces a global consistency constraint on the supports of all locally estimated coefficients. Since the identification procedure described above is applied independently to each connected component of $\Omega \setminus \widetilde{\mathcal{C}}$, the process can be implemented in parallel to improve computational efficiency.

\subsection{Local evolution error across phase boundary}\label{evolution section}
In this work, we assume the underlying phase boundary $\Gamma$ is a simple curve that can be parameterized as a continuous map from $\mathcal{D}$ to $\Omega$. For each fixed $x \in \mathcal{D}$, define $L_x(t) = (x, t)$ for $t \in [0,T]$. We expect the dynamics along $L_x([0,T])$ to change noticeably near its intersection with $\Gamma$, which can be used to locate the phase boundary more precisely.

Based on this intuition, we propose to numerically evolve the identified PDEs along the lines $L_x([0,T])$ for sampled spatial points $x \in \mathcal{D}$. From Section~\ref{sec_PDE_ident}, we have the identified PDEs
\begin{equation}\label{twoidentfied}
u_t (x,t) = \cF_\delta(x,t) := \sum_{k=1}^K \widehat{a}_k^{\delta} f_k(x,t), \quad \text{for } (x,t) \in R^\delta, \quad \delta \in \{\ell, r\},
\end{equation}
where $\{\widehat{a}_k^{\delta}\}_{k=1}^K$ denotes the estimated coefficients in $R^\delta$.

To localize the phase boundary, for each $x$ we discretize $\cF_\ell$ and $\cF_r$ into $\widehat{\cF}_\ell$ and $\widehat{\cF}_r$, respectively, and compare their single‑step evolution with the data. Given a sampled spatial point $x \in \mathcal{D}$, we consider a uniform sequence of time grid points
\[
\tau_m(x) = \gamma^\ell(x) + (m-1) \Delta t, \quad m = 0,1,\dots, M_x,
\]
where $M_x = (\gamma^r(x) - \gamma^\ell(x)) / \Delta t$. We compute
\begin{equation}\label{numerical evo}
 \widetilde{U}^{\delta}_{m+1}(x) := \widehat{U}_m(x) + \Delta t \cdot \widehat{\cF}_\delta\bigl(\bm{U}_m(x)\bigr), \quad \delta \in \{\ell, r\},
 \end{equation}
where $\bm{U}_m(x) = (U_{m}(x-h\Delta x),\dots, U_{m}(x+h\Delta x)) \in \mathbb{R}^{2h+1}$ for some positive integer $h \geq 1$ denotes denoised data at time $\tau_m(x)$, and $\widehat{U}_m(x)$ approximates $u(x,\tau_m(x))$ by a linear combination of the data in $\bm U_m(x)$. Periodic boundary conditions in space are assumed.

For each $\delta \in \{\ell, r\}$, we define the single‑step evolution error along $L_x([0,T])$ at $\tau_m(x)$ as
 \begin{equation}\label{edelta}
 e^{\delta}(m, x) := \bigl|\widetilde{U}^{\delta}_{m+1}(x) - U_{m+1}(x)\bigr|, \qquad m = 0, 1, \dots, M_x.
 \end{equation}
For each $x \in \mathcal{D}$, we thus obtain two sequences
\begin{equation}\label{eq_sequences}
\be^{\ell}(x) := \bigl( e^{\ell}(0, x), \dots, e^{\ell}(M_x, x) \bigr)
\quad\text{and}\quad
\be^r(x) := \bigl( e^{r}(0, x), \dots, e^{r}(M_x, x) \bigr),
\end{equation}
each recording the single‑step evolution errors of the PDEs identified in $R^\ell$ and $R^r$, respectively, within the interval $t \in (\gamma^\ell(x), \gamma^r(x))$.

We expect large errors when numerically evolving the incorrect model. That is, $e^\ell(m,x)$ stays small if $\tau_m(x) < \Gamma(x)$ and becomes large if $\tau_m(x) > \Gamma(x)$, whereas $e^r(m,x)$ stays large if $\tau_m(x) < \Gamma(x)$ and becomes small if $\tau_m(x) > \Gamma(x)$.

\begin{figure}
    \centering
    \begin{tabular}{c@{\vspace{2pt}}c@{\vspace{2pt}}c}
    \includegraphics[width=0.33\textwidth]{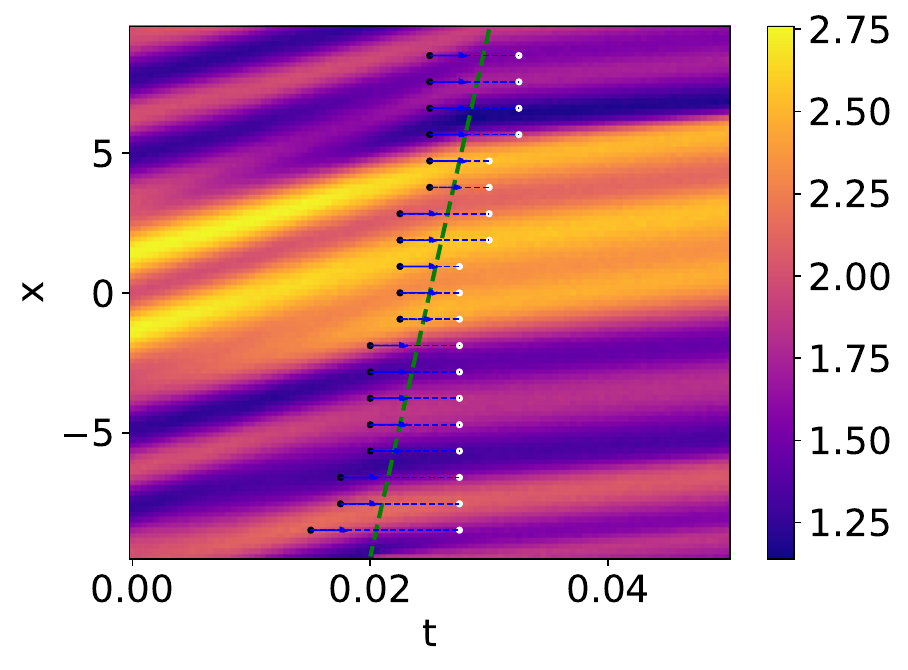}&
    \includegraphics[width=0.33\textwidth]{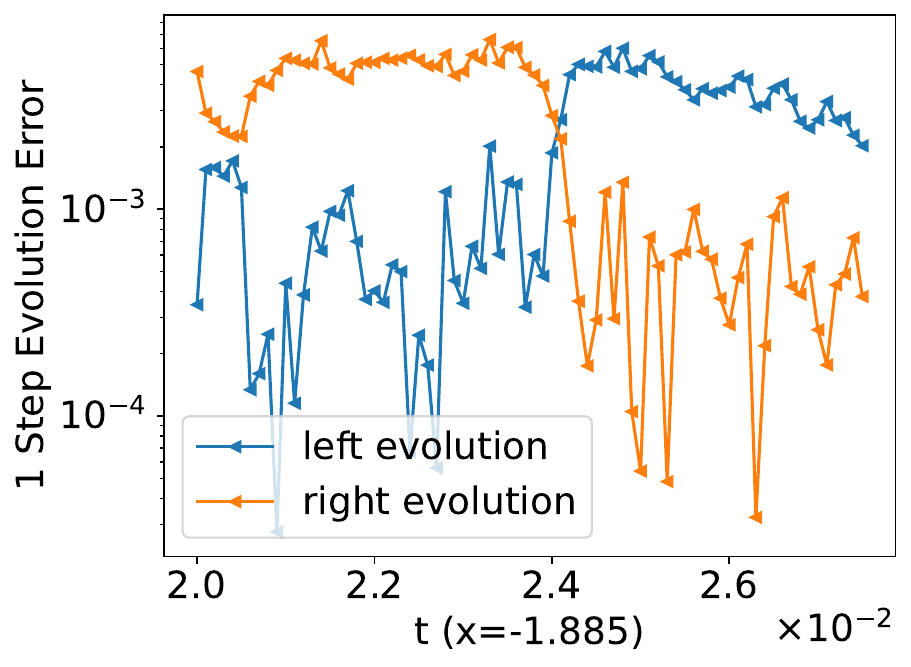}&
    \includegraphics[width=0.33\textwidth]{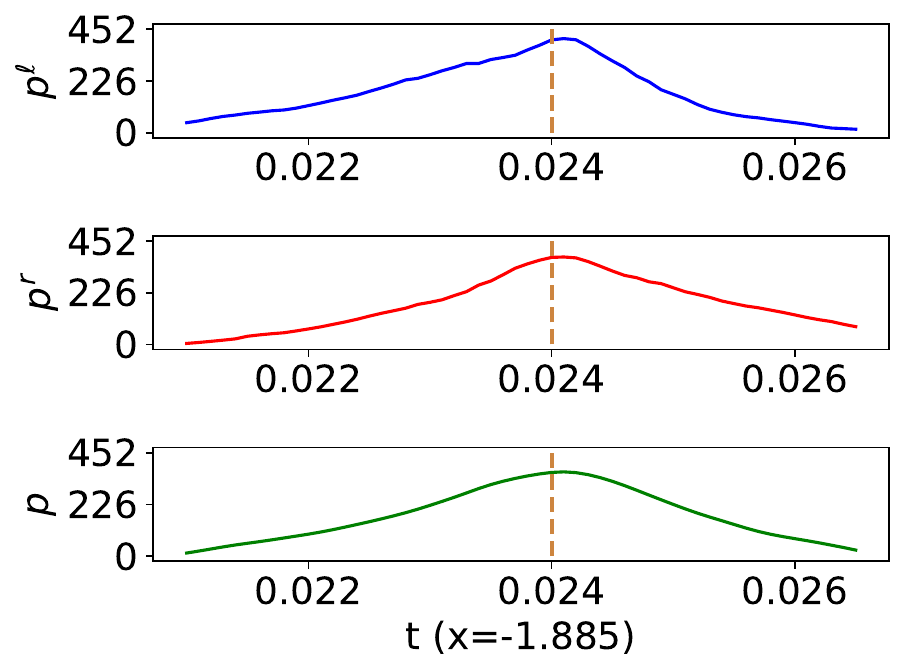}\\
    (a)&(b)&(c)
    \end{tabular}
    \caption{(a) Noisy observations of a trajectory governed first by the KdV equation and then by the Burgers equation. The green dashed line indicates the true phase boundary. The PDEs identified on both sides are numerically evolved along horizontal slices (blue dashed lines). (b) Simulation error sequences $\be^{\ell}(x)$ (blue) and $\be^{r}(x)$ (orange) defined in~\eqref{eq_sequences}. (c) Density functions $p^\ell_{x}$ (blue) and $p^r_{x}$ (red) defined in~\eqref{pdelta}, and their Wasserstein barycenter $p_{x}$ (green) defined in~\eqref{ot}. The vertical brown dashed line in (c) marks the ground truth location of the boundary point. Here $x = -1.885$ and the data contain $2\%$ additive noise.}

    \label{fig_demo_cumsum}
\end{figure}

We analyze the variation of $\mathbf{e}^\ell(x)$ in~\eqref{eq_sequences} near the phase boundary; the analysis for $\mathbf{e}^r(x)$ is analogous. For any $(x,\tau_m) \in \Omega_1$, we assume that the numerical scheme~\eqref{numerical evo} is consistent whenever $\bm U_m(x) \subseteq \Omega_1$, i.e.,
\begin{equation}\label{assumption}
    \begin{aligned}
    u(x,\tau_m) &= \widehat{U}_m(x) + \mathcal{O}(\Delta x^p), \quad p > 1,\\
    \mathcal{F}_\ell(x,\tau_m) &= \widehat{\mathcal{F}}_\ell(\bm{U}_m(x)) + \mathcal{O}(\Delta x^q), \quad q > 1.
    \end{aligned}
\end{equation}
If $(x,\tau_{m+1}) \in \Omega_1$, substituting the Taylor expansion of $u(x,\tau_{m+1})$ about $(x,\tau_m)$ together with~\eqref{assumption} into~\eqref{numerical evo} yields the truncation error
\begin{equation}\label{T1}
\begin{aligned}
    e^\ell(m,x) &= u(x,\tau_{m+1}) - u(x,\tau_{m}) - \Delta t\bigl(u_t(x,\tau_{m}) + \mathcal{O}(\Delta x^q)\bigr) + \mathcal{O}(\Delta x^p) \\
    &= \Delta t\, u_t(x,\tau_{m}) + \mathcal{O}(\Delta t^2) - \Delta t\bigl(u_t(x,\tau_{m}) + \mathcal{O}(\Delta x^q)\bigr) + \mathcal{O}(\Delta x^p) \\
    &= \mathcal{O}(\Delta t^2) + \mathcal{O}(\Delta t\,\Delta x^q) + \mathcal{O}(\Delta x^p).
\end{aligned}
\end{equation}

Next, we examine the $m$-th component of $\mathbf{e}^\ell(x)$ when $\tau_m < \gamma(x) < \tau_{m+1}$. In this case, the Taylor expansion used in~\eqref{T1} cannot be applied, because the smoothness of $u(x,t)$ between $(x,\tau_m)$ and $(x,\tau_{m+1})$ is not guaranteed.

\begin{proposition}\label{proposition}
For $x \in \mathcal{D}$, if $\bm{U}_m(x) \subseteq \Omega_1$ and $\tau_m < \gamma(x) < \tau_{m+1}$, let  
\[
\partial_t^+ u(x,\gamma(x)) := \lim_{h \to 0^+} \frac{u(x,\gamma(x)+h) - u(x,\gamma(x))}{h}
\]
and  
\[
\partial_t^- u(x,\gamma(x)) := \lim_{h \to 0^-} \frac{u(x,\gamma(x)+h) - u(x,\gamma(x))}{h}.
\]
Then
\begin{equation}\label{T2}
e^\ell(m,x) = \bigl(\partial_t^+ u(x,\gamma(x)) - \partial_t^- u(x,\gamma(x))\bigr) \mathcal{O}(\Delta t) + \mathcal{O}(\Delta x^p).
\end{equation}
\begin{proof}
See Appendix~\ref{proof of propo}.
\end{proof}
\end{proposition}

From Proposition~\ref{proposition}, we see that when $\bm U_m(x) \subseteq \Omega_1$, the order of $\Delta t$ in $e^\ell(m,x)$ reduces to $1$ as the trajectory crosses the phase boundary, and the error is proportional to the difference between the left and right derivatives at $(x,\gamma(x))$. This indicates a significant growth in $e^\ell(m,x)$.

When $\bm{U}_m(x)$ contains data from both $\Omega_1$ and $\Omega_2$, the order of $e^\ell(m,x)$ depends on the accuracy of the evolution scheme~\eqref{numerical evo} and on how the data in $\bm{U}_m(x)$ are partitioned between the two phases. If all $\bm{U}_m(x)$ are from $\Omega_2$, the magnitude of entries of $\be^\ell(x)$ around the phase boundary is governed by the difference between $\widehat{\cF}_\ell(\bm{U}_m(x))$ and $\widehat{\cF}_r(\bm{U}_m(x))$. In these last two scenarios, $\widehat{\cF}_\ell$ in~\eqref{numerical evo} is computed using data that are partially or entirely from $\Omega_2$. Consequently, $\be^\ell(x)$ is expected to remain small initially, rise sharply when crossing $\Gamma$, and then fluctuate at a larger magnitude. For such a pattern, the time point that maximizes the difference between the preceding and subsequent average error values is close to the location where the error jumps; that is, the point where the trajectory crosses the phase boundary.

The experimental results support the analysis above. In Figure~\ref{fig_demo_cumsum}(a), we show a set of evolution trajectories governed first by the KdV equation and then by the Burgers equation, with $2\%$ additive noise. The blue horizontal lines indicate $L_x([0,T])\cap \widetilde{\mathcal{C}}$ for a collection of spatial points $x \in \mathcal{D}$. In (b), we display the sequences~\eqref{eq_sequences} for $x = -1.855$, where the blue and orange curves represent $\be^\ell(x)$ and $\be^r(x)$, respectively.  We observe that $\be^{\ell}(x)$ remains relatively low for $t \leq 0.024$, increases sharply near $t = 0.024$, and stays high thereafter. This suggests that data prior to $t \approx 0.024$ are consistent with the PDE identified in $R^\ell$, whereas data after $t \approx 0.024$ deviate from it. In contrast, $\be^r(x)$ stays relatively high for $t \leq 0.024$, drops abruptly around $t = 0.024$, and remains low afterward. This indicates that the data become compatible with the PDE identified for $R^r$ only after $t \approx 0.024$.

\subsection{Change point detection and uncertainty quantification}\label{boundary estimation}
The discussion in Section~\ref{evolution section} suggests that abrupt changes in the error sequences~\eqref{eq_sequences} associated with a spatial point $x \in \mathcal{D}$ indicate where the evolution line $L_x([0,T])$ crosses the underlying phase boundary $\Gamma(x)$. Our goal is to detect such a change point and provide uncertainty quantification for the resulting estimate.

For each $x \in \mathcal{D}$, we compute the cumulative sum (CUSUM) statistics~\cite{robbins2011mean} associated with the error sequence $\be^\delta$:
\begin{equation}
C^\delta_x(m) = m\Bigl(1 - \frac{m}{M_x+1}\Bigr)\bigl(A_m^\delta (x) - B_m^\delta (x)\bigr)^2, \qquad \delta \in \{\ell, r\},
\end{equation}
for $m = m_0, m_0+1, \dots, M_x - m_0$, where
\begin{equation}\label{Abdelta}
A^\delta_m(x) := \frac{1}{m}\sum_{j=0}^{m-1} e^{\delta}(j, x), \quad 
B^\delta_m(x) := \frac{1}{M_x - m+1}\sum_{j=m}^{M_x} e^\delta(j, x), \qquad \delta \in \{\ell, r\},
\end{equation}
with $M_x = (\gamma^r(x) - \gamma^\ell(x)) / \Delta t$. Here $m_0$ is a positive integer satisfying $0 \leq m_0 < M_x/2$; its role is explained in Section~\ref{bspline}. The statistic $C^\delta_x(m)$ compares the average of $\be^\delta$ before and after index $m$, weighted by the sequence length. Under suitable assumptions on the stochastic structure of the observed sequence, the distribution of the maximum value of the CUSUM statistics can be shown to be closely approximated by that of a Brownian bridge~\cite[Theorem 1]{robbins2011mean}. This property is used to test for the presence of a change point. For further details and extensions, we refer readers to~\cite{macneill1974tests,robbins2011mean}.

We take the index with the maximal CUSUM value as the estimated change point:
\begin{equation}\label{eq_change_point_estimate_one_side}
\widehat{\gamma}^\delta(x) = \gamma^\ell + m^\delta_x\Delta t,
\end{equation}
where $m^\delta_x = \argmax_{m \in \{m_0,\dots,M_x - m_0\}} C^\delta_x(m)$. However, we note that $\widehat{\gamma}^\ell(x)$ and $\widehat{\gamma}^r(x)$ generally differ. To synthesize them into a single estimated boundary point, we propose a novel strategy.

For each simulation direction, we normalize the CUSUM statistics and define
\begin{equation}\label{pdelta}
p^\delta_x(m) := \frac{C_x^\delta(m)}{\sum_{k=m_0}^{M_x-m_0} C_x^{\delta}(k)\Delta t}, \qquad \delta \in \{\ell, r\},
\end{equation}
where $\Delta t$ is the time step. Note that $p_x^\delta$ defines a probability mass function (PMF) over the index set $\{m_0,\dots, M_x - m_0\}$.  If we treat the change point as a random variable following the distribution $p_x^\delta$, then  the estimate~\eqref{eq_change_point_estimate_one_side} corresponds to the maximum likelihood estimator. We combine these two distributions by computing their Wasserstein barycenter~\cite{sturm2003probability}:
\begin{equation}\label{ot}
p_x := \argmin_{p \in \mathcal{P}} \bigl\{W_2^2(p, p_x^{\ell}) + W_2^2(p, p_x^r)\bigr\},
\end{equation}
where $W_2$ is the Wasserstein distance between two probability distributions, and $\mathcal{P}$ is the space of all probability measures on $\{m_0,\dots,M_x-m_0\}$. The minimization in~\eqref{ot} can be performed efficiently via iterative Bregman projection~\cite{benamou2015iterative}. Figure~\ref{fig_demo_cumsum}(c) shows the density functions $p^\ell_x$, $p^r_x$, and $p_x$ derived from the error sequences in (b). Consequently, we estimate $\Gamma(x)$ by the maximum likelihood principle applied to $p_x$, setting
\begin{equation}\label{eq_gamma_estimate}
\widehat{\gamma}(x) = \gamma^\ell + \argmax_{m \in \{m_0,\dots,M_x - m_0\}} p_x(m) \cdot \Delta t.
\end{equation}

The construction of $p_x$ also provides a natural framework for uncertainty quantification of $\widehat{\Gamma}(x)$. Specifically, for any $p \in [0,1]$, a confidence interval for $\widehat{\gamma}(x)$ at level $p$ is $[\widehat{\gamma}(x) - \epsilon(x,p),\; \widehat{\gamma}(x) + \epsilon(x,p)]$, where
\begin{equation}\label{confidence level}
\epsilon(x,p) := \Delta t \cdot \min\Bigl\{h \in \mathbb{Z}_+ \mid \mathbb{P}_x\bigl(|\widehat{\gamma}(x)-\gamma(x)| < h\Delta t\bigr) > p\Bigr\},   
\end{equation}
and $\PP_x: \{m_0,\dots,M_x-m_0\}\to[0,1]$ is the discrete probability measure associated with the PMF $p_x$.  Note that the above constructions are computed for a finite sample of points in $\mathcal{D}$. In Section~\ref{bspline}, we extend them to a continuous setting.

\begin{remark}
    The multi-step numerical evolution scheme used in \cite{he2022robust} could also be applied for accumulating evolution errors. However, we find that single-step evolution not only reduces computational cost but also leads to more confident and accurate change point location estimates. A comparison between single-step and multi-step evolution is provided in Appendix~\ref{multistep}.
\end{remark}

\section{Implementation Details}\label{sec4}
This section provides implementation details for the proposed Phase-IDENT method. In Section~\ref{initial cover}, we describe the selection of the percentile threshold for high-CEE patches from~\eqref{psi} and the formation of a region $\mathcal{C}$ as an initial estimate of the underlying phase boundary. In Section~\ref{Sec_convexification}, we discuss the construction of a larger cover $\widetilde{\mathcal{C}}$ (see~\eqref{eq_Ca_hat}) by extending $\mathcal{C}$. In Section~\ref{bspline}, we present a continuous parametric representation of the phase boundary.
\subsection{Construction of the initial cover}\label{initial cover}

In this work, we initially cover $\Omega$ with rectangular patches. For some integer $N_P \geq 1$, we define the following collection of $N_P^2$ patches that cover $\Omega$:
\begin{equation}\label{eq_Pij}
\begin{aligned}
\mathcal{P}_{i,j} = \bigl\{(x, t) \in \Omega \mid {} &x_{\min} + (i-1)w \leq x \leq x_{\min} + (i+1)w,\\
&(j-1)h \leq t \leq (j+1)h\bigr\}, \qquad i,j = 1,\dots,N_P,
\end{aligned}
\end{equation}
where $w = (x_{\max} - x_{\min})/(N_P+1)$ and $h = T/(N_P+1)$ denote half the spatial width and half the temporal width of each patch, respectively. Collectively, we set $\Pi = \{\mathcal{P}_{i,j} \mid i,j = 1,\dots,N_P\}$.

The patch size must be chosen appropriately. On the one hand, if the patches are too large, all of them would intersect the underlying phase boundary, leading to uniformly high CEE values. As a result, high‑CEE patches would become ineffective for detecting phase changes or constructing an initial cover. Therefore, we require the half‑width $h$ to satisfy
\begin{equation}
h < \frac{\varepsilon}{2},
\end{equation}
where $\varepsilon > 0$ is the assumed separation between the minimal and maximal observation times and the underlying phase boundary in~\eqref{eq_true_phase_boundary}.

On the other hand, excessively small patches may contain insufficient data for reliable model identification. Accordingly, we impose the condition
\begin{equation}
\left(\frac{2w}{\Delta x}+1 \right)\cdot \left(\frac{2h}{\Delta t}+1 \right) > 2K,
\end{equation}
where $K$ is the size of the feature dictionary, and $\Delta x$ and $\Delta t$ are the spatial and temporal sampling intervals, respectively.

With the covering constructed in~\eqref{eq_Pij}, we present an adaptive method for selecting the parameter $\alpha$ used to choose high‑CEE patches so that the complement of their union contains two connected components of $\Omega$. Unlike the phase‑transition detection step (Section~\ref{patch detection}), where patches intersecting spatial boundaries are excluded, here we select patches from the temporally interior region $\Psi$ defined in~\eqref{psi}. For each $\alpha \in \{0,1,\dots,100\}$, we set
\begin{equation}
\Psi_{\alpha} := \{\mathcal{P} \in \Psi \mid \operatorname{CEE}(\mathcal{P}) \geq P_\alpha\},
\end{equation}
where $P_\alpha$ is the $(100-\alpha)$-th percentile of the CEE values of patches in $\Psi$. Let $\mathcal{C}_{\alpha}$ be the topological interior of the largest connected component of $\bigcup_{\mathcal{P} \in \Psi_{\alpha}}\mathcal{P}$. We define the optimal $\alpha^*$ as the smallest $\alpha$ for which
\begin{equation}\label{eq_complement}
\Omega \setminus \mathcal{C}_{\alpha}
\end{equation}
consists of two connected components. Such an optimal value always exists. Notice that $\mathcal{C}_{\alpha_1} \subseteq \mathcal{C}_{\alpha_2}$ whenever $0 \leq \alpha_1 \leq \alpha_2 \leq 100$. Moreover, by construction, $\Omega \setminus \mathcal{C}_{100}$ is the disjoint union of $\mathcal{D} \times [0,2h]$ and $\mathcal{D} \times [T-2h,T]$. For noise‑free data, we begin the search at $\alpha = 0$; for noisy data, we start at $\alpha = 10$. We then increment $\alpha$ by one until~\eqref{eq_complement} splits into two connected components. Denote the resulting parameter by $\alpha^*$ and set $\mathcal{C} = \mathcal{C}_{\alpha^*}$.

\subsection{Construction of $\widetilde{\cC}$ via convexification}\label{Sec_convexification}
As discussed in Section~\ref{patch detection}, if a phase transition exists, the covering $\mathcal{C}$ constructed in Section~\ref{initial cover} is expected to concentrate around the phase boundary. To precisely locate the boundary, we use the numerical evolution strategy of Section~\ref{evolution section}. However, $\mathcal{C}$ may not be directly suitable for this purpose. Firstly, for some $x\in \cD$, the path for numerical evolution $L_x([0,T])\cap \cC$ may be disconnected. Secondly, there might exist $x\in \cD$ such that $\Gamma(x)\notin \cC$ due to the noise in our data.

\begin{figure}
    \centering
    \begin{tabular}{ccc}
    \includegraphics[width=0.23\textwidth]{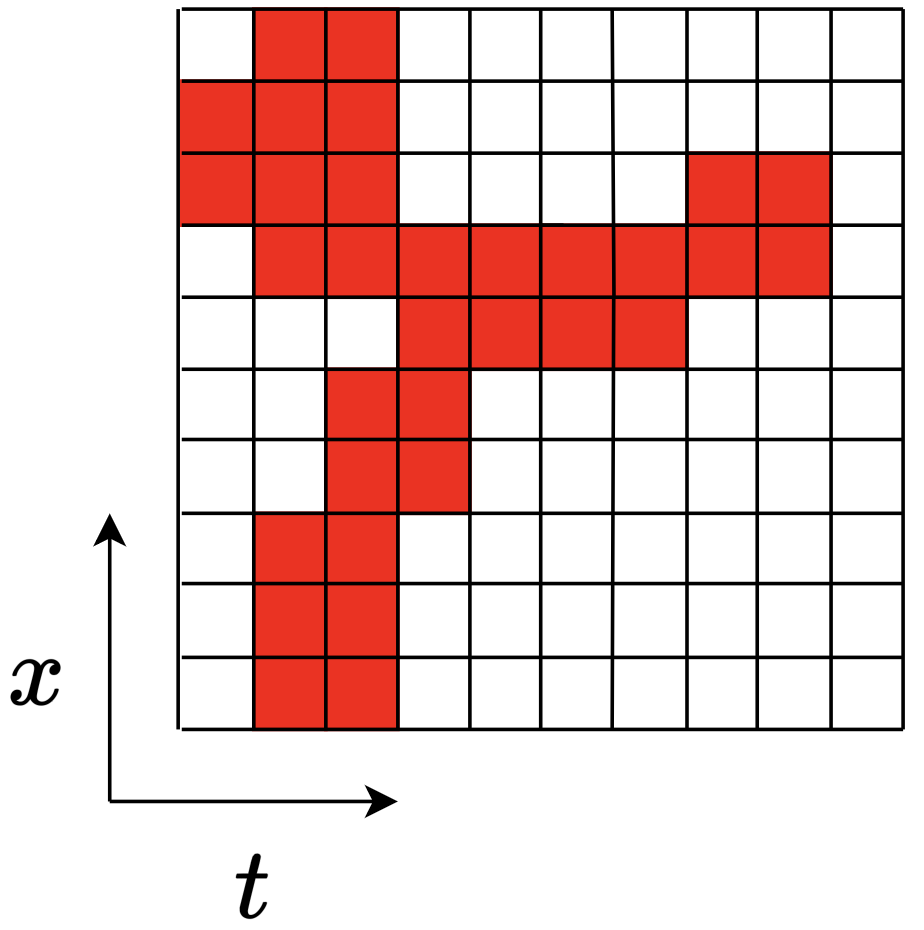}&
    \includegraphics[width=0.23\textwidth]{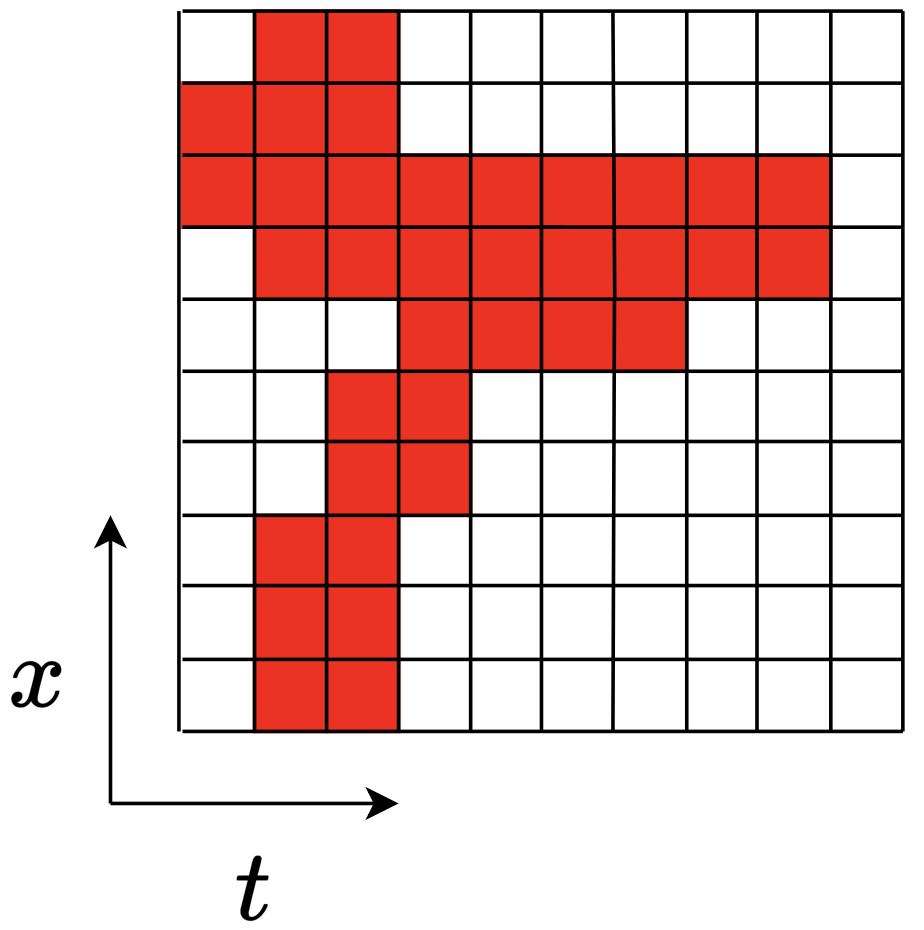}&
    \includegraphics[width=0.23\textwidth]{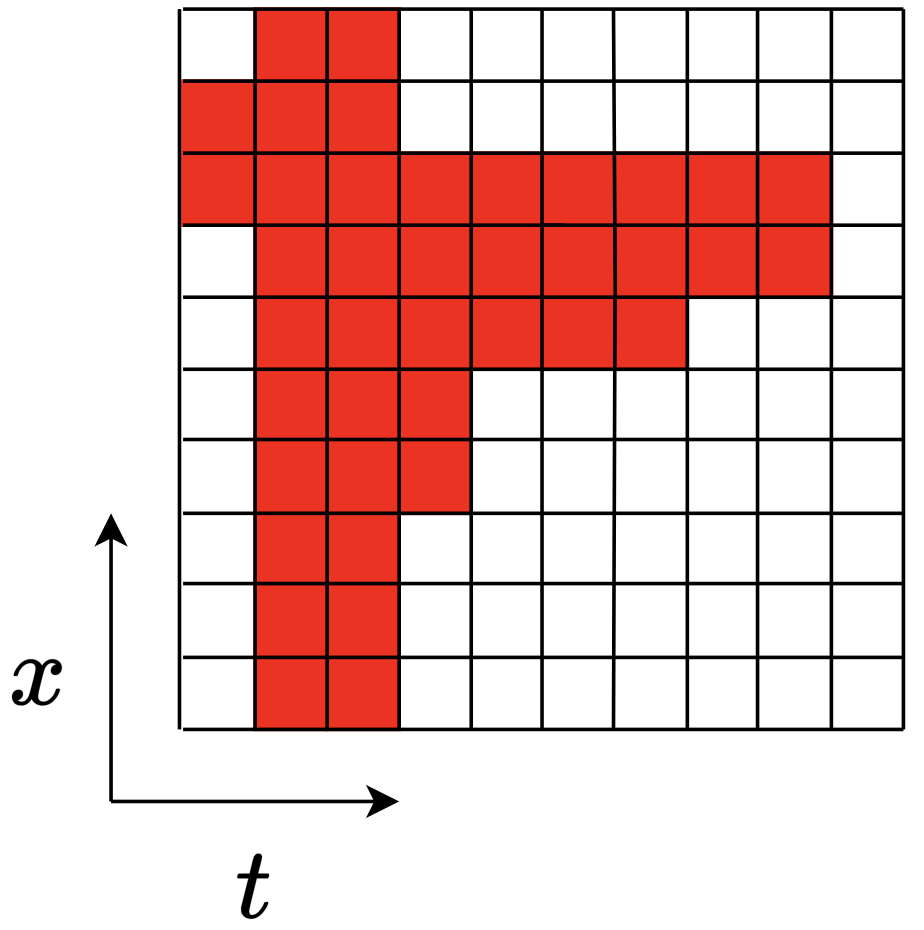}\\
    (a)&(b)&(c)
    \end{tabular}
    \caption{Schematic diagram of row and column padding in Section~\ref{Sec_convexification}. (a) The initial cover $\cC$. (b) Convexification in time, $\cC^\prime$ as defined in~\eqref{rowpadding}. (c) Convexification in space,  $\widetilde{\cC}$ as defined in~\eqref{columnpadding}.
    }

    \label{patches}
\end{figure}

To address these limitations, we propose to enlarge \(\mathcal{C}\) to obtain a set \(\widetilde{\mathcal{C}} \subset \Omega\) such that
\[
L_x([0,T]) \cap \widetilde{\mathcal{C}} = \bigl( \gamma^{\ell}(x),\, \gamma^{r}(x) \bigr)
\]
for some \(\gamma^{\ell}(x), \gamma^{r}(x) \in (0,T)\) with \(\gamma^{\ell}(x) < \gamma^{r}(x)\). The construction proceeds by successive convexification in time and then in space. For each \(x \in \mathcal{D}\), let  
\[
\mathcal{T}_x = \{ t \in [0,T] \mid (x,t) \in \mathcal{C} \}.
\]  
We define
\begin{equation}\label{rowpadding}
\mathcal{C}' := \bigcup_{x \in \mathcal{D}} \operatorname{Conv}_{\mathbb{R}}\bigl(\mathcal{T}_x\bigr),
\end{equation}
where \(\operatorname{Conv}_{\mathbb{R}}(A)\) denotes the convex hull of a set \(A \subset \mathbb{R}\).
Next, for each \(t \in [0,T]\) let  
\[
\mathcal{S}_t = \{ x \in \mathcal{D} \mid (x,t) \in \mathcal{C}' \}.
\]  
The final cover is constructed as
\begin{equation}\label{columnpadding}
\widetilde{\mathcal{C}} := \operatorname{int}\Bigl( \bigcup_{t \in [0,T]} \operatorname{Conv}_{\mathbb{R}}(\mathcal{S}_t) \Bigr),
\end{equation}
where \(\operatorname{int}(\cdot)\) denotes the topological interior. Figure~\ref{patches} illustrates this extension procedure. 

\subsection{Continuous representation for the phase boundary estimation}\label{bspline}
We select the vertical lines \(L_x\) described in Section~\ref{evolution section} so that they pass through the centers of the patches constructed in~\eqref{eq_Pij}. That is,
\[
x \in \mathcal{X} := \{x_{\min} + nw \mid n = 1,\dots,N_P\},
\]
where \(w\) is the half‑width of each patch and \(N_P\) is the number of patches in the spatial dimension. A pair of PMFs~\eqref{pdelta} and their Wasserstein barycenter~\eqref{ot} are then computed over the grid points
\begin{equation}
G = \{(x,\gamma^{\ell}(x) + m\Delta t) \in \widetilde{\mathcal{C}} \mid m = m_0,m_0+1,\dots, M_x - m_0,\; x \in \mathcal{X}\},
\end{equation}
with notation consistent with Section~\ref{boundary estimation}, where a positive integer $m_0$ is introduced to exclude the points near two ends of error sequences~\eqref{eq_sequences} to be candidate change points due to insufficient information. In this work, we set \(m_0 = \bigl\lfloor \frac{2h}{5\Delta t}\bigr\rfloor\), where \(h\) is the half‑width of each patch in the temporal dimension. Consequently, the collection of estimated phase‑boundary points given by~\eqref{eq_gamma_estimate} is a subset of \(G\).

To obtain a continuous representation of the phase boundary and its uncertainty, we proceed as follows. For each sampled spatial point \(x_i \in \mathcal{X}\), we have a discrete probability mass function \(p_{x_i}\) defined on a set of time indices within \(\widetilde{\mathcal{C}}\). To unify the domains of all \(p_{x_i}\) for interpolation, we first extend each \(p_{x_i}\) to the common interval
\begin{equation}\label{DT}
\mathcal{T} := [\gamma_{\min} + m_0\Delta t,\; \gamma_{\max} - m_0\Delta t],
\end{equation}
where \begin{equation}
\gamma_{\min} = \min_{1\le i \le N_P} \gamma^{\ell}(x_i), \qquad 
\gamma_{\max}= \max_{1\le i \le N_P} \gamma^{r}(x_i).
\end{equation}
The extension is performed by padding with zeros outside the original support of each \(p_{x_i}\). After this alignment, for every \(t \in \mathcal{T}\) we obtain a set of probability values \(\{p_{x_i}(t)\}_{i=1}^{N_P}\) over the spatial points \(x_i\).

We then interpolate these values along the spatial axis to construct a two‑dimensional probability density function 
\begin{equation}\label{eq_continuous_p}
\widetilde{p}: \mathcal{D} \times \mathcal{T}\to\mathbb{R}
\end{equation}
of the change point location. In this work we use piecewise cubic Hermite interpolating polynomial (PCHIP)~\cite{fritsch1984method}, which preserves monotonicity and avoids spurious oscillations~\cite{rabbath2019comparison}. Because extrapolation near the spatial boundaries \(\partial\mathcal{D} \times \mathcal{T}\) can occasionally produce small negative values, we truncate the interpolant at zero. This simple approach already yields satisfactory results in our experiments; more sophisticated schemes, such as positivity‑preserving piecewise rational cubic interpolation~\cite{butt1993preserving, hussain2008positivity}, can also be employed.

With the continuous density \(\widetilde{p}\) obtained, we fit a smooth parametric curve to the high‑probability region using B‑splines. Specifically, we represent the estimated phase boundary as
\begin{equation}\label{NB}
\widehat{\gamma}(x) := \sum_{j=1}^J \widehat{\beta}_j B_{j}(x), \qquad x \in \mathcal{D},
\end{equation}
where each \(B_j\) is a cubic B‑spline basis defined on a knot sequence \(x_{\min}=\tau_1 < \tau_2 < \dots < \tau_{J+4}=x_{\max}\) via the Cox–de Boor recurrence~\cite{de1972calculating}:
\begin{equation}
\begin{aligned}
B_{j,0}(x) &= 
\begin{cases}
1, & \tau_j < x \le \tau_{j+1},\\
0, & \text{otherwise},
\end{cases}\\[6pt]
B_{j,k}(x) &= \frac{x-\tau_j}{\tau_{j+k}-\tau_j} B_{j,k-1}(x) + \frac{\tau_{j+k+1}-x}{\tau_{j+k+1}-\tau_{j+1}} B_{j+1,k-1}(x), \quad k=1,2,3.
\end{aligned}
\end{equation}

Given a uniform grid \(\{(x_m,t_n)\}\) over \(\mathcal{D} \times \mathcal{T}\), the coefficients \(\widehat{\beta}_1,\dots,\widehat{\beta}_J\) are obtained by solving
\begin{equation}\label{Bspline fit}
\min_{\beta_1,\dots,\beta_J \in \mathbb{R}^J}
\sum_{m=1}^M \sum_{n=1}^N \widetilde{p}(x_m, t_n)\,
\Bigl(\sum_{j=1}^J \beta_j B_j(x_m) - t_n\Bigr)^2
+ \lambda N \sum_{m=1}^M \Bigl(\sum_{j=1}^J \beta_j \frac{dB_j}{d x}(x_m)\Bigr)^2,
\end{equation}
where \(\lambda > 0\) is a regularization parameter that controls the smoothness of the estimated boundary.

\section{Numerical Experiments}\label{sec5}
In this section, we present numerical experiments to validate the effectiveness and analyze the behavior of Phase‑IDENT. We test the proposed method on numerical solutions obtained from various combinations of linear and nonlinear PDEs. Let \(\mathcal{U} = \{U_{j}^n\}_{1 \leq j \leq M,\;1 \leq n \leq N}\) denote the discrete solution data, where the indices \(j\) and \(n\) correspond to spatial and temporal grid points, respectively. To simulate observational noise, for a given noise level \(\alpha > 0\), we add independent, mean‑zero Gaussian noise with standard deviation \(\sigma\), where
\begin{equation}
\sigma = \alpha\% \cdot \sqrt{\frac{1}{MN} \sum_{j=1}^{M}\sum_{n=1}^{N}\bigl(U_j^n - \tfrac{1}{2}(\max\mathcal{U}+\min\mathcal{U})\bigr)^2}\;.
\end{equation}
In our experiments, the PDE feature library contains all monomials of derivatives up to order 3 and multiplicative degree up to 2, resulting in a total of 15 candidate terms. Throughout the tests we use the following fixed parameters: the number of patches per dimension is \(N_P = 19\) (Section~\ref{initial cover}), the confidence level in~\eqref{confidence level} is \(p = 80\%\), the number of B‑spline bases in~\eqref{NB} is \(J = 26\), and the default value of the regularization parameter $\lambda$ in~\eqref{Bspline fit} is $0.1$.

To evaluate PDE identification performance, we employ the Jaccard similarity coefficient (JSC)~\cite{jaccard1912distribution}. Let \(\widehat{\mathcal{S}}\) and \(\mathcal{S}^*\) be the index sets of active features in the identified PDE and the ground‑truth PDE, respectively. The JSC is computed as
\begin{equation}\label{single jsc}
    \operatorname{JSC}\bigl(\widehat{\mathcal{S}},\mathcal{S}^*\bigr) := \frac{|\widehat{\mathcal{S}} \cap \mathcal{S}^*|}{|\widehat{\mathcal{S}} \cup \mathcal{S}^*|},
\end{equation}
where \(|\cdot|\) denotes set cardinality.

We also aim to precisely locate the phase boundary. Because information near the spatial boundaries \(\partial\mathcal{D} = \{x_{\min},x_{\max}\}\) may be less reliable, we quantify the accuracy of the phase boundary estimate on a restricted region
\begin{equation}\label{omegaprime}
    \Omega' := [x_{\min}+2w,\; x_{\max}-2w] \times \mathcal{T},
\end{equation}
where \(2w\) is the spatial width of a patch defined in~\eqref{eq_Pij} and \(\mathcal{T}\) is given in~\eqref{DT}. Given the ground‑truth boundary \(\Gamma: x \mapsto (x,\gamma(x))\)~\eqref{eq_true_phase_boundary} and its estimate \(\widehat{\Gamma}: x \mapsto (x,\widehat{\gamma}(x))\), the boundary error is computed over a discrete set \(\mathcal{X}'\) of spatial points in \(\Omega'\):
\begin{equation}\label{egammaeq} 
e_\Gamma := \frac{1}{|\mathcal{X}'|}\sum_{x \in \mathcal{X}'} \bigl|\widehat{\gamma}(x) - \gamma(x)\bigr|,
\end{equation}
where \(\mathcal{X}'\) is a uniform grid with spacing \(\Delta x' = (x_{\max} - x_{\min} - 4w)/80\).

\subsection{General performances}
\begin{figure}
    \centering
    \begin{tabular}{c@{\vspace{2pt}}c@{\vspace{2pt}}c@{\vspace{2pt}}c}
    &\includegraphics[width=0.33\linewidth]{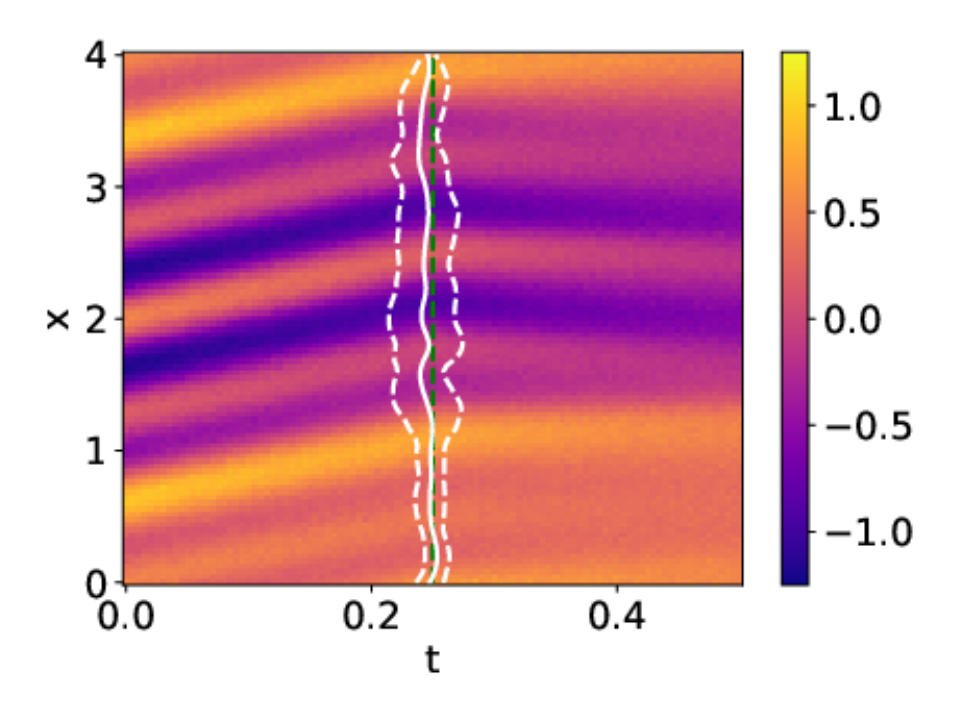}&
    \includegraphics[width=0.33\linewidth]{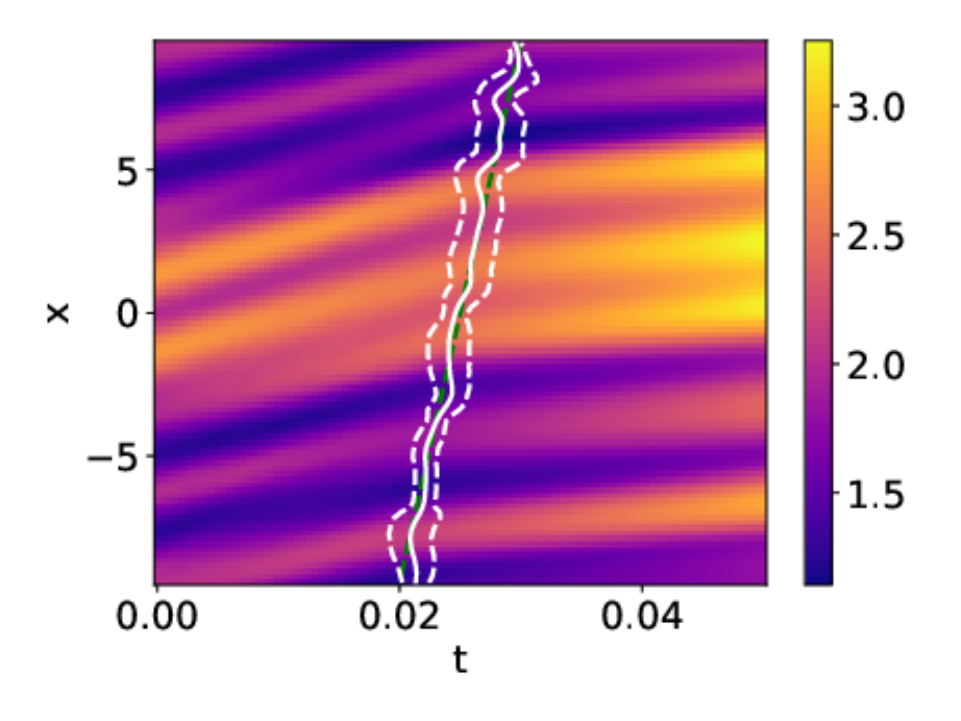}&
    \includegraphics[width=0.33\linewidth]{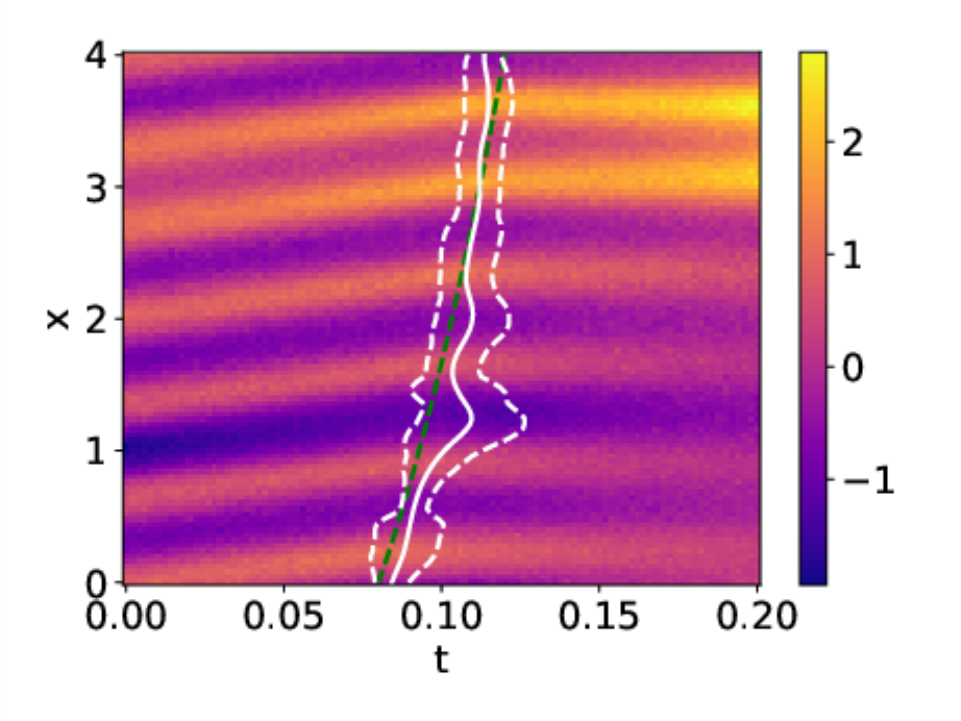}\\
    &(a)&(b)&(c)
    \end{tabular}
    \caption{Phase boundary estimation for three two‑phase PDE examples: (a) T$\to$VB~\eqref{gttVB} with $5\%$ additive noise; (b) KdV$\to$AR~\eqref{gtkdvar} with $1\%$ additive noise; (c) AR$\to$RD~\eqref{gtarrd} with $10\%$ additive noise. The green dashed line shows the true phase boundary $\Gamma$, the white solid line is the estimated boundary $\widehat{\Gamma}$, and the white dashed lines indicate the $80\%$ confidence intervals. In all three cases, the true $\Gamma$ lies within the corresponding confidence band.  }
    \label{boudaimage}
\end{figure}

We test Phase‑IDENT with the following two‑phase PDE:
\begin{equation}\label{gttVB}
\begin{cases}
\partial_t u(x,t)  = -2\,\partial_x u(x,t), &\text{for } 0 < x < 4,\; 0 < t < 0.25,\\[2pt]
\partial_t u(x,t)  = -u\,\partial_xu(x,t) + 0.1\,\partial_{xx}u(x,t), &\text{for } 0 < x < 4,\; 0.25 < t < 0.5,
\end{cases}
\end{equation}
which consists of a transport equation followed by a viscous Burgers equation (abbreviated as T$\to$VB). The initial condition is
\begin{equation}
    u(x,0)=\sin^2 (\pi x)\cos(0.5\pi x)+0.5 \cos (3\pi x).
\end{equation}
We compute the solution using Godunov’s method for the Burgers part and the Lax–Friedrichs scheme with explicit central differencing for the viscous part on a $2001\times2001$ grid. Then we add $5\%$ additive noise to the data. Figure~\ref{boudaimage}(a) shows the observed trajectory, with the green dashed line indicating the true phase boundary. Applying Phase‑IDENT yields the estimated boundary $\widehat{\Gamma}$ (white solid curve), which is close to the true boundary $\Gamma$. The white dashed curves illustrate the $80\%$ confidence interval, meaning that the probability of the true boundary lying between them exceeds $80\%$. The identified two‑phase PDE is
\begin{equation}\label{tvb}
\begin{cases}
\partial_t u(x,t)  = -1.954\,\partial_x u(x,t), &\text{for } 0 < x < 4,\; 0 < t < \widehat{\gamma}(x),\\[2pt]
\partial_t u(x,t)  = -0.963\,u\,\partial_xu(x,t) + 0.116\,\partial_{xx}u(x,t), &\text{for } 0 < x < 4,\; \widehat{\gamma}(x) < t < 0.5.
\end{cases}
\end{equation}
Comparing~\eqref{tvb} with~\eqref{gttVB} shows that Phase‑IDENT accurately recovers both the form and the coefficients of the underlying PDEs from noisy data.

Next, we test the method on a case where the phase boundary has a non‑zero slope. Consider
\begin{equation}\label{gtkdvar}
\begin{cases}
\partial_t u(x,t) = -30\,u\,\partial_x u(x,t) + 9\,\partial_{xxx} u(x,t), &\text{for } -3\pi < x < 3\pi,\; 0 < t < 0.025 + \frac{x}{600\pi},\\[2pt]
\partial_t u(x,t) = -37\,\partial_xu(x,t) + 10\,u(x,t), &\text{for } -3\pi < x < 3\pi,\; 0.025 + \frac{x}{600\pi} < t < 0.05,
\end{cases}
\end{equation}
which consists of a KdV equation followed by an advection‑reaction equation (abbreviated as KdV$\to$AR). The initial condition is
\begin{equation}\label{initialkdv}
\begin{aligned}
u(x,0) = &\;3A^2\sec^2\bigl(0.18A(x+10)\bigr) + 3B^2\sec^2\bigl(0.18B(x-9)\bigr)\\
         &+ 3C^2\sec^2\bigl(0.18C(x-3)\bigr) + 3D^2\sec^2\bigl(0.18D(x+3)\bigr),
\end{aligned}
\end{equation}
with $A=2.9$, $B=2.7$, $C=3.3$, $D=3.1$. The solution is computed using a spectral method for the KdV part and an upwind scheme for the advection‑reaction part on a $501\times501$ grid. The data are contaminated with $1\%$ additive noise. The estimated boundary $\widehat{\Gamma}$ is shown in Figure~\ref{boudaimage}(b); it again stays close to the true boundary, and the confidence interval remains relatively narrow. Because we employ $L^2$ regularization on the derivative in the boundary construction, $\widehat{\Gamma}$ exhibits slight oscillations around the true curve. Other regularizations, such as total variation, could be used to impose different priors. The identified PDE is
\begin{equation}\label{kdvar}
\begin{cases}
\partial_t u(x,t) = -31.533\,u\,\partial_x u(x,t) + 8.453\,\partial_{xxx} u(x,t), &\text{for } -3\pi < x < 3\pi,\; 0 < t < \widehat{\gamma}(x),\\[2pt]
\partial_t u(x,t) = -38.295\,\partial_xu(x,t) + 9.296\,u(x,t), &\text{for } -3\pi < x < 3\pi,\; \widehat{\gamma}(x) < t < 0.05.
\end{cases}
\end{equation}
Unlike the T$\to$VB case, although the correct features are selected, the coefficient estimates are less accurate. For the KdV component, similar sensitivity to noise for higher‑order terms has been observed in~\cite{he2022robust}. For the advection‑reaction component, the reaction‑term coefficient is inherently harder to recover precisely; nevertheless, the relative error of the obtained coefficient stays below $10\%$.

Phase‑IDENT remains effective even under higher noise levels. We consider an equation composed of an advection‑reaction part followed by a reaction‑diffusion part:
\begin{equation}\label{gtarrd}
\begin{cases}
\partial_t u(x,t) = -3\,\partial_x u(x,t) + 4\,u^2(x,t), &\text{for } 0 < x < 4,\; 0 < t < \gamma(x),\\[2pt]
\partial_t u(x,t) = 0.15\,\partial_{xx}u(x,t) + 6\,u^2(x,t), &\text{for } 0 < x < 4,\; \gamma(x) < t < 0.2,
\end{cases}
\end{equation}
where $\gamma(x)=\dfrac{2+\sqrt{3x+4}}{50}$, abbreviated as AR$\to$RD. The initial condition is
\begin{equation}
    u(x,0)=-\cos^2(\pi x)\sin(0.5\pi x)+0.8\cos(3\pi x),
\end{equation}
for $x \in [-3\pi, 3\pi)$. We solve the advection‑reaction part with an upwind scheme and the reaction‑diffusion part with explicit forward Euler time stepping, obtaining $1001\times1001$ discrete observations over $\Omega$. The data are then corrupted with $10\%$ additive noise. The estimated boundary $\widehat{\Gamma}$ and the ground truth are displayed in Figure~\ref{boudaimage}(c). Although the localization accuracy degrades at this noise level, the confidence interval provided by Phase‑IDENT still successfully encloses the true boundary.

The identified PDE is
\begin{equation}\label{arrd}
\begin{cases}
\partial_t u(x,t) = -2.885\,\partial_x u(x,t) + 4.026\,u^2(x,t), &\text{for } 0 < x < 4,\; 0 < t < \widehat{\gamma}(x),\\[2pt]
\partial_t u(x,t) = 0.143\,\partial_{xx}u(x,t) + 5.611\,u^2(x,t), &\text{for } 0 < x < 4,\; \widehat{\gamma}(x) < t < 0.2.
\end{cases}
\end{equation}
These results demonstrate that even at $10\%$ noise, the active features of the equations, including the reaction terms, can still be correctly identified. However, similar to the KdV$\to$AR case, the accuracy of the reaction‑term coefficient estimates is compromised under noisy conditions.

\subsection{Performances under different noise levels}\label{experiment}

In this section, we present a detailed study of the performance of Phase‑IDENT under different noise levels. We consider the two‑phase PDE consisting of the KdV equation followed by the inviscid Burgers equation:
\begin{equation}\label{gtKB}
\begin{cases}
\partial_t u(x,t) = -30\,u\,\partial_x u(x,t) + 9\,\partial_{xxx} u(x,t), & \text{for } -3\pi < x < 3\pi,\; 0 < t < \gamma(x),\\[6pt]
\partial_t u(x,t) = -21\,\partial_xu(x,t), & \text{for } -3\pi < x < 3\pi,\; \gamma(x) < t < 0.05,
\end{cases}
\end{equation}
where $\gamma(x) = 0.025 + \dfrac{x}{600\pi}$, 
abbreviated as KdV$\to$B. The numerical grid is $501\times501$, and the initial condition is again given by~\eqref{initialkdv}. We apply Phase‑IDENT to data corrupted by $1\%$, $2\%$, and $4\%$ additive Gaussian noise, performing 20 independent trials for each noise level. We evaluate the performance with respect to the following three aspects.\\
\begin{table}
\centering
\small
\renewcommand{\arraystretch}{1.3}
\begin{tabular}{|c|c|c|c|c|}
\hline
\multirow{2}{*}{} & \multicolumn{2}{c|}{$\widehat{\Omega}_1$} & \multicolumn{2}{c|}{$\widehat{\Omega}_2$}\\ 
\cline{2-5} 
 & Identified  & Average JSC & Identified & Average JSC \\
\hline
$1\%$ noise& $\begin{array}{cc}
     & \partial_tu=-30.560_{\pm0.184}u\partial_xu \\
     & +8.839_{\pm0.074}\partial_{xxx}u
\end{array}$ & $1.0$ & $\partial_tu=-22.096_{\pm0.299}u\partial_xu$ & $1.0$ \\ \hline
$2\%$ noise& $\begin{array}{cc}
     & \partial_tu=-32.337_{\pm0.341}u\partial_xu \\
     & +8.025_{\pm0.108}\partial_{xxx}u
\end{array}$ & $1.0$ & $\partial_tu=-21.536_{\pm0.377}u\partial_xu$& $1.0$ \\ \hline
$4\%$ noise& $\begin{array}{cc}
     & \partial_tu=-37.069_{\pm0.663}u\partial_xu \\
     & +5.633_{\pm0.257}\partial_{xxx}u
\end{array}$& $1.0$ & $\partial_tu=-23.549_{\pm0.625}u\partial_xu$& $0.975$ \\ \hline
\end{tabular}
\caption{PDE identification results for KdV$\to$B~\eqref{gtKB} generated by 20 random noisy data. Ground truth PDE in $\Omega_1$ and $\Omega_2$ are $\partial_tu=-30u\partial_xu+9\partial_{xxx}u$ and $\partial_tu=-21u\partial_xu$, respectively. The results show mean and standard deviation of the most frequently identified features coefficients.}
\label{identifiedcase3}
\end{table}

\noindent\textbf{Model identification.} In Table~\ref{identifiedcase3}, we report the most frequently identified domain‑specific models together with the mean and standard deviation of their estimated coefficients. We observe that for noise levels ranging from $1\%$ to $4\%$, the most frequently identified models consistently match the ground‑truth model. For the KdV component, Phase‑IDENT always recovers the correct model as the average JSC remains exactly $1.0$. The inviscid Burgers component is correctly identified in all trials at $1\%$ and $2\%$ noise; at $4\%$ noise the average JSC drops to $0.975$. As the noise level increases, the recovered coefficients generally become less accurate and exhibit slightly larger variability. These results indicate that identifying the true feature terms is more robust to noise than precisely estimating their coefficients.\\[5pt]
\begin{figure}
    \centering
    \begin{subfigure}[b]{0.3\textwidth}
        \centering
        \includegraphics[width=\linewidth]{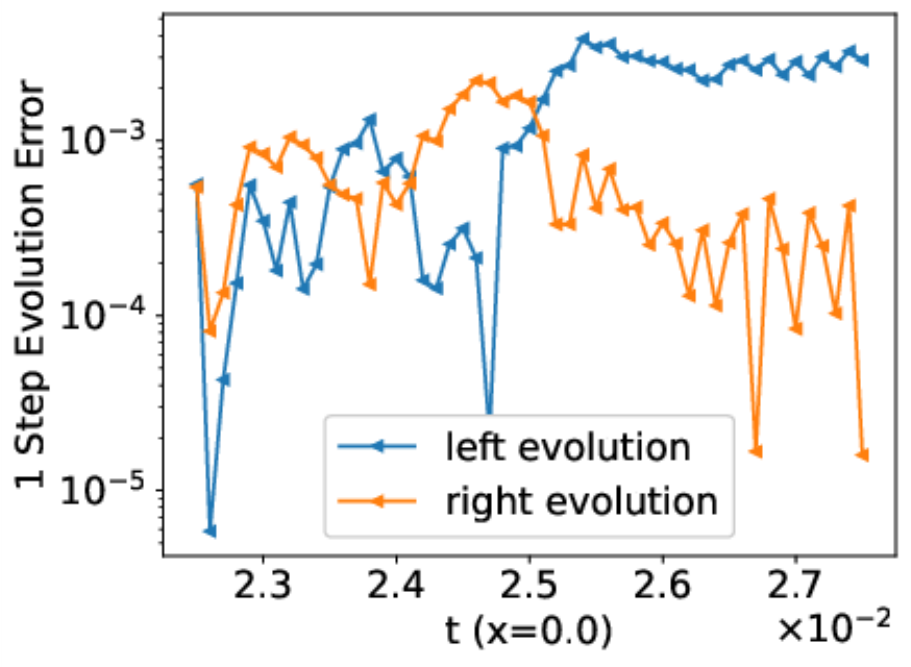}
        \caption{} 
        \label{fig:subfig2}
    \end{subfigure}
    \hfill
    \begin{subfigure}[b]{0.3\textwidth}
        \centering
        \includegraphics[width=\linewidth]{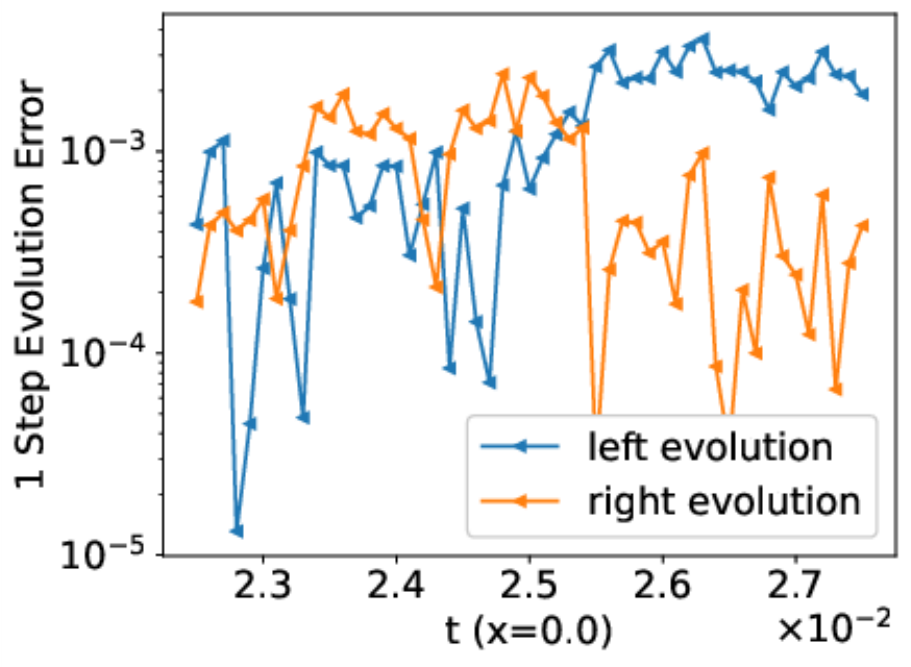}
        \caption{}
        \label{fig:subfig1}
    \end{subfigure}
    \hfill
    \begin{subfigure}[b]{0.3\textwidth}
        \centering
        \includegraphics[width=\linewidth]{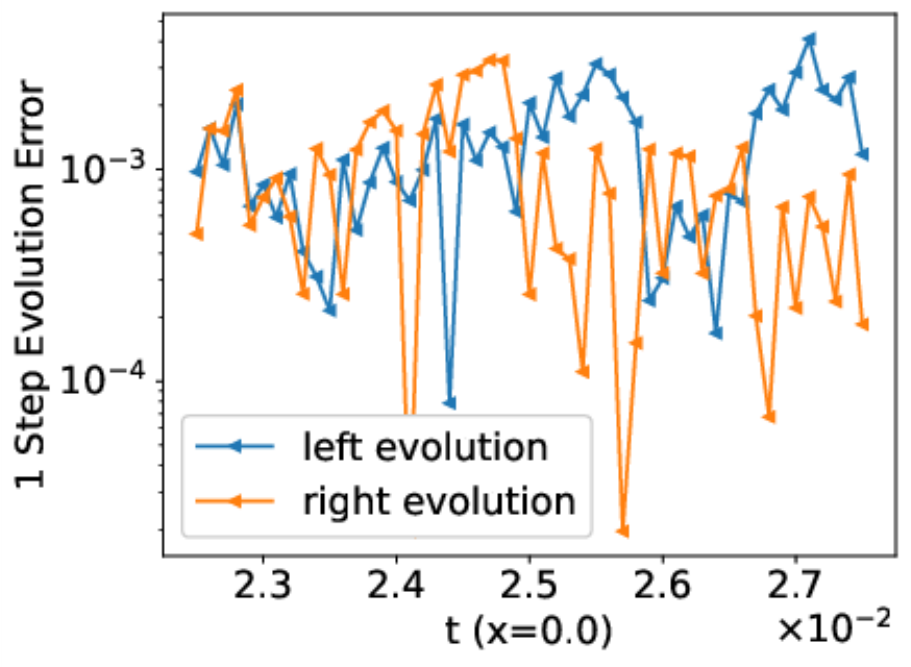}
        \caption{} 
        \label{fig:subfig3}
    \end{subfigure}

    \begin{subfigure}[b]{0.3\textwidth}
        \centering
        \includegraphics[width=\linewidth]{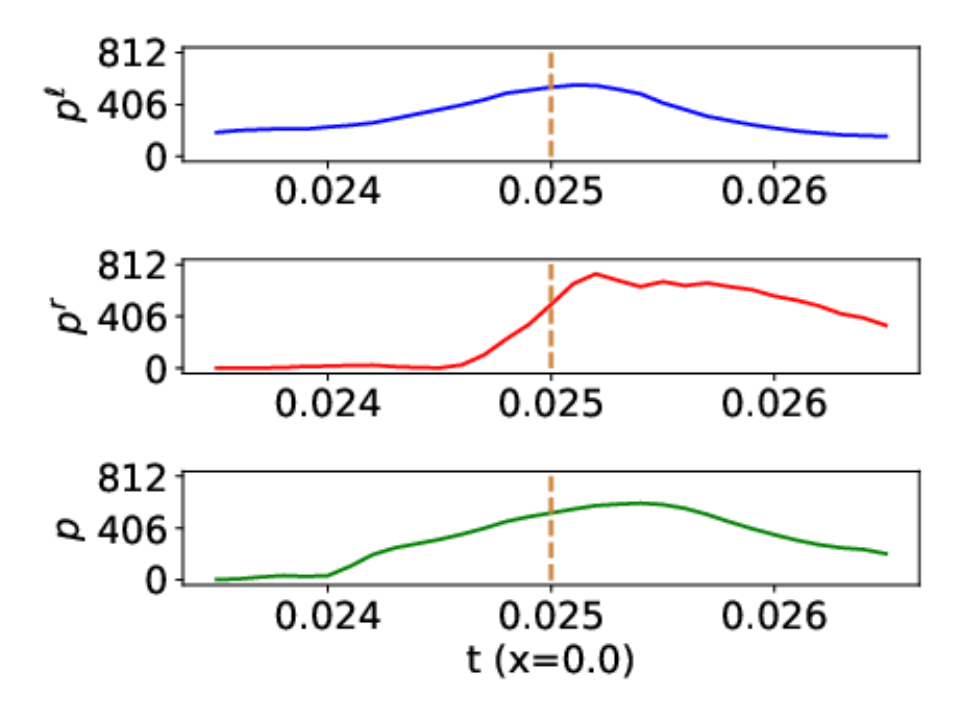}
        \caption{}
        \label{fig:subfig5}
    \end{subfigure}
    \hfill
    \begin{subfigure}[b]{0.3\textwidth}
        \centering
        \includegraphics[width=\linewidth]{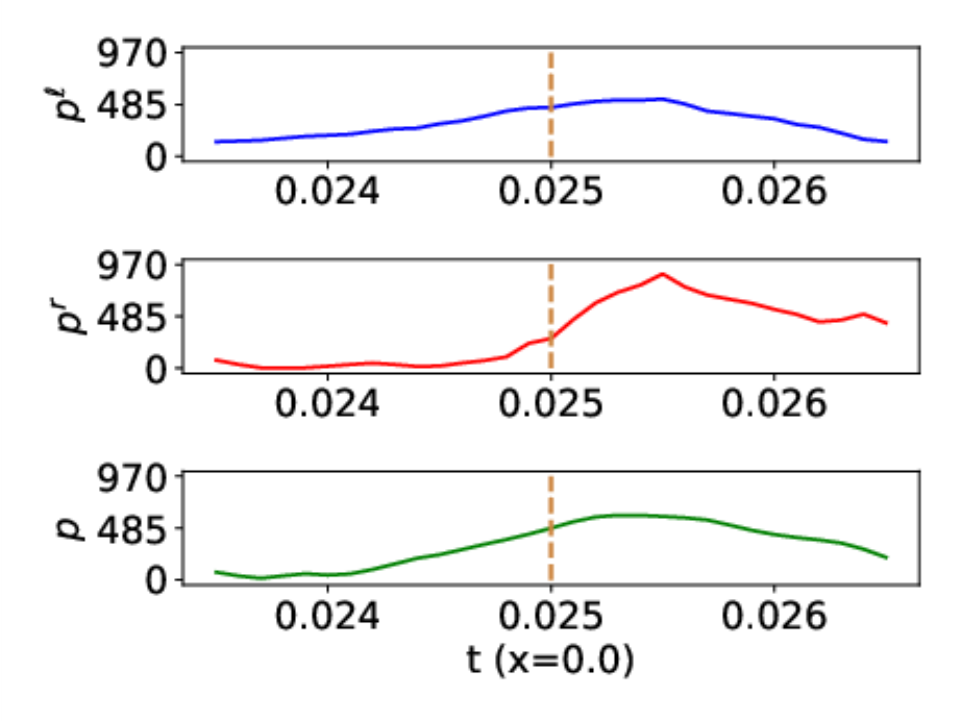}
        \caption{} 
        \label{fig:subfig4}
    \end{subfigure}
    \hfill
    \begin{subfigure}[b]{0.3\textwidth}
        \centering
        \includegraphics[width=\linewidth]{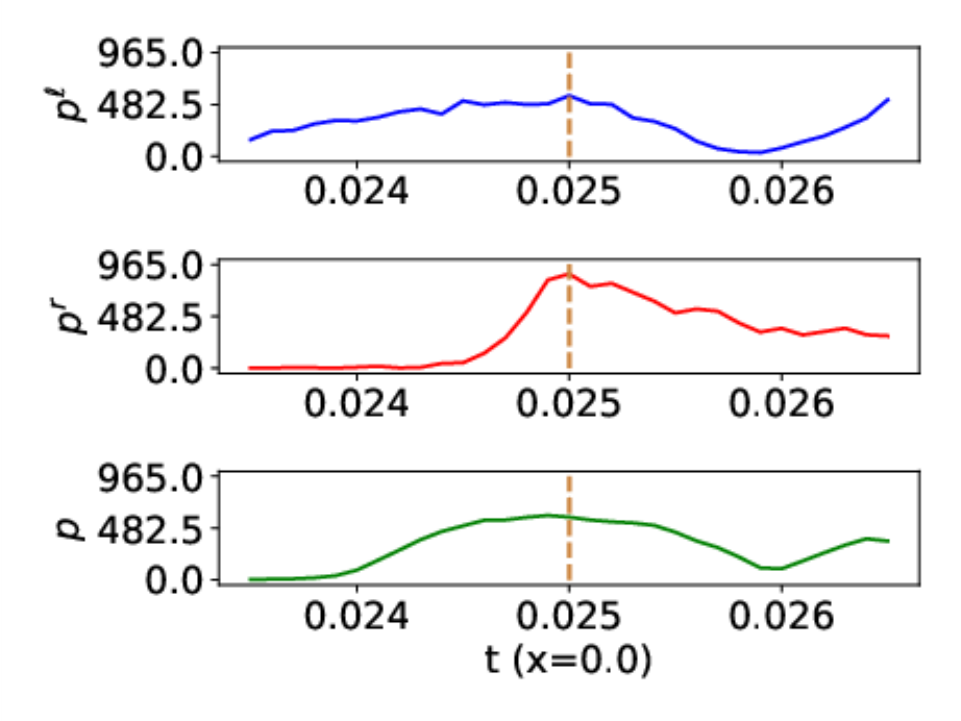}
        \caption{} 
        \label{fig:subfig6}
    \end{subfigure}
    
    \caption{Change point detection. Each column displays a pair of evolution error sequences~\eqref{eq_sequences} and the corresponding PMFs~\eqref{pdelta} and~\eqref{ot}. The first row shows the evolution errors for KdV$\to$B~\eqref{gtKB} under three noise levels: (a) $1\%$; (b) $2\%$; (c) $4\%$. The second row presents the change point detection results under the same noise levels: (d) $1\%$; (e) $2\%$; (f) $4\%$. The brown vertical lines indicate the true phase boundary location. Blue and red curves represent the discrete probability mass functions \(p_0^\ell\) and \(p_0^r\), respectively; the green curve is the Wasserstein barycenter \(p_0\) of these two distributions.}
    \label{CPD CUSUM}
\end{figure}
\noindent\textbf{Change point detection.} In Figure~\ref{CPD CUSUM}, we illustrate how noise affects the detection of change points as described in Section~\ref{boundary estimation}. In the first row (a)–(c), we show the single‑step evolution errors~\eqref{eq_sequences} for the spatial point \(x = (x_{\min} + x_{\max})/2 = 0\); in the second row (d)–(f), we display the corresponding probability mass functions \(p^\ell_0\) and \(p^r_0\) (defined in~\eqref{pdelta}) and their Wasserstein barycenter~\eqref{ot}.

We observe that at the $1\%$ noise level, the error curves obtained by evolving the identified domain‑specific models exhibit a noticeable jump near \(t = 2.5 \times 10^{-2}\). This abrupt change is captured by both probability mass functions \(p_0^\ell\) and \(p_0^r\), with the peak of \(p_0^r\) appearing slightly to the right of the true boundary. Consequently, the barycenter distribution shows a single, distinct peak close to the true phase‑boundary location.

At the $2\%$ noise level, both error curves exhibit more pronounced fluctuations, reflecting the sensitivity of numerical evolution to noise, yet the abrupt change remains visible in panel (b). In panel (e), the peak of \(p_0^\ell\) shifts slightly to the right of the true location because of oscillations in the evolution error around the change point; the peak of \(p_0^r\) shifts further, reaching its maximum near $0.0255$. Despite these individual shifts, the barycenter distribution obtained by combining \(p_0^\ell\) and \(p_0^r\) still peaks at the true phase boundary. This demonstrates the value of incorporating evolution information from both domains and synthesizing it through the barycenter.

When the noise increases to $4\%$, the change points in panel (c) become less distinct. Nevertheless, the peaks of the probability mass functions in panel (f) still successfully locate the true phase boundary. In particular, \(p_0^\ell\) displays two prominent peaks, which roughly correspond to two sudden increments in the blue curve of (c) that are induced by data noise. The change point detection based on the barycenter distribution~\eqref{ot} shows desirable robustness against such noise‑induced fluctuations.\\[2pt]

\begin{figure}
    \centering
    \subfloat[]{\includegraphics[width=0.24\textwidth]{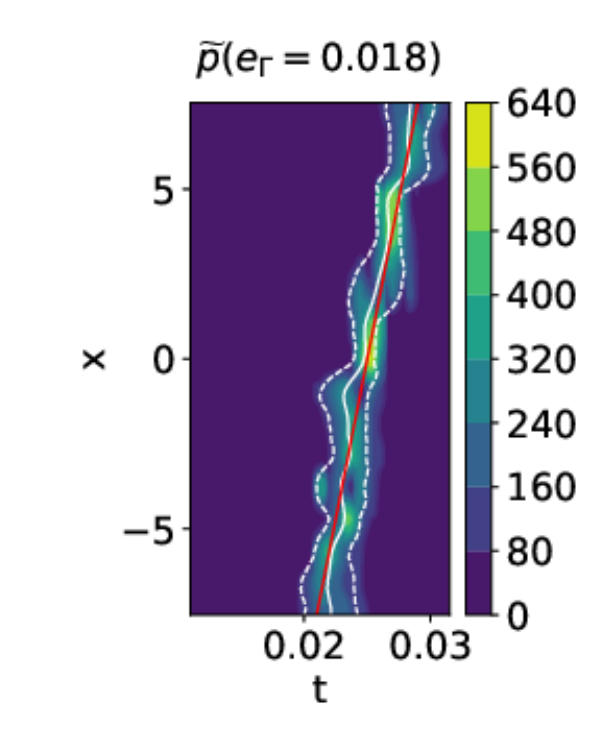}} 
    \subfloat[]{\includegraphics[width=0.24\textwidth]{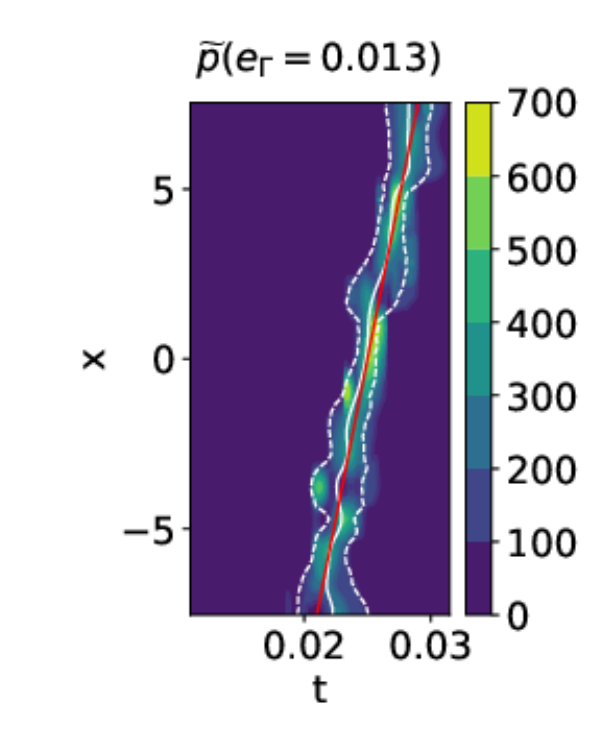}}
    \subfloat[]{\includegraphics[width=0.24\textwidth]{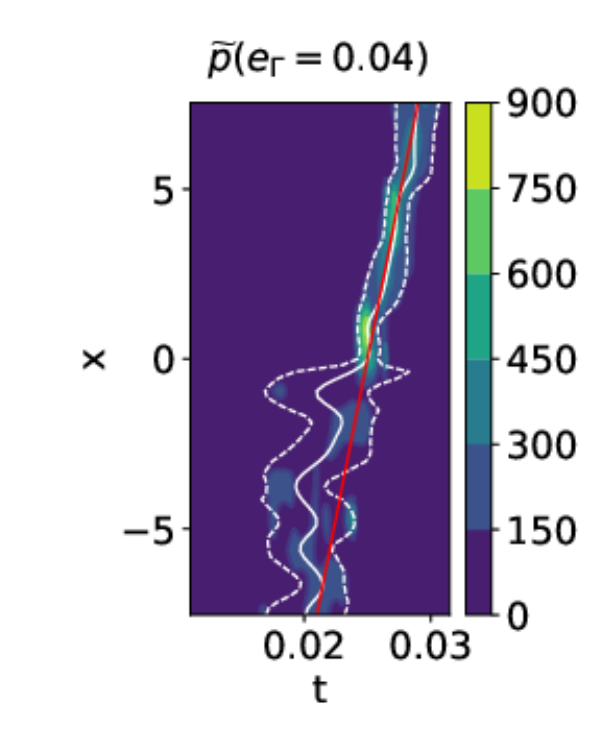}}
    \caption{Estimated phase boundary \(\widehat{\Gamma}\) and probability density function \(\widetilde{p}\)~\eqref{eq_continuous_p} over the region \(\Omega^\prime\) defined in~\eqref{omegaprime}: (a) KdV$\to$B with $1\%$ noise; (b) KdV$\to$B with $2\%$ noise; (c) KdV$\to$B with $4\%$ noise. The red line shows the true phase boundary \(\Gamma\); the solid white line is the estimated boundary \(\widehat{\Gamma}\); the dashed white lines mark the $80\%$ confidence interval; and \(e_\Gamma\) is the boundary‑estimation error computed via~\eqref{egammaeq}.}
    \label{Distribution}
\end{figure}

\noindent\textbf{Phase boundary estimation.} In Figure~\ref{Distribution}, we examine the phase boundary estimated by Phase‑IDENT under different noise levels. When the noise level is low ($1\%$ or $2\%$), the continuous estimate $\widehat{\Gamma}$ (white solid curve) stays close to the true phase boundary (red line). Furthermore, the $80\%$ confidence intervals (white dashed curves) remain relatively narrow while fully containing the true boundary.

As the noise level increases to $4\%$, $\widehat{\Gamma}$ displays visible oscillations and deviates from the ground truth. In particular, the upper part of $\widehat{\Gamma}$ still aligns well with the true boundary $\Gamma$, but in the lower part, where the probability density becomes more diffuse and is predominantly distributed to the left of $\Gamma$, the estimated boundary shifts significantly leftward.

Under higher noise levels, the high‑CEE patches are less concentrated, which leads to a wider covering $\widetilde{\mathcal{C}}$ (see~\eqref{eq_Ca_hat}). Consequently, the evolution sequences defined in~\eqref{eq_sequences} extend over longer intervals in some regions, making it harder to distinguish between fluctuations caused by noise and genuine change points induced by the phase boundary. Despite these challenges, the proposed $80\%$ confidence interval still covers most of the true phase boundary.\\[2pt]

\subsection{Effects of the phase boundary slope}
 We examine the effect of the phase boundary slope on the performance of Phase‑IDENT. Consider a parametric family of two‑phase PDEs
\begin{equation}\label{eq_effect_phase_boundary_slope}
\begin{cases}
\partial_t u = -\,u\,\partial_xu(x,t), &\text{for } 0 < x < 4,\; 0 < t < \gamma_s(x),\\[2pt]
\partial_t u = 4\,\partial_xu(x,t), &\text{for } 0 < x < 4,\; \gamma_s(x) < t < 0.5,
\end{cases}
\end{equation}
where \(\gamma_s(x) = s\,(x-2) + 0.25\) and the parameter \(s \in \mathbb{R}\) represents the slope of the phase boundary. Prior to \(\Gamma\), the system is governed by an inviscid Burgers equation; afterward, it follows a transport equation. We use the initial condition
\begin{equation}\label{inicond}
    u(x,0) = \sin^2(\pi x)\,\cos(0.5\pi x)\,\cos(0.3\pi x)
\end{equation}
for this Burgers‑to‑transport (B$\to$T) model. The solution is computed on a $1001\times1001$ grid and then uniformly downsampled in space to $201\times1001$. We compare results for noise levels of $1\%$ and $5\%$. Throughout these experiments, all hyper‑parameters of Phase‑IDENT are held fixed. For each noise level we perform $20$ independent trials and summarize the performance as follows.

\begin{table}
    \centering
\begin{tabular}{|c|c|c|c|c|}
\hline
\multirow{2}{*}{} & \multicolumn{2}{c|}{$1\%$ noise} & \multicolumn{2}{c|}{$5\%$ noise}\\ 
\cline{2-5} 
 & $\widehat{\Omega}_1$ & $\widehat{\Omega}_2$ & $\widehat{\Omega}_1$ & $\widehat{\Omega}_2$ \\
    \hline
    $s=0.03$ & $-0.947_{\pm 0.003}u\partial_xu$ & $3.992_{\pm 0.012}\partial_xu$ & $-0.796_{\pm 0.010}u\partial_xu$ & $3.901_{\pm 0.016}\partial_xu$ \\
    \hline
    $s=0.06$ & $-0.955_{\pm 0.002}u\partial_xu$ & $3.973_{\pm 0.012}\partial_xu$ & $-0.816_{\pm 0.016}u\partial_xu$ & $3.860_{\pm 0.017}\partial_xu$ \\
    \hline
    $s=0.08$ & $-0.962_{\pm 0.004}u\partial_xu$ & $3.990_{\pm 0.009}\partial_xu$ & $-0.824_{\pm 0.014}u\partial_xu$ & $3.868_{\pm 0.015}\partial_xu$ \\
    \hline
    $s=0.10$ & $-0.965_{\pm 0.007}u\partial_xu$ & $3.991_{\pm 0.005}\partial_xu$ & $-0.854_{\pm 0.024}u\partial_xu$ &$3.867_{\pm 0.016}\partial_xu$ \\
    \hline
\end{tabular}
\caption{PDE identification results for B$\to$T~\eqref{eq_effect_phase_boundary_slope} under $1\%$ and $5\%$ noise. The ground‑truth PDEs are $\partial_t u = -u\,\partial_x u$ in $\Omega_1$ and $\partial_t u = 4\,\partial_x u$ in $\Omega_2$. The results were obtained from 20 independent runs with randomly generated noisy data.}
\label{btgamma5}
\end{table}

In Table~\ref{btgamma5}, we report the most frequently identified models in the estimated phase domains as the slope \(s\) of the underlying phase boundary varies from \(0.03\) to \(0.1\). We observe that although increased noise worsens the accuracy of coefficient recovery, the correct features are still identified. These results indicate that the slope \(s\) does not strongly affect the model‑identification performance of Phase‑IDENT.

\begin{figure}
    \centering
    \begin{tabular}{c@{\vspace{2pt}}c@{\vspace{2pt}}c}
    \includegraphics[width=0.3\textwidth]{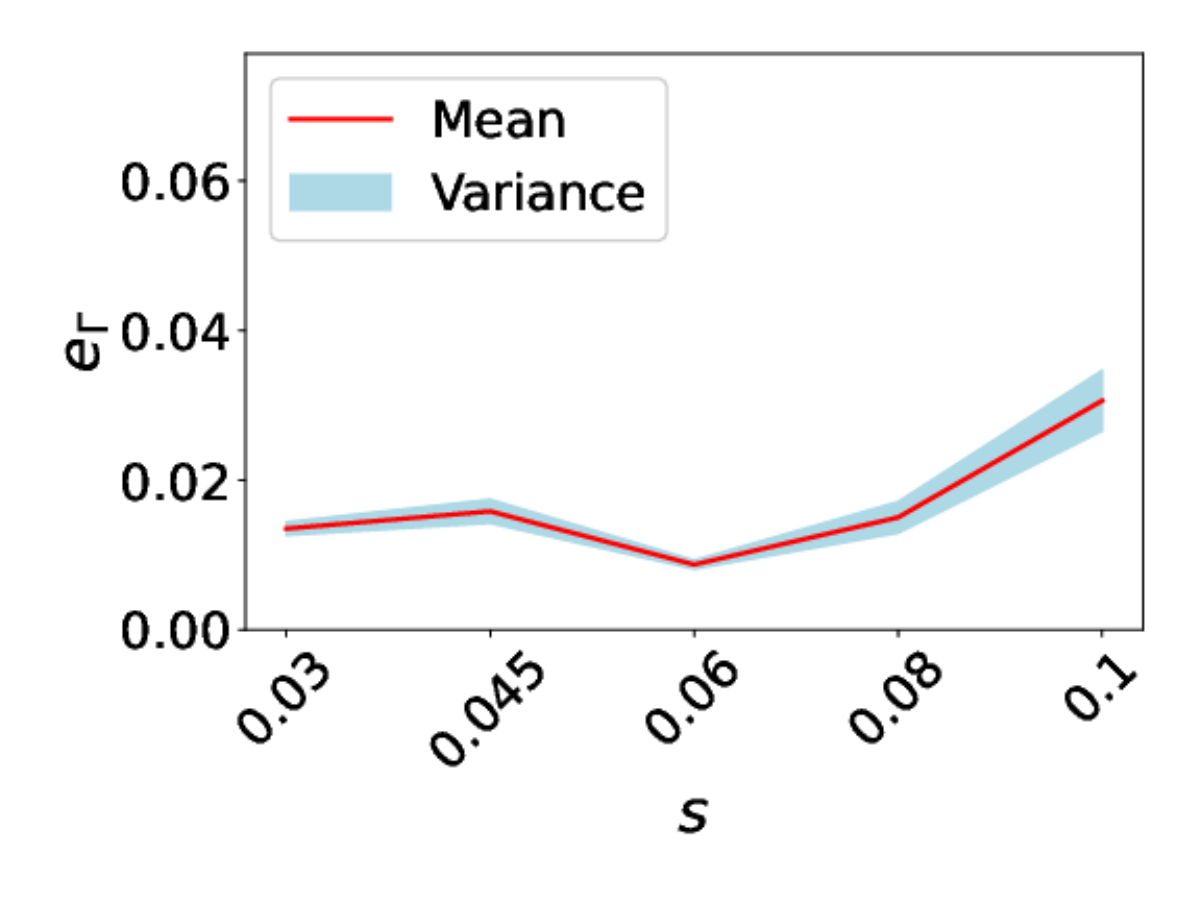}&
    \includegraphics[width=0.3\linewidth]{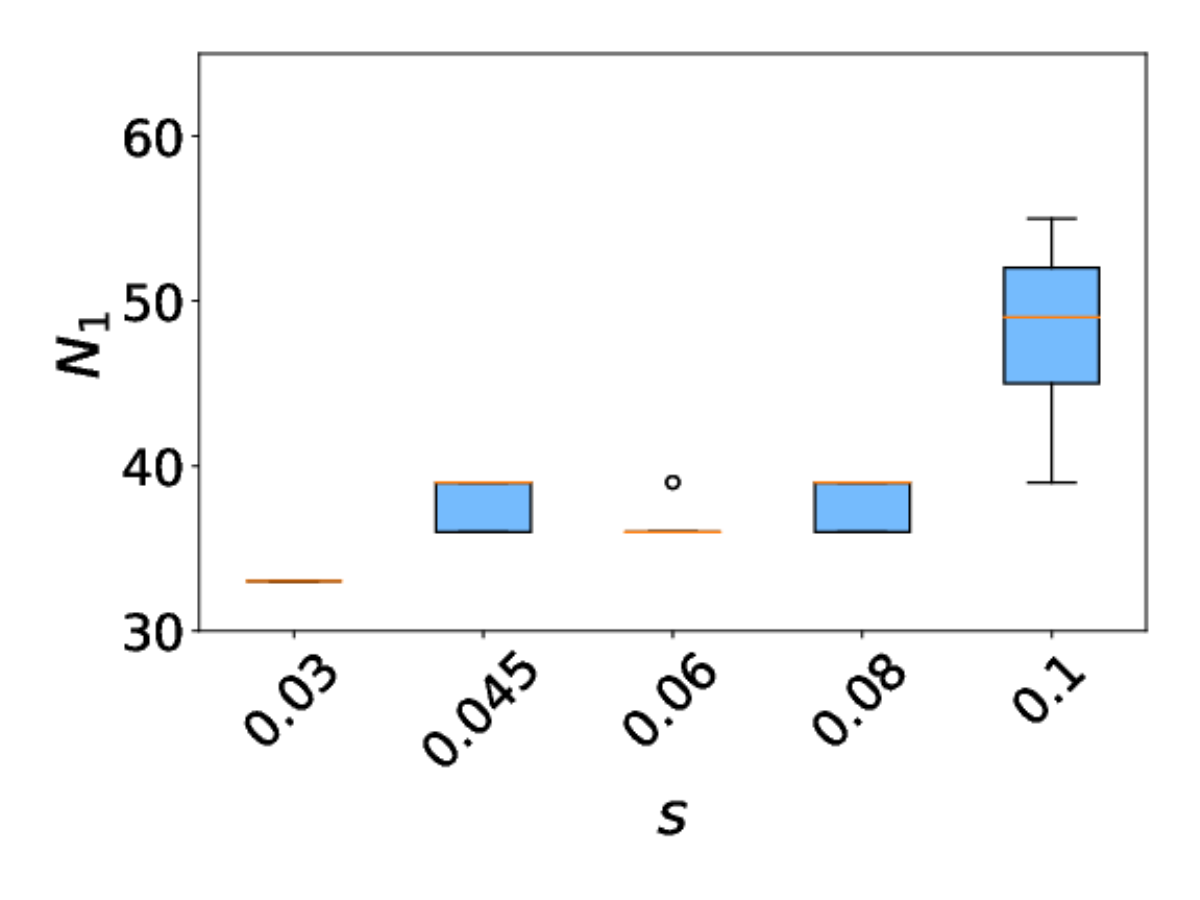}&

    \includegraphics[width=0.3\linewidth]{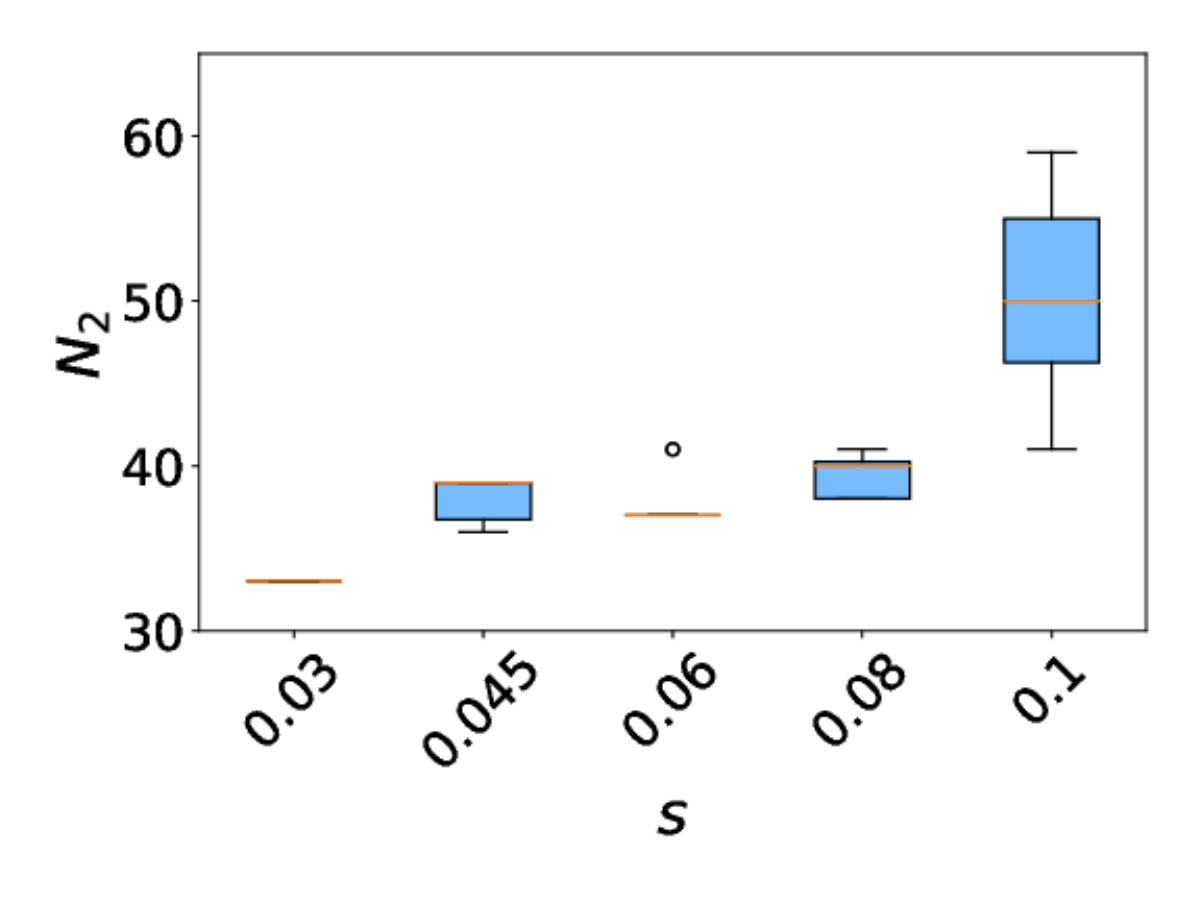}\\
    (a)&(b)&(c)\\
    \includegraphics[width=0.3\textwidth]{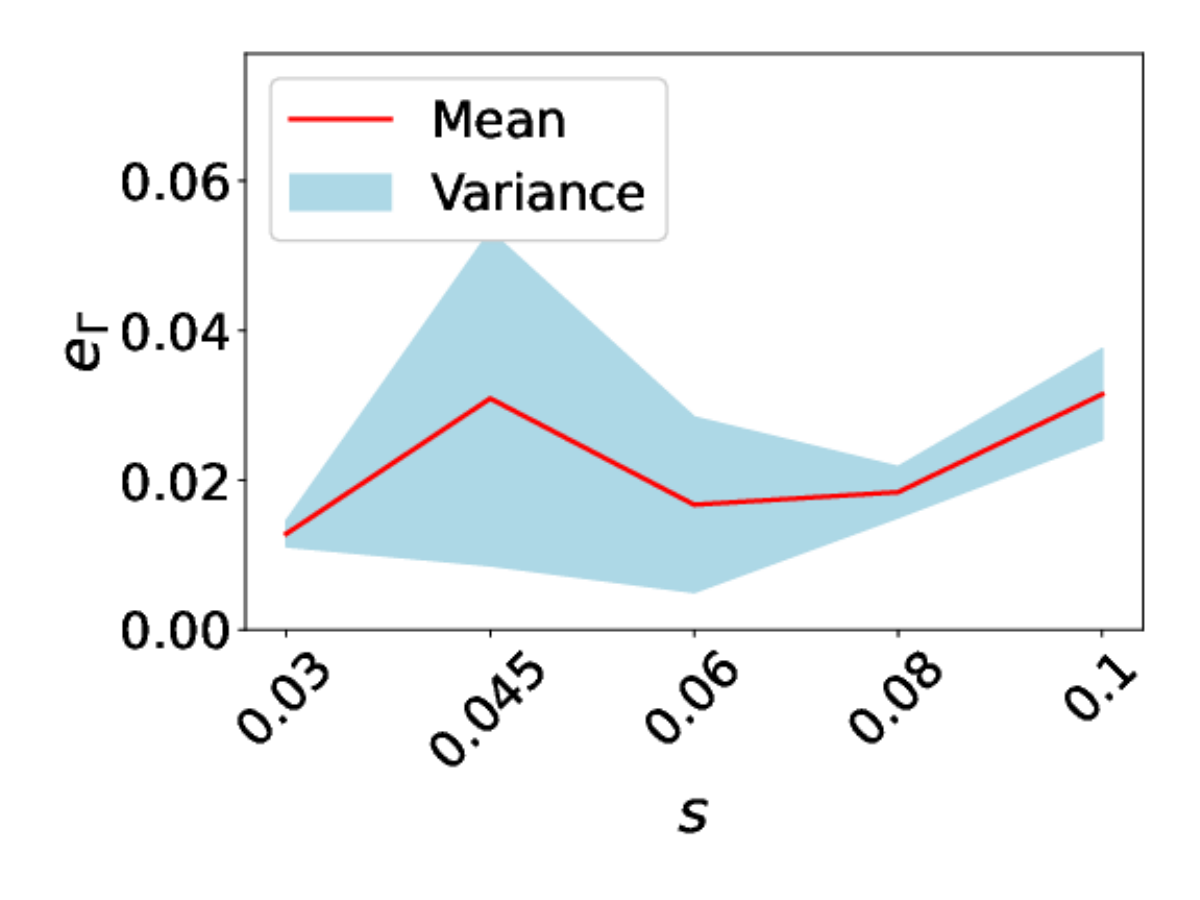}&
    \includegraphics[width=0.3\linewidth]{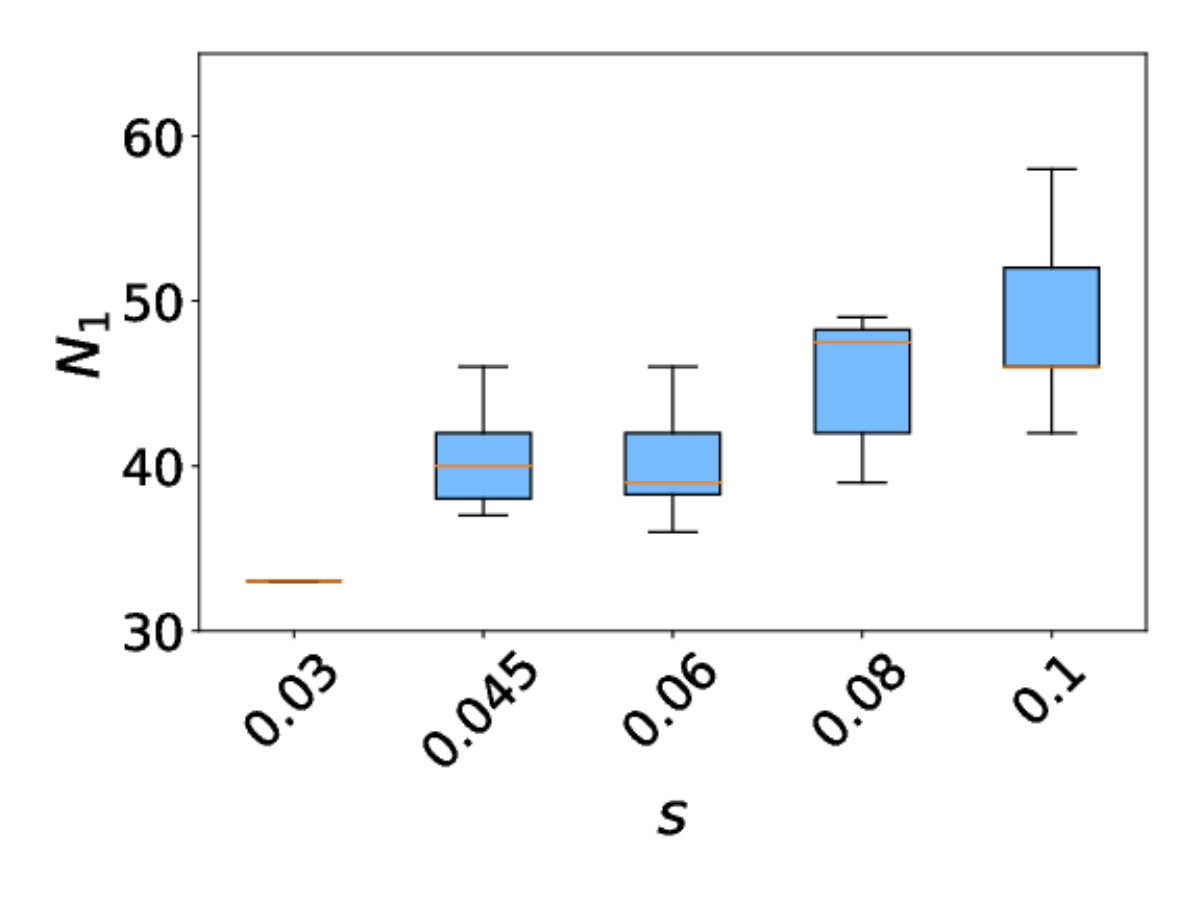}&
    
    \includegraphics[width=0.3\linewidth]{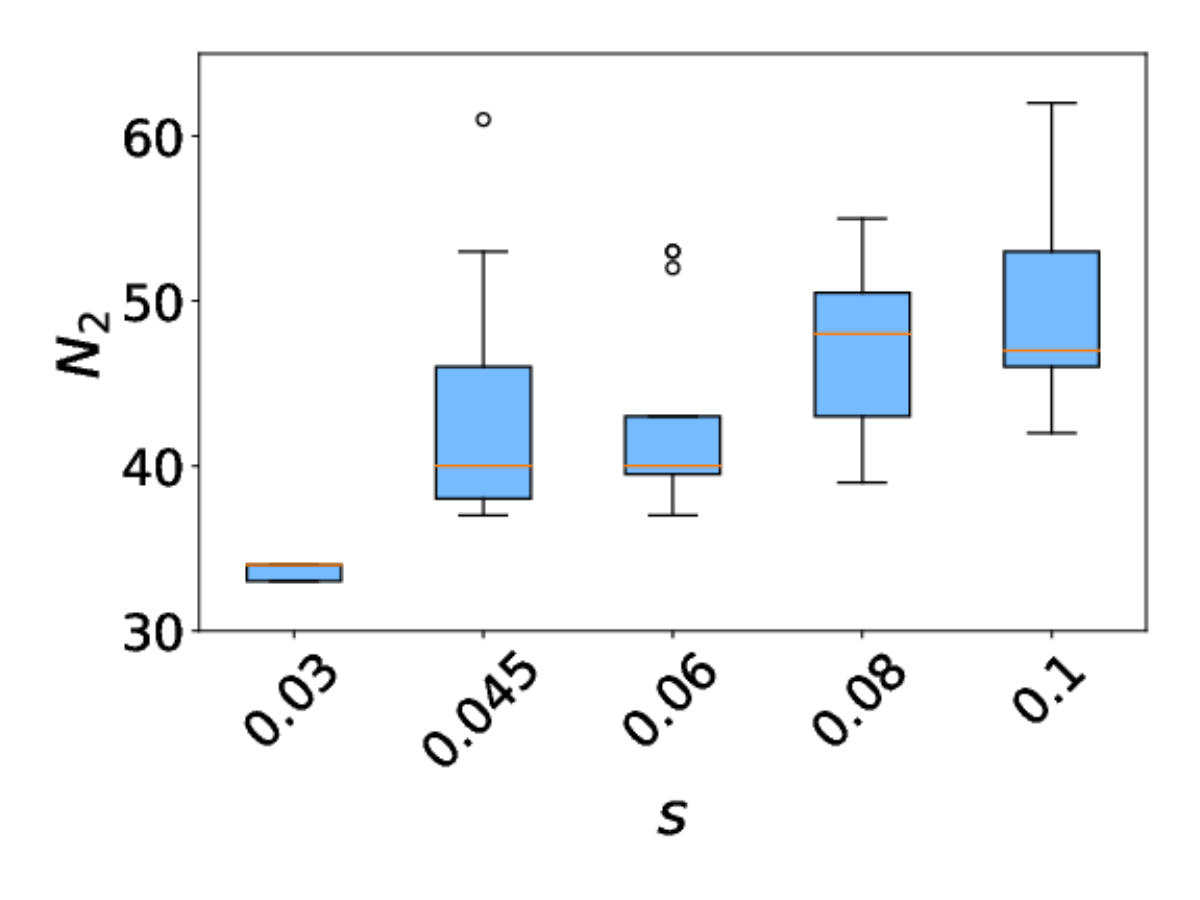}\\
    (d)&(e)&(f)
    
    \end{tabular}
    \caption{Effects of the phase‑boundary slope \(s\) on the accuracy of phase‑boundary estimation.  
(a) Boundary error \(e_\Gamma\) for different slopes under $1\%$ noise;   (b) Number of patches \(N_1\) in \(\mathcal{C}\) for different slopes under $1\%$ noise;  
(c) Number of patches \(N_2\) in \(\widetilde{\mathcal{C}}\) for different slopes under $1\%$ noise;  
(d) Boundary error \(e_\Gamma\) for different slopes under $5\%$ noise;  
(e) Number of patches \(N_1\) in \(\mathcal{C}\) for different slopes under $5\%$ noise;  
(f) Number of patches \(N_2\) in \(\widetilde{\mathcal{C}}\) for different slopes under $5\%$ noise.   The results are averaged over $20$ independent trials.}  
    \label{egammag}
\end{figure}

In Figure~\ref{egammag}(a) and (d), we show the mean and standard deviation of the boundary‑location error \(e_{\Gamma}\) defined in~\eqref{egammaeq} as the slope \(s\) varies. In general, the average error increases with \(s\). Additionally, we record the number of patches that constitute \(\mathcal{C}\) and \(\widetilde{\mathcal{C}}\) (constructed in~\eqref{columnpadding} of Section~\ref{sec4}), denoted by \(N_1\) and \(N_2\), respectively. Comparing Figure~\ref{egammag}(a,d) with (c,f), we see a consistent relationship between \(N_2\) and the boundary‑estimation error \(e_\Gamma\). Under \(1\%\) noise, only the case \(s = 0.1\) exhibits relatively large variance in \(e_\Gamma\), which corresponds to the widest spread of \(N_2\) in Figure~\ref{egammag}(c). When the noise level rises to \(5\%\), the estimation becomes more variable, especially for \(s = 0.045\). Accordingly, Figure~\ref{egammag}(f) shows the largest variance of \(N_2\) at \(s = 0.045\). From Figure~\ref{egammag}(b,c) we note that padding adds a relatively small number of patches under \(1\%\) noise, whereas under $5\%$ noise the increase from \(N_1\) to \(N_2\) in Figure~\ref{egammag}(e,f) is more pronounced. This occurs because high‑CEE patches become more scattered when the data are contaminated by stronger noise. Both a wider initial cover \(\mathcal{C}\) and the subsequent padding contribute to the marked increase in the variance of \(e_\Gamma\) from $1\%$ to $5\%$ noise.

These results numerically illustrate how the boundary‑estimation accuracy of Phase‑IDENT depends on both the boundary slope and the size of the constructed cover \(\widetilde{\mathcal{C}}\). From the viewpoint of PDE evolution, the slope of the phase boundary reflects the homogeneity of the transition time across space. When points near a given \(x \in \mathcal{D}\) transition at roughly the same time, i.e., when \(|\Gamma'(x)| \approx 0\), Phase‑IDENT provides robust and accurate boundary estimates. When the transition occurs more abruptly, the performance can be compromised.

\subsection{Effects of dynamics of the observed data}
The observed data strongly influence PDE identification, a relationship that has been systematically analyzed in~\cite{he2024much}. Here we investigate how the dynamics of the data, measured by a Sobolev‑type norm, affect phase boundary detection as described in Section~\ref{boundary estimation}.

Specifically, we examine the peaks of the probability mass functions \(p^\ell_x\) and \(p^r_x\) defined in~\eqref{pdelta}. When these peaks coincide, the estimated boundary location is unambiguous; otherwise, their discrepancy can lead to less accurate localization. To quantify this effect, we compute the normalized peak difference
\begin{equation}\label{epsilonx}
\epsilon(x):=\frac{\bigl|\argmax_{m} p^\ell_x(m)-\argmax_{m} p^r_x(m)\bigr|}{M_x-2m_0+1},
\end{equation}
for each \(x\in\mathcal{D}\).

For each \(x\), let \(\mathcal{I}_x\) be the set of indices of grid points inside the thin strip \([x-\Delta x, x+\Delta x]\times[\gamma^\ell(x),\gamma^r(x)]\). We then compute a Sobolev‑type norm as
\begin{equation}\label{ns}
n_s(x)\mathrel{\mathop:}=
\frac{1}{|\mathcal{I}_x|}
\sum_{(j,n)\in\mathcal{I}_x}
\bigl[(D_x U_{j}^n)^2+(D_t U_{j}^n)^2\bigr],
\end{equation}
where \(|\mathcal{I}_x|\) is the number of points in \(\mathcal{I}_x\), and \(D_x U_{j}^n\), \(D_t U_{j}^n\) denote spatial and temporal partial derivatives approximated via the 5‑point ENO scheme~\cite{harten1997uniformly} and forward differencing, respectively.

\begin{figure}
    \centering
    \begin{subfigure}[b]{0.3\textwidth}
        \centering
        \includegraphics[width=\linewidth]{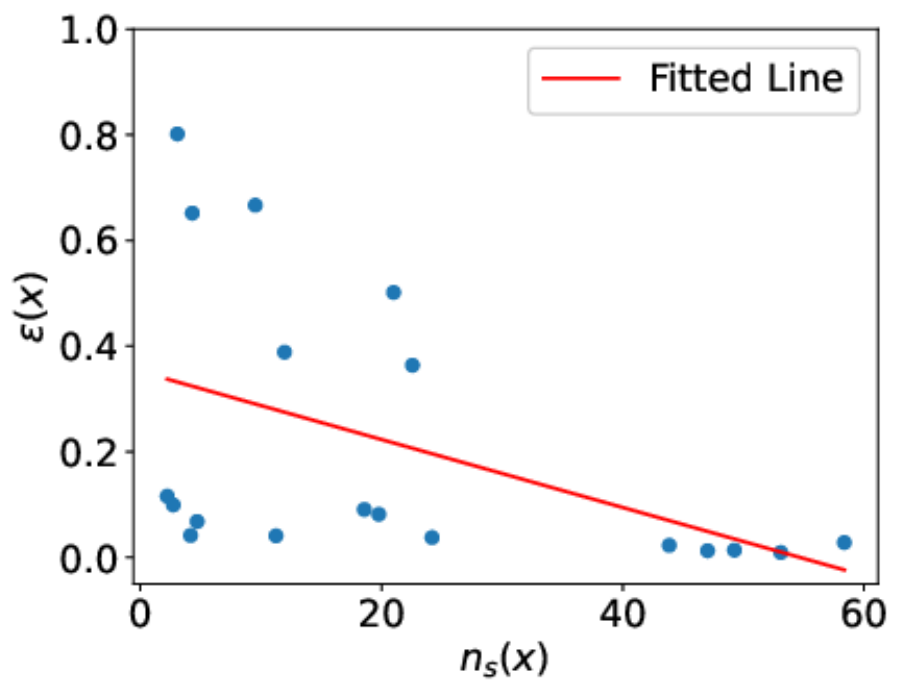}
        \caption{} 
        \label{fig:subfig9}
    \end{subfigure}
    \hfill
    \begin{subfigure}[b]{0.3\textwidth}
        \centering
        \includegraphics[width=\linewidth]{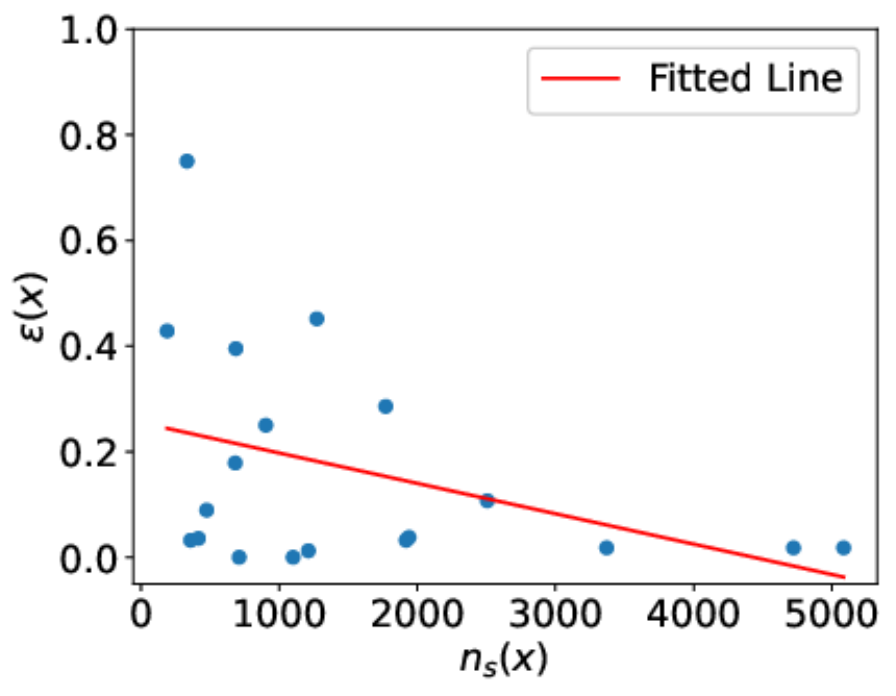}
        \caption{} 
        \label{fig:subfig9}
    \end{subfigure}
    \hfill    
    \begin{subfigure}[b]{0.3\textwidth}
        \centering
        \includegraphics[width=\linewidth]{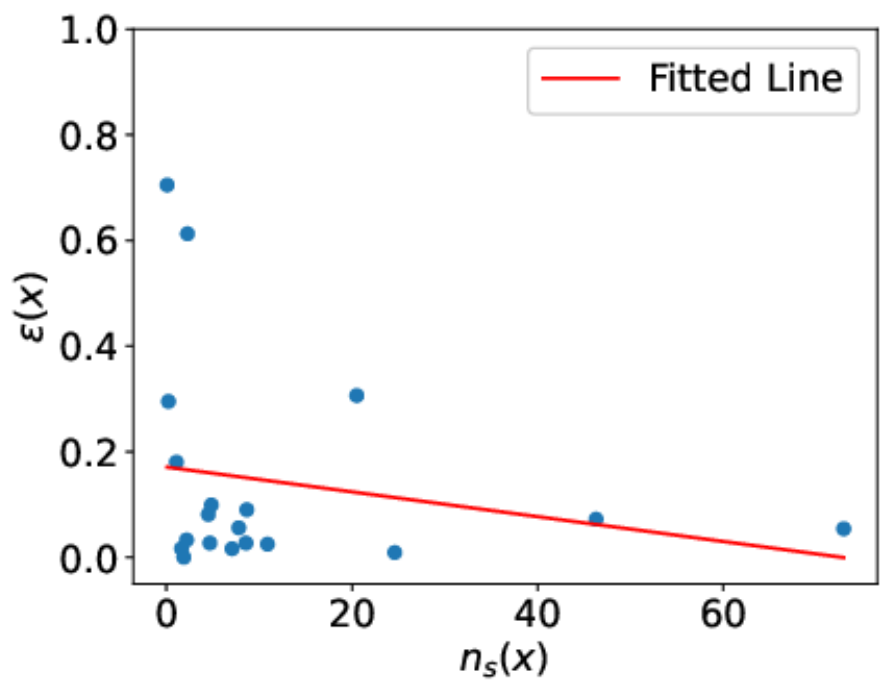}
        \caption{} 
        \label{fig:subfig9}
    \end{subfigure}
    
    \caption{Correlation between the relative change point detection discrepancy \(\epsilon(x)\) defined in~\eqref{epsilonx} and the approximate Sobolev norm \(n_s(x)\) defined in~\eqref{ns} under $1\%$ additive noise for three two‑phase PDEs: (a) T$\to$VB~\eqref{gttVB}; (b) KdV$\to$B~\eqref{gtKB}; and (c) B$\to$T~\eqref{gammacurve}. Regions where the discrepancy \(\epsilon(x)\) is large generally correspond to smaller values of the Sobolev norm \(n_s(x)\).}
    \label{Relation}
\end{figure}

Figure~\ref{Relation}(a)–(c) show a negative correlation between the change point detection discrepancy \(\epsilon(x)\) and the Sobolev norm \(n_s(x)\) for three test cases: T$\to$VB~\eqref{gttVB}, KdV$\to$B~\eqref{gtKB}, and the following B$\to$T
\begin{equation}\label{gammacurve}
\begin{cases}
\partial_t u(x,t) = -\,u\,\partial_x u(x,t), & \text{for } 0<x<4,\;0<t<\gamma(x),\\[2pt]
\partial_t u(x,t) = 4\,\partial_{x}u(x,t), & \text{for } 0<x<4,\;\gamma(x)<t<0.2,
\end{cases}
\end{equation}
with \(\gamma(x)=(2+\sqrt{3x+4})/50\) (a curved boundary used here to study the effect of spatial variation) and the initial condition~\eqref{inicond}. This result is expected: when the trajectory data contain less dynamical information indicated by a lower Sobolev norm, the magnitudes of the derivative features are generally smaller, making the change point detection more sensitive to noise. 

\subsection{Effects of patch size}\label{patchsizeeffect}
We investigate how patch size affects two aspects: (i) the CEE value on each patch, and (ii) the recovery of active features from each patch. These experiments are performed on the T$\to$VB case defined in~\eqref{gttVB}. For a patch $\mathcal{P}$, let $N_{\mathcal{P}}$ denote the number of grid points it contains. We define the average CEE as
\begin{equation}\label{averagecee}
    \overline{\operatorname{CEE}}(\mathcal{P}) := \frac{\operatorname{CEE}(\mathcal{P})}{N_{\mathcal{P}}},
\end{equation}
where $\operatorname{CEE}(\mathcal{P})$ is computed as in~\eqref{CEE formu}.

Because a patch may lie entirely inside $\Omega_1$, entirely inside $\Omega_2$, or intersect both regions, we assess its PDE identification performance as follows. Let $\mathcal{S}_{\Omega_1}$ and $\mathcal{S}_{\Omega_2}$ be the index sets of active features in the ground‑truth PDEs for $\Omega_1$ and $\Omega_2$, respectively, and let $\mathcal{S}_{\mathcal{P}}$ denote the indices of active features identified from $\mathcal{P}$. The patch‑wise JSC is defined by
\begin{equation}\label{patchwise jsc}
\operatorname{JSC}\bigl(\mathcal{P}\bigr) := \max\!\Bigl\{
\operatorname{JSC}(\mathcal{S}_{\Omega_1}, \mathcal{S}_{\mathcal{P}}),\;
\operatorname{JSC}(\mathcal{S}_{\Omega_2}, \mathcal{S}_{\mathcal{P}})\Bigr\},
\end{equation}
where the JSC between two index sets is computed via~\eqref{single jsc}. A patch‑wise JSC of $0$ means that none of the correct features were identified on that patch, whereas a value of $1$ indicates that either the PDE type of $\Omega_1$ or that of $\Omega_2$ was perfectly recovered.

We consider coverings with patch sizes $10\times 10$, $40\times 40$, $100\times 100$, $200\times 200$, and $400\times 400$. For each size, we randomly select $20$ patches that intersect the phase boundary; these are called \textit{inconsistent patches}. For comparison, we also select $10$ patches from $\Omega_1$ and $10$ from $\Omega_2$ that do not intersect the phase boundary; these are referred to as \textit{consistent patches}. For every sampled patch we compute the average CEE~\eqref{averagecee} and the patch‑wise JSC~\eqref{patchwise jsc}.

\begin{figure}
    \centering
    \begin{tabular}{c@{\vspace{2pt}}c@{\vspace{2pt}}c}
    \includegraphics[width=0.33\linewidth]{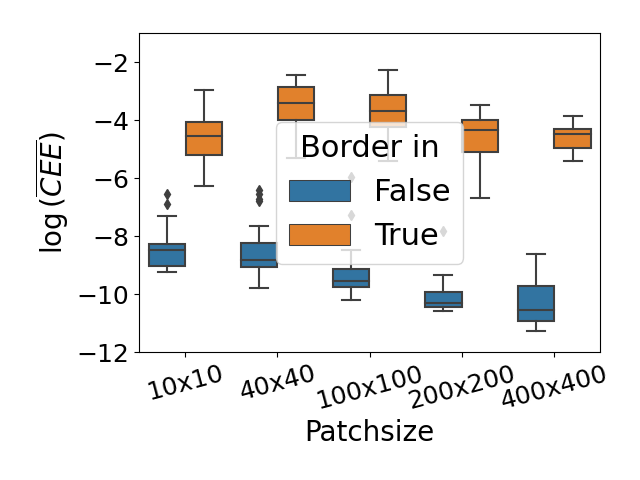}&
    \includegraphics[width=0.33\linewidth]{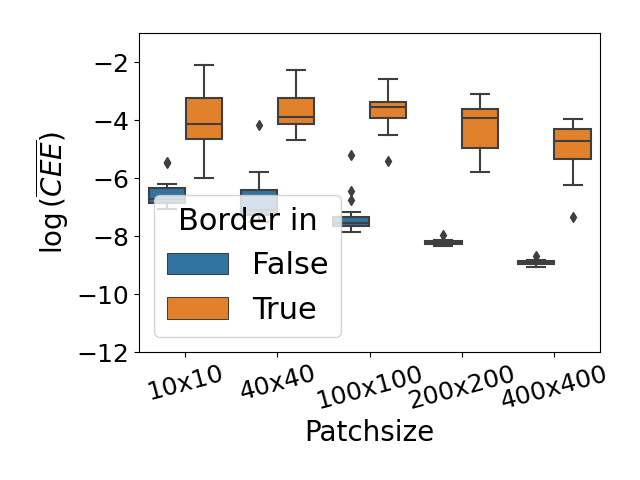}&
    \includegraphics[width=0.33\linewidth]{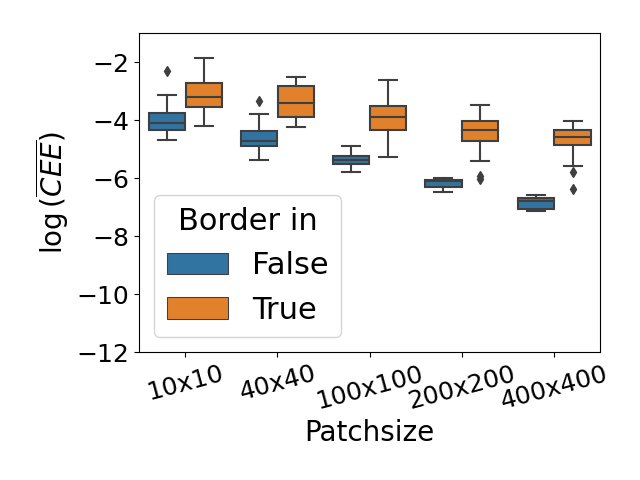}\\
    (a)&(b)&(c)\\
    \end{tabular}
    \caption{Patch‑wise average CEE values~\eqref{averagecee} versus patch size for different noise levels: (a) $0.1\%$, (b) $1\%$, (c) $10\%$. Orange and blue boxes show the distributions of average CEE over inconsistent and consistent patches, respectively. At a fixed patch size, higher noise reduces the CEE gap between consistent and inconsistent patches. For a fixed noise level, larger patch sizes widen this CEE difference, which aids the construction of the phase boundary covering $\widetilde{\mathcal{C}}$.
    }
    \label{CEE PATCH}
\end{figure}

Figure~\ref{CEE PATCH} displays the relationship between average CEE (in log scale) and patch size for different noise levels in panels (a)–(c). The results show that larger patch sizes generally amplify the difference in average CEE values between consistent and inconsistent patches. As the noise level rises, the average CEE values increase overall, while the distinction between the two patch types becomes less pronounced.

\begin{figure}
    \centering
    \begin{subfigure}[b]{0.98\textwidth}
        \centering
        \includegraphics[width=\linewidth]{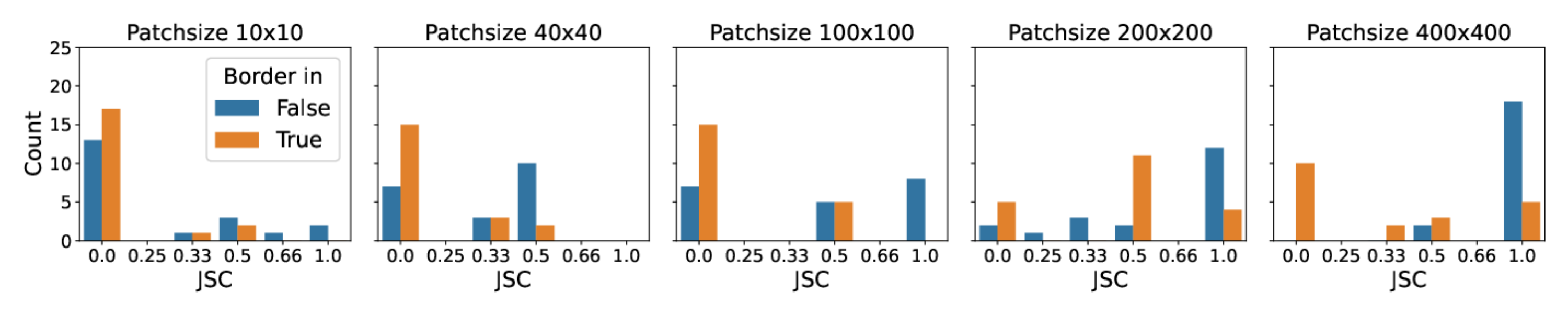}
        \caption{}
        \label{fig:subfig1}
    \end{subfigure}
    \hfill
    \begin{subfigure}[b]{0.98\textwidth}
        \centering
        \includegraphics[width=\linewidth]{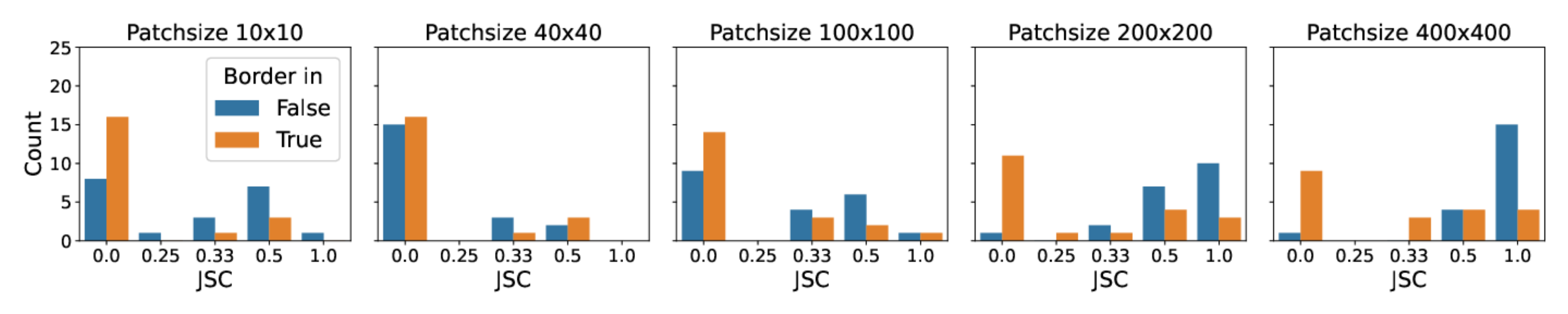}
        \caption{}
        \label{fig:subfig2}
    \end{subfigure}
    \hfill
    \begin{subfigure}[b]{0.98\textwidth}
        \centering
        \includegraphics[width=\linewidth]{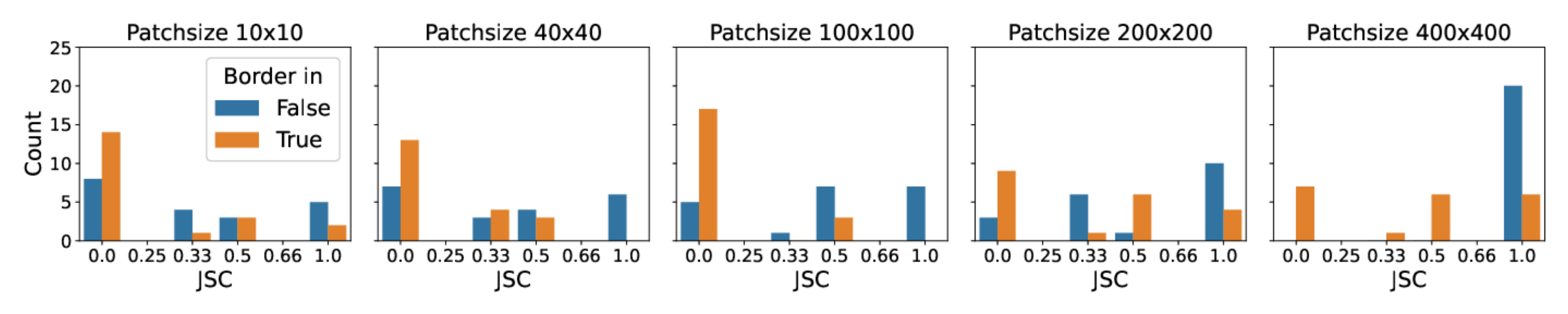}
        \caption{} 
        \label{fig:subfig4}
    \end{subfigure}
    \hfill
    \caption{Histograms of patch‑wise JSC values~\eqref{patchwise jsc} for different patch sizes under three noise levels: (a) $0.1\%$, (b) $1\%$, (c) $10\%$. The test case is T$\to$VB~\eqref{gttVB}. Orange and blue bars represent the histograms of patch‑wise JSC for inconsistent and consistent patches, respectively. Identification performance on consistent patches improves as patch size increases, whereas performance on inconsistent patches remains largely unchanged.}
    \label{JSC PATCH}
\end{figure}

Figure~\ref{JSC PATCH} shows the distribution of patch‑wise JSC values across different patch sizes and noise levels. As the patch size increases, the identification performance on consistent patches improves markedly. The JSC values on inconsistent patches also rise when more data are used, but a fraction of them remain $0$ even for the largest patch size ($400\times400$). The effect of noise on identification for a fixed patch size is relatively small, confirming the noise‑robustness of the identification algorithm~\cite{he2022robust}. It should be noted, however, that excessively large patches should be avoided because they would cause most patches to become inconsistent; see the discussion in Section~\ref{initial cover}.

\subsection{Validation of CUSUM for phase boundary detection}\label{other CPD} 
As introduced in Section~\ref{boundary estimation}, we employ the CUSUM method to detect the underlying phase boundary. We validate this choice by comparing it with two alternative change point detection approaches: the Gaussian kernel‑based method of Truong et al.~\cite{truong2020selective}, which evaluates splits of the sequence using a kernel similarity measure, and the Bayesian parametric method of Martínez and Mena~\cite{martinez2014nonparametric}, which models data on each side of a proposed split with Ornstein–Uhlenbeck processes~\cite{uhlenbeck1930theory}.

\begin{figure}
\centering
\begin{tabular}{c|c|c}
\toprule
T$\to$VB& KdV$\to$B&B$\to$T\\\midrule
\includegraphics[width=0.3\linewidth]{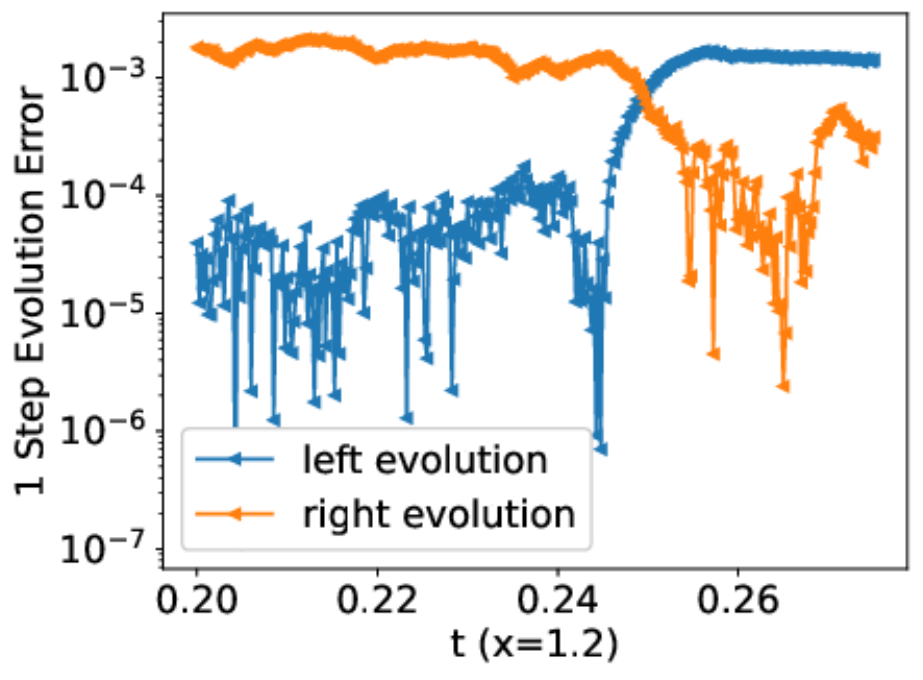}&
\includegraphics[width=0.3\linewidth]{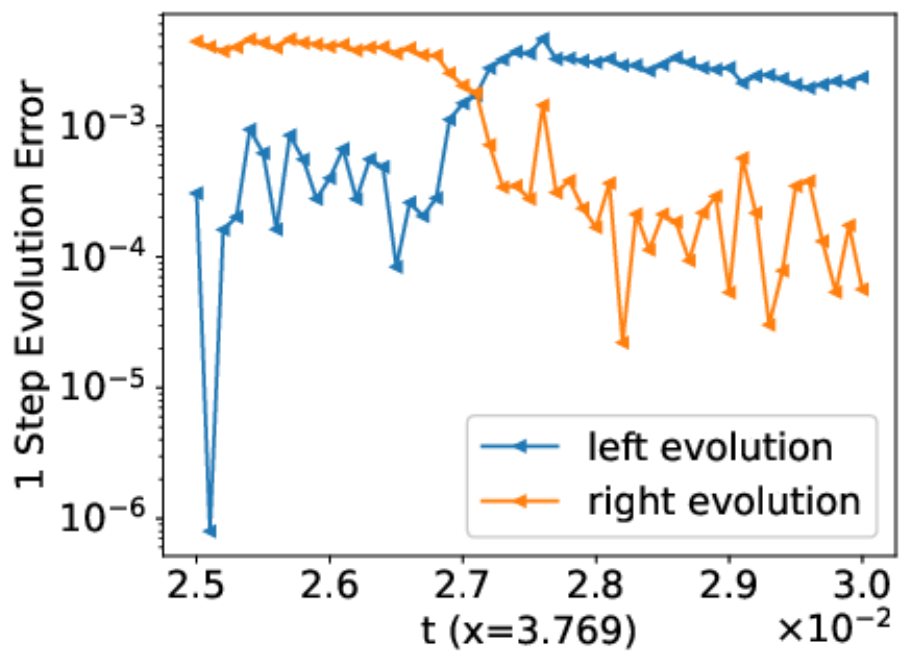}&
\includegraphics[width=0.3\linewidth]{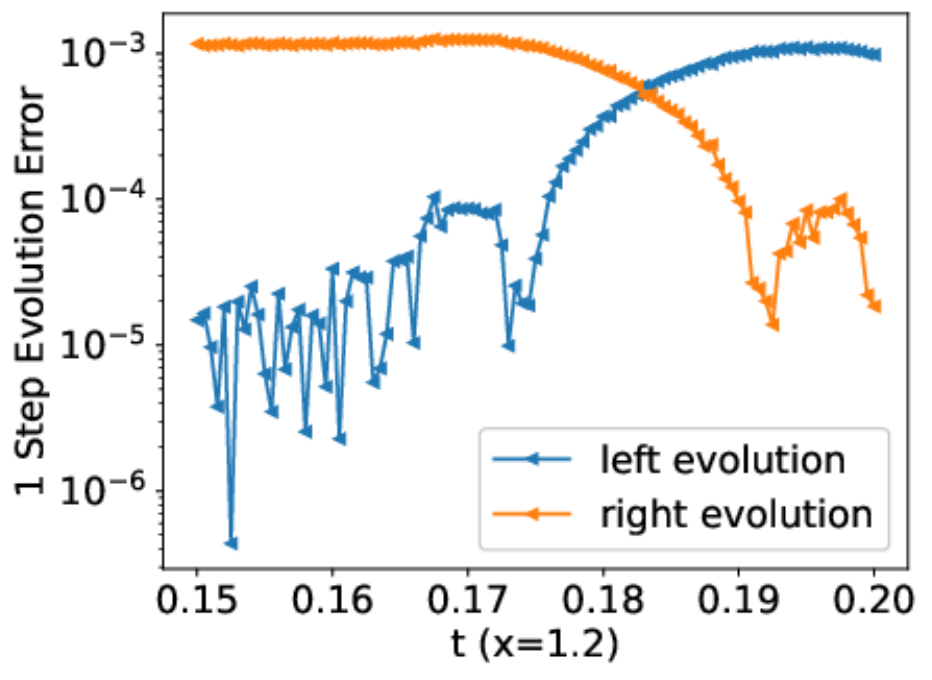}\\\midrule
\multicolumn{3}{c}{CUMSUM}\\\midrule
\includegraphics[width=0.3\linewidth]{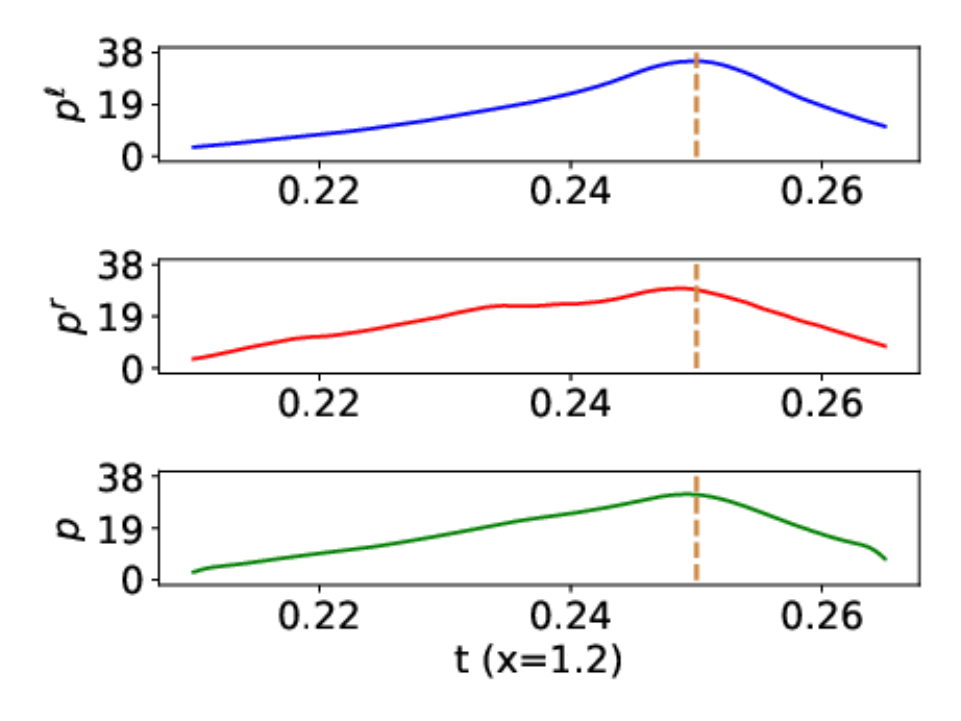}&
\includegraphics[width=0.3\linewidth]{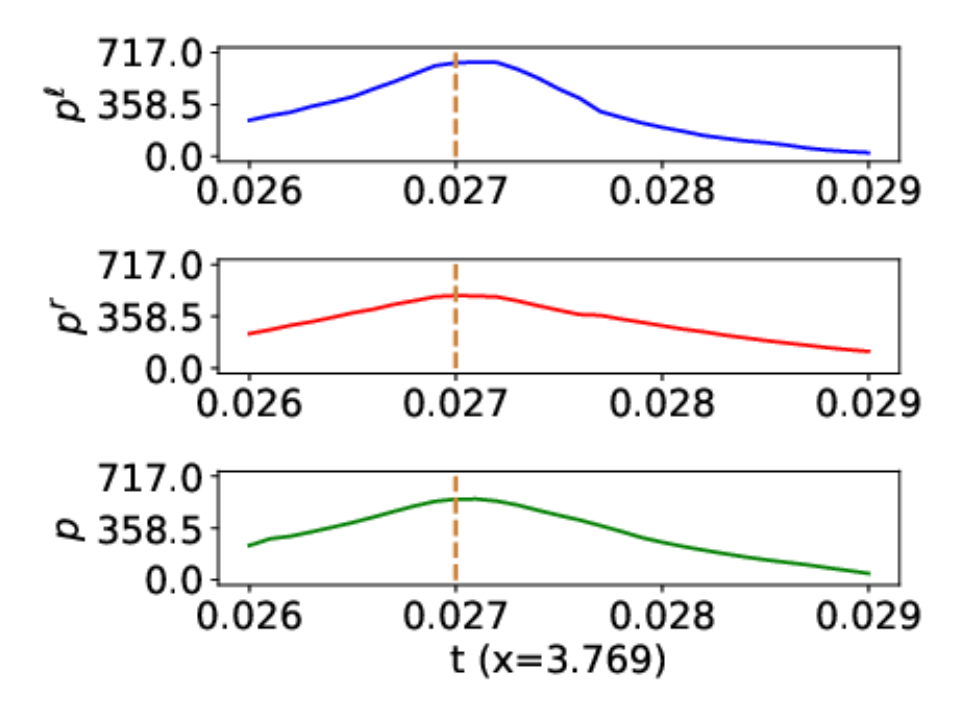}&
\includegraphics[width=0.3\linewidth]{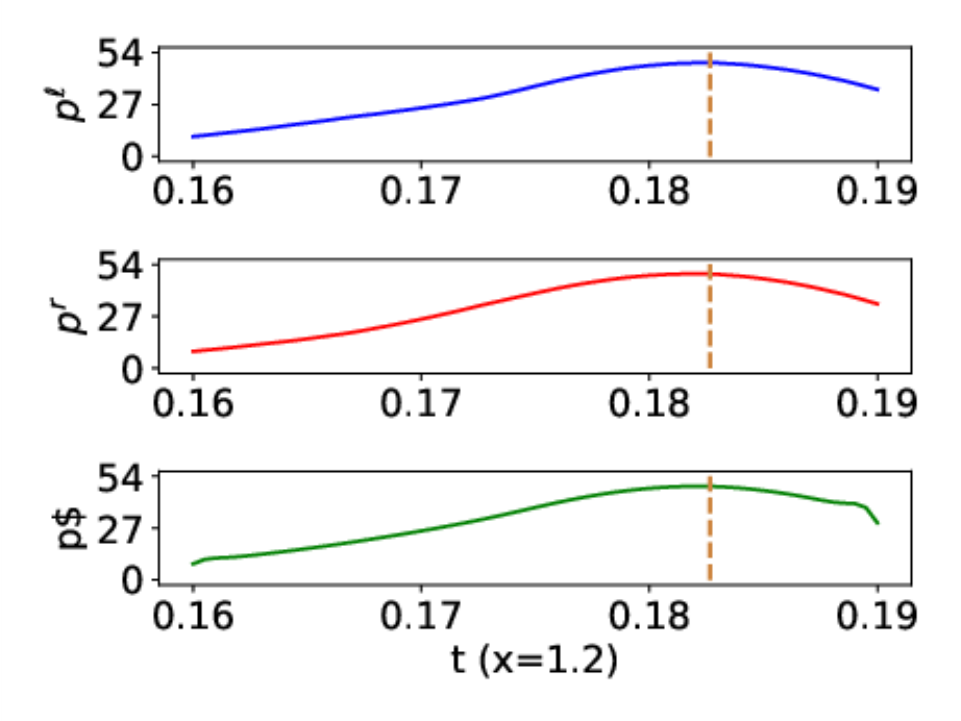}\\\midrule
\multicolumn{3}{c}{Gaussian}\\\midrule
\includegraphics[width=0.3\linewidth]{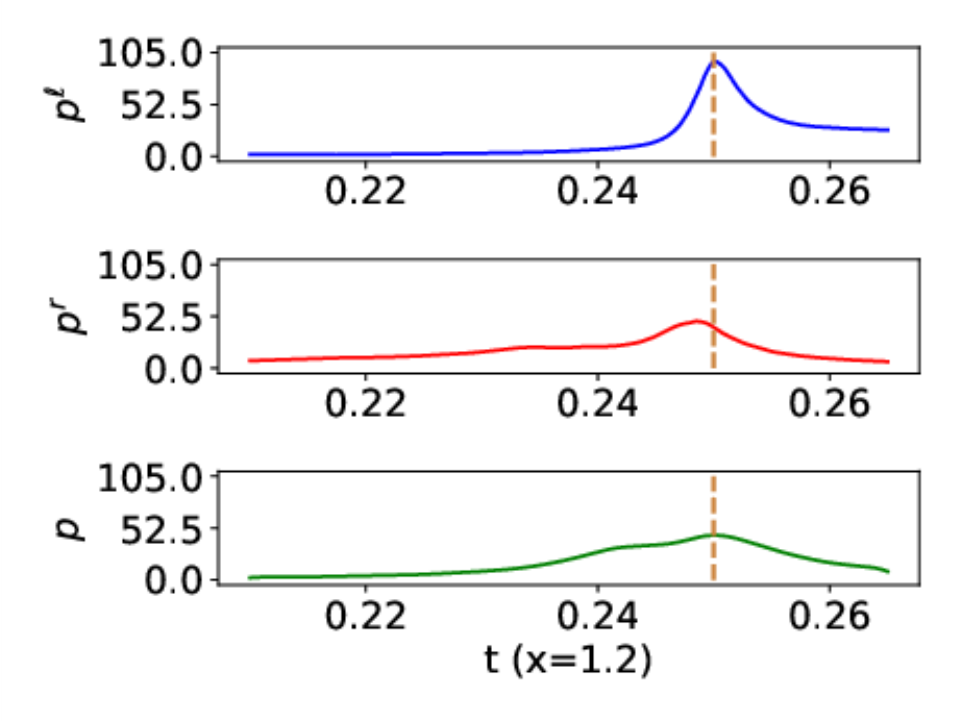}&
\includegraphics[width=0.3\linewidth]{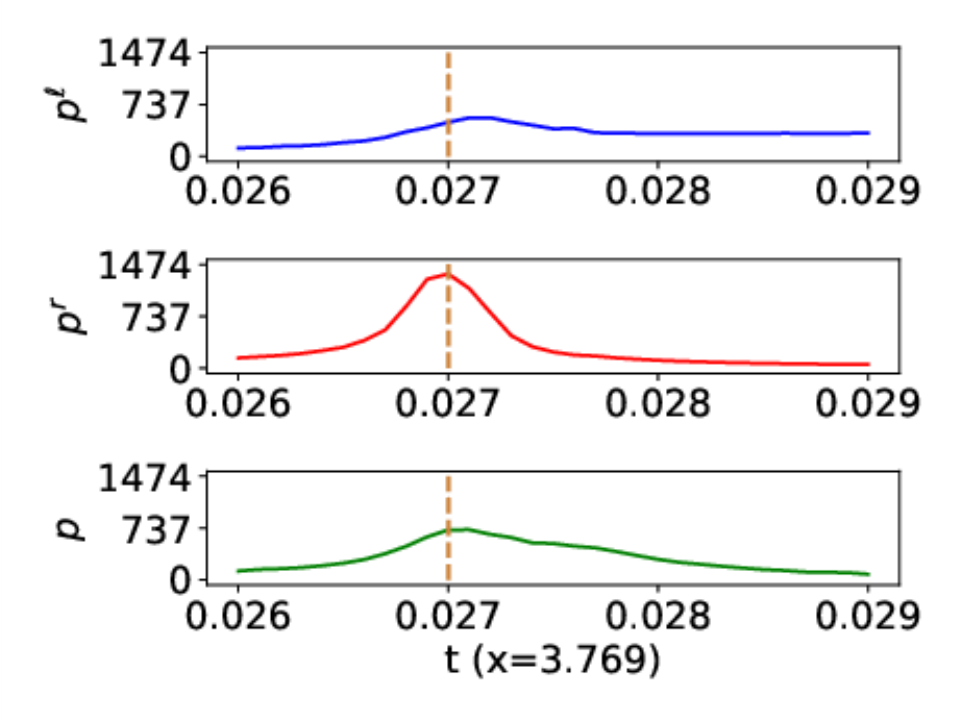}&
\includegraphics[width=0.3\linewidth]{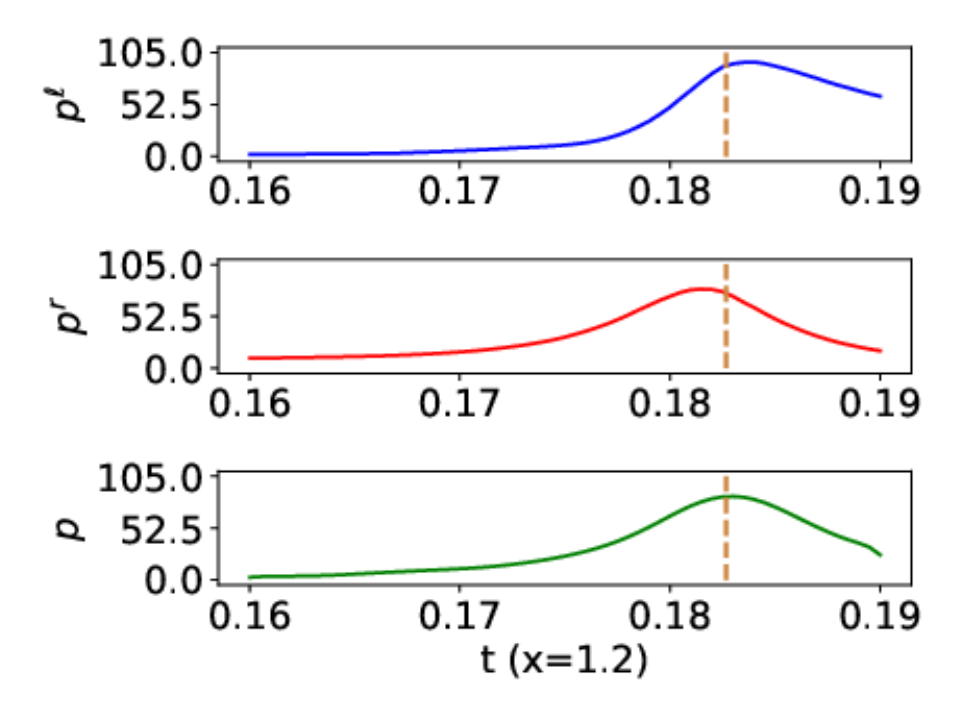}\\\midrule
\multicolumn{3}{c}{Bayesian}\\\midrule
\includegraphics[width=0.3\linewidth]{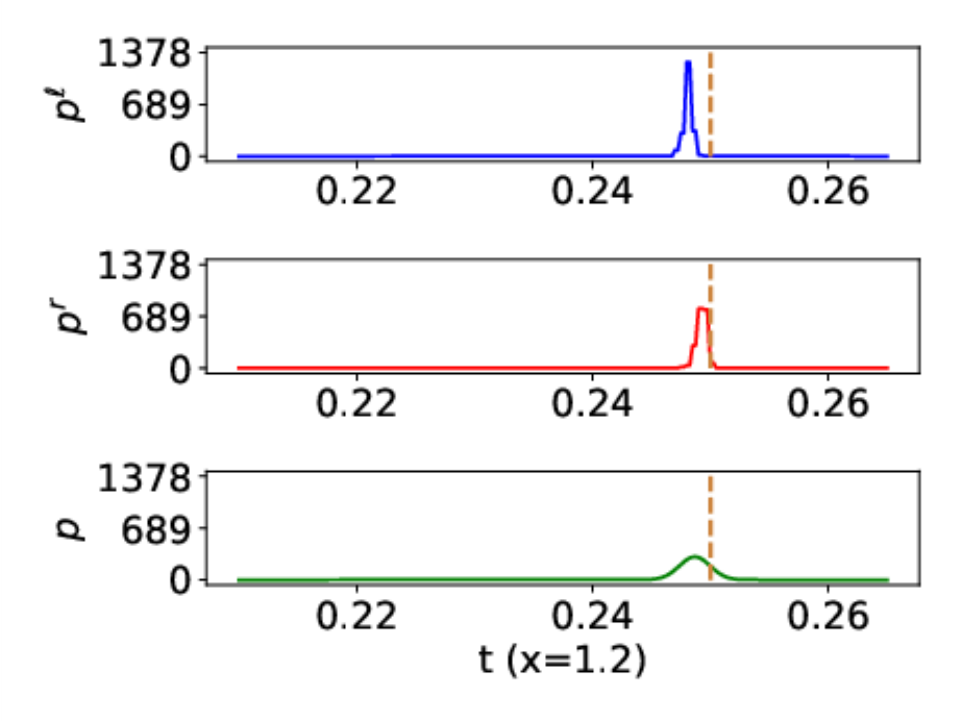}&
\includegraphics[width=0.3\linewidth]{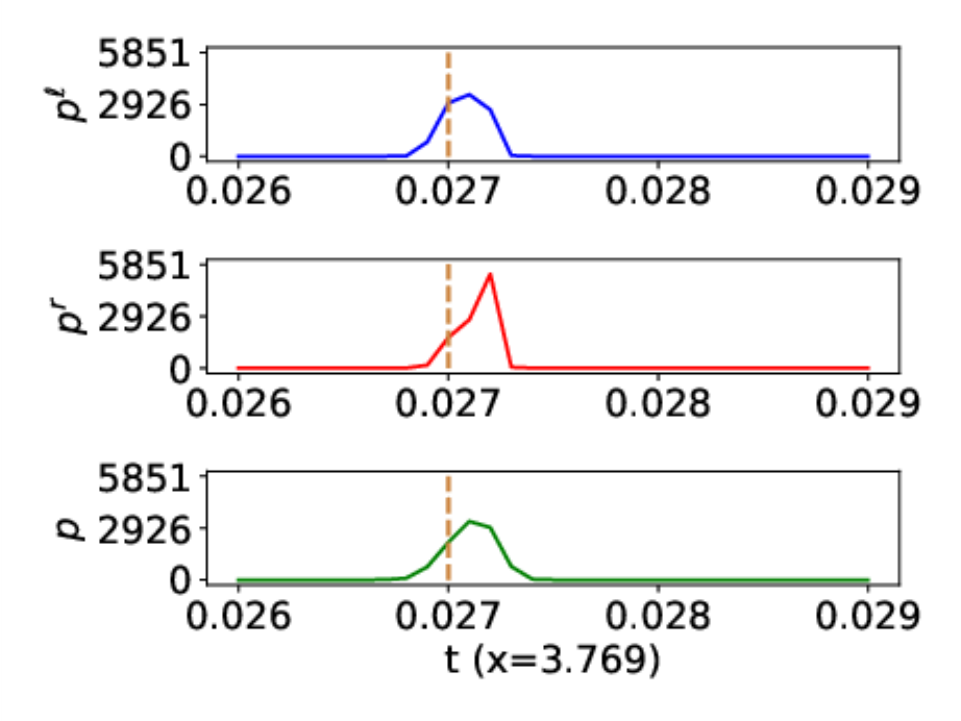}&
\includegraphics[width=0.3\linewidth]{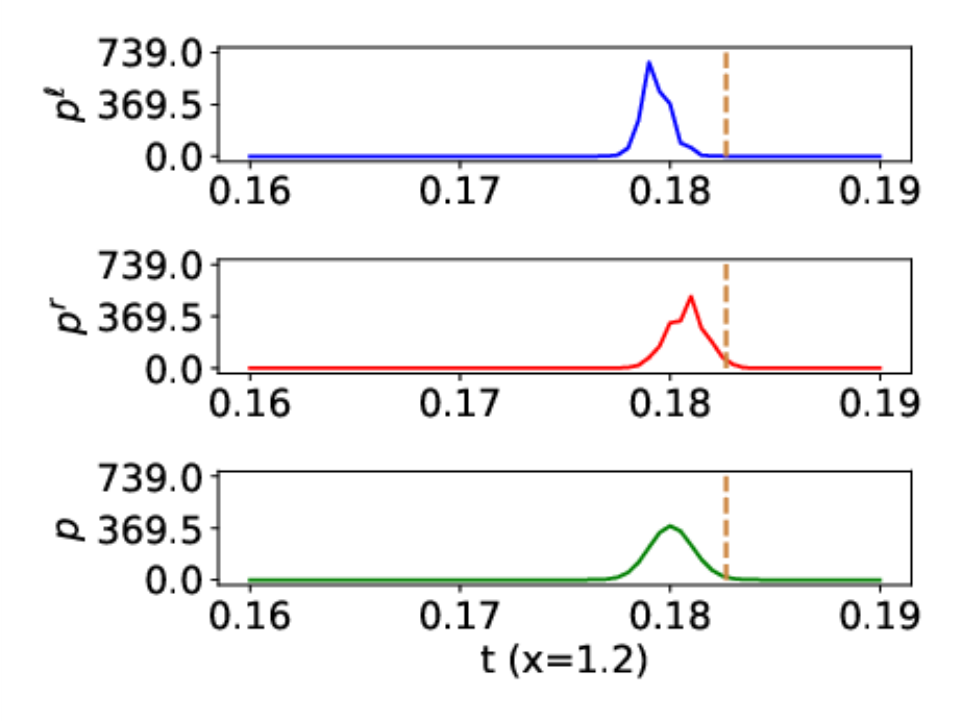}\\\bottomrule
\end{tabular}
\caption{Comparison of the change point detection methods: CUSUM (used in this paper), Gaussian kernel method~\cite{truong2020selective}, and Bayesian method~\cite{martinez2014nonparametric} tested on T$\to$VB~\eqref{gttVB}, KdV$\to$B~\eqref{gtKB} and B$\to$T~\eqref{gammacurve} under $1\%$ noise. The first row represents a pair of evolution error curves~\eqref{eq_sequences} of three different two-phase PDEs: T$\to$VB ($x=1.2$), KdV$\to$B ($x=3.769$) and B$\to$T ($x=1.2$). Each column shows respective change point detection results by three methods. The brown vertical line marks the precise location of the boundary point. The blue and red curves represent the results on $\be^\ell$ and $\be^r$, respectively. The green curve is the Wasserstein barycenter~\eqref{ot} of the two probability mass functions.}\label{changepoint detection}
\end{figure}

We compare the performances of CUSUM, the Gaussian kernel method, and the Bayesian method on the three test cases T$\to$VB~\eqref{gttVB}, KdV$\to$B~\eqref{gtKB}, and B$\to$T~\eqref{gammacurve}. As shown in Figure~\ref{changepoint detection}, the methods exhibit different characteristics and yield varying results. Across all cases, CUSUM consistently produces a sharp peak near the true change point, and the results from \(p^\ell_x\) and \(p^r_x\) are well aligned. The Gaussian kernel method also detects the change point accurately, but the consistency between its results on \(\be^\ell\) and \(\be^r\)~\eqref{eq_sequences} is  lower. The Bayesian model tends to concentrate the estimated change point distribution over a very narrow interval.

\begin{figure}
    \centering
    \begin{tabular}{c}
        \includegraphics[width=0.65\linewidth]{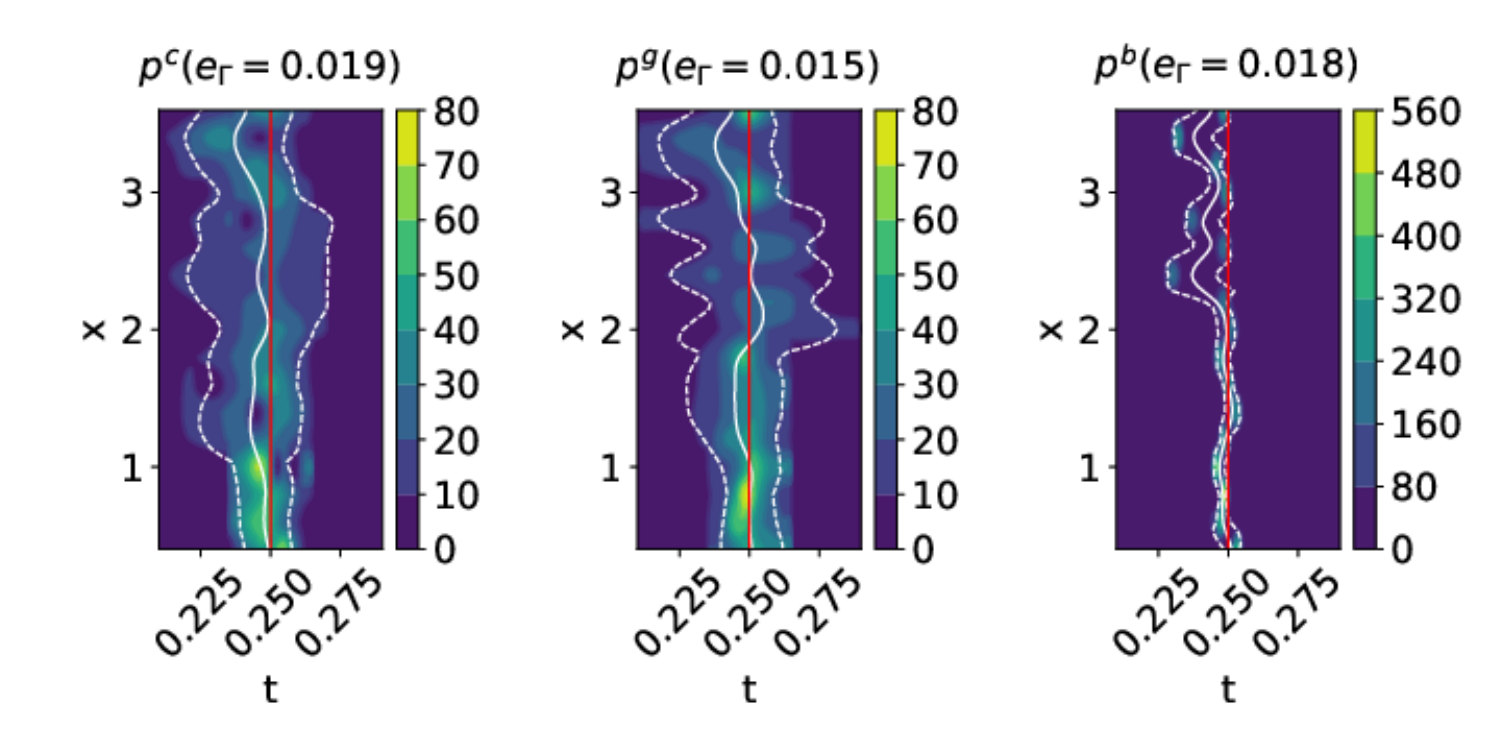}\\
        (a)\\
        \includegraphics[width=0.65\linewidth]{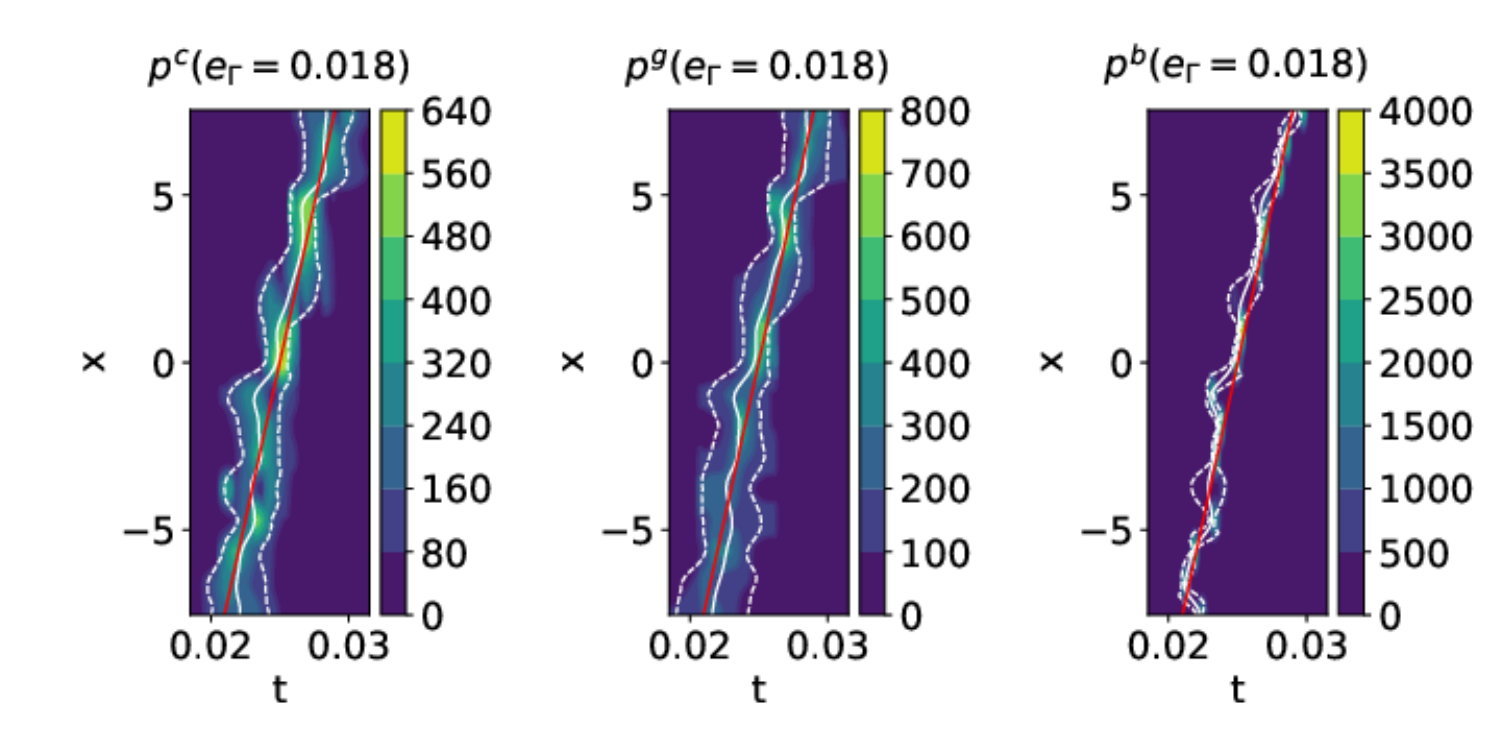}\\
        (b)\\
        \includegraphics[width=0.65\linewidth]{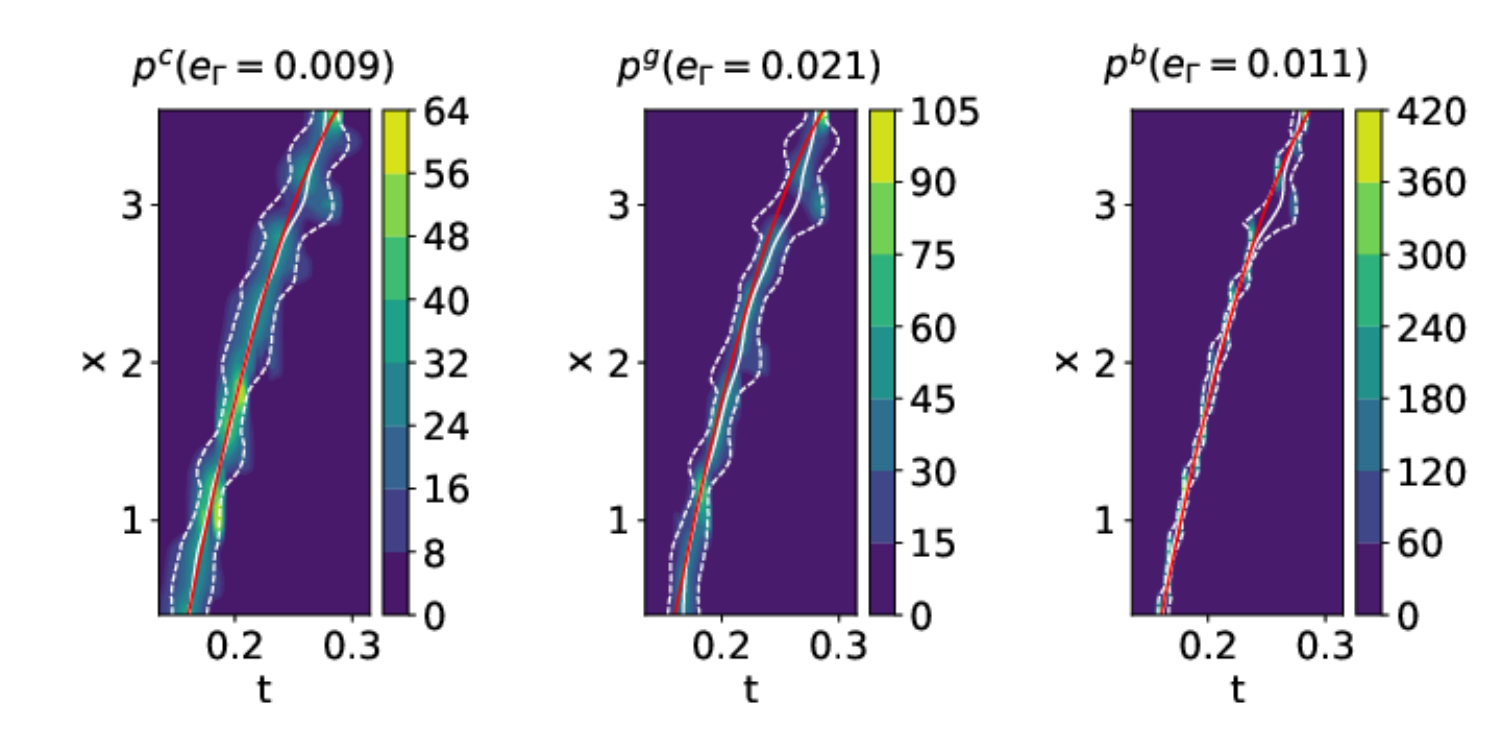}\\
        (c)
        
    \end{tabular}
    \caption{Phase boundary estimation by three different change point detection methods tested on three two-phase PDEs: (a) T$\to$VB~\eqref{gttVB}; (b) KdV$\to$B~\eqref{gtKB}; (c) B$\to$T~\eqref{gammacurve}. $p^c$, $p^g$ and $p^b$ are the approximate density of change point distribution given by CUSUM~\cite{robbins2011mean} (used in this work), Gaussian kernel-based method~\cite{truong2020selective} and Bayesian model~\cite{martinez2014nonparametric}, respectively. Compared with other two methods, Bayesian model tends to generate  narrower confidence interval. While Gaussian kernel-based method performs comparably to CUSUM on T$\to$VB and KdV$\to$B, it yields the largest error on B$\to$T.}
    \label{TB Distribution}
\end{figure}

In Figure~\ref{TB Distribution}, we display the estimated boundaries and confidence intervals produced by the three methods. The Bayesian model yields boundary estimates with the shortest confidence intervals, whereas the other two methods produce more moderate intervals. Although the Gaussian kernel method achieves the best boundary‑recovery accuracy when \(\Gamma\) is vertical (i.e., parallel to the time axis), its performance degrades when the slope of \(\Gamma\) increases. CUSUM demonstrates stable localization across all cases, and the confidence intervals it produces consistently cover the true phase boundary.

\section{Conclusion}\label{sec6}
In this paper, we propose a novel method, Phase‑IDENT, for identifying two‑phase PDE models and localizing the underlying phase boundary from noisy observational data. Our approach builds on and extends an existing PDE identification algorithm to recover models within each phase domain. By incorporating change point detection techniques, we simultaneously reconstruct the phase boundary and provide rigorous uncertainty quantification. We present comprehensive numerical experiments on simulated data, demonstrating the effectiveness and robustness of the proposed method. The uncertainty quantification we provide for detected phase boundaries offers valuable insight, enhances the interpretability of the results, and can facilitate future applications to real‑world data. Although this work focuses on two‑phase systems, the framework of Phase‑IDENT can be readily extended to handle multiple phases. Extensions to higher‑dimensional phase boundaries will be explored in future work.

\appendix
\section{Pseudo-code for Subspace Pursuit (SP)}\label{alg_SP}
For completeness, we present in Algorithm~\ref{alg:subspace_pursuit} the pseudo-code for the Subspace Pursuit~\cite{dai2009subspace}.
\begin{algorithm}
\caption{\textbf{:Subspace Pursuit SP}$(k; F, \mathbf{b})$}
\label{alg:subspace_pursuit}
\SetAlgoLined
\KwIn{$F \in \mathbb{R}^{NM \times K}$, $\mathbf{b} \in \mathbb{R}^{NM}$ and sparsity $k \in \mathbb{N}$.}

\BlankLine
\textbf{Initialization:} \\
$j \gets 0$ \\
$G \gets$ column-normalized version of $F$ \\
$\mathcal{I}^0 \gets \{k \text{ indices corresponding to the largest magnitude entries in the vector } G^*\mathbf{b}\}$ \\
$\mathbf{b}_{\text{res}}^0 \gets \mathbf{b} - G_{\mathcal{I}^0}G_{\mathcal{I}^0}^\dagger \mathbf{b}$ \\

\BlankLine
\While{True}{
    \textbf{Step 1:} 
    $\tilde{\mathcal{I}}^{j+1} \gets \mathcal{I}^j \cup \{k \text{ indices corresponding to the largest magnitude entries in the vector } G^*\mathbf{b}_{\text{res}}^j\}$ \\

    \textbf{Step 2:} 
    $\hat{\mathbf{c}} \gets G_{\tilde{\mathcal{I}}^{j+1}}^\dagger \mathbf{b}$ \\

    \textbf{Step 3:} 
    $\mathcal{I}^{j+1} \gets \{k \text{ indices corresponding to the largest elements of } \hat{\mathbf{c}}\}$ \\

    \textbf{Step 4:} 
    $\mathbf{b}_{\text{res}}^{j+1} \gets \mathbf{b} - G_{\mathcal{I}^{j+1}}G_{\mathcal{I}^{j+1}}^\dagger \mathbf{b}$ \\

    \textbf{Step 5:} 
    \If{$\|\mathbf{b}_{\text{res}}^{j+1}\|_2 > \|\mathbf{b}_{\text{res}}^j\|_2 $}{
        $\mathcal{I}^{j+1} \gets \mathcal{I}^j$ \\
        \textbf{Terminate the algorithm}
    }
    \Else{
        $j \gets j + 1$ \\
        \textbf{Iterate}
    }
}

\KwOut{$\hat{\mathbf{c}} \in \mathbb{R}^K$ satisfying $\hat{\mathbf{c}}_{\mathcal{I}_j} = F_{\mathcal{I}_j}^\dagger \mathbf{b}$ and $\hat{\mathbf{c}}_{(\mathcal{I}_j)^C} = 0$, here $(\mathcal{I}_j)^C$ denotes the complement of $\mathcal{I}_j$.}
\end{algorithm}

\section{Single versus multiple step evolution for change point detection}\label{multistep}
\begin{figure}
    \centering
    \begin{tabular}{c@{\vspace{2pt}}c@{\vspace{2pt}}c@{\vspace{2pt}}c}

    \includegraphics[width=0.24\textwidth]{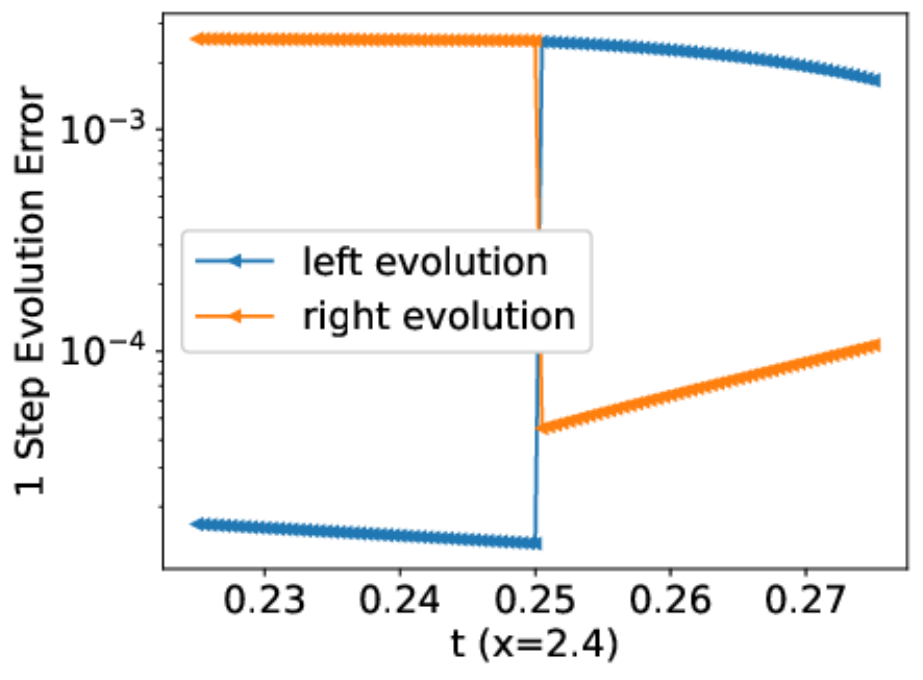}&
    \includegraphics[width=0.24\linewidth]{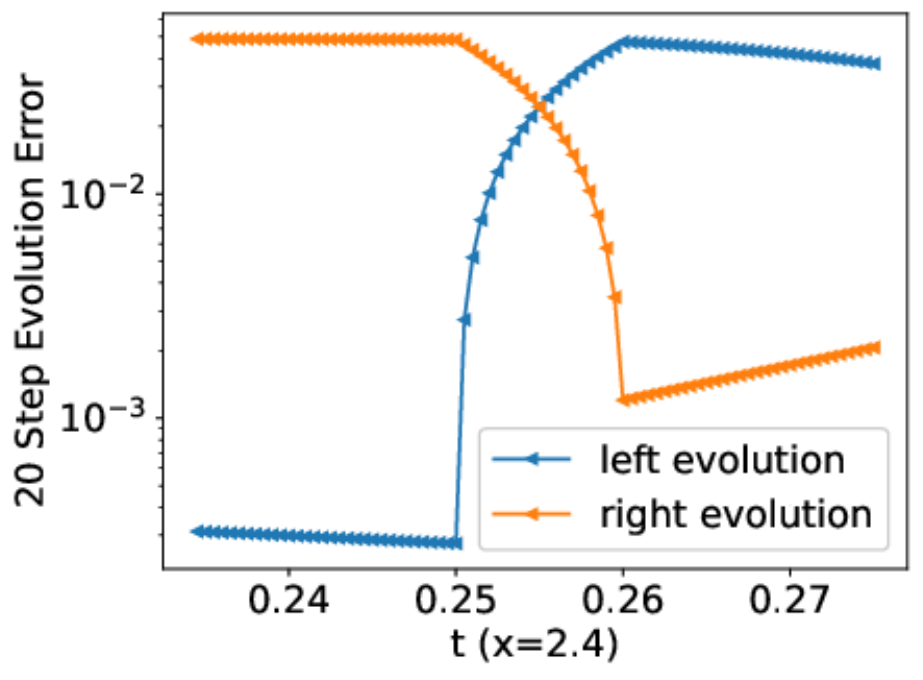}&
    \includegraphics[width=0.24\textwidth]{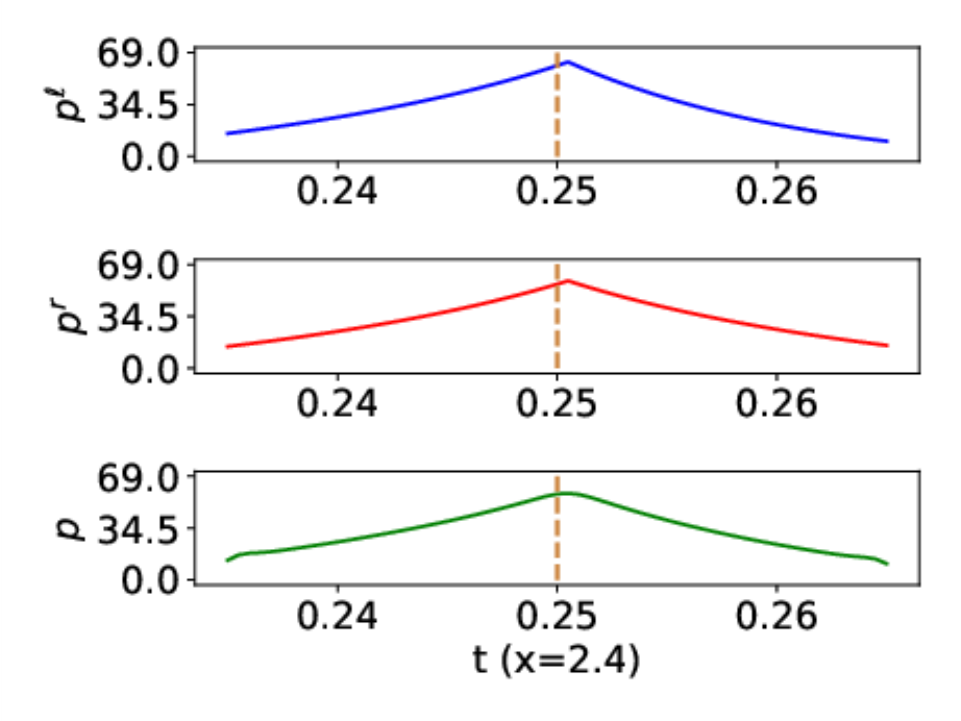}&
    \includegraphics[width=0.24\linewidth]{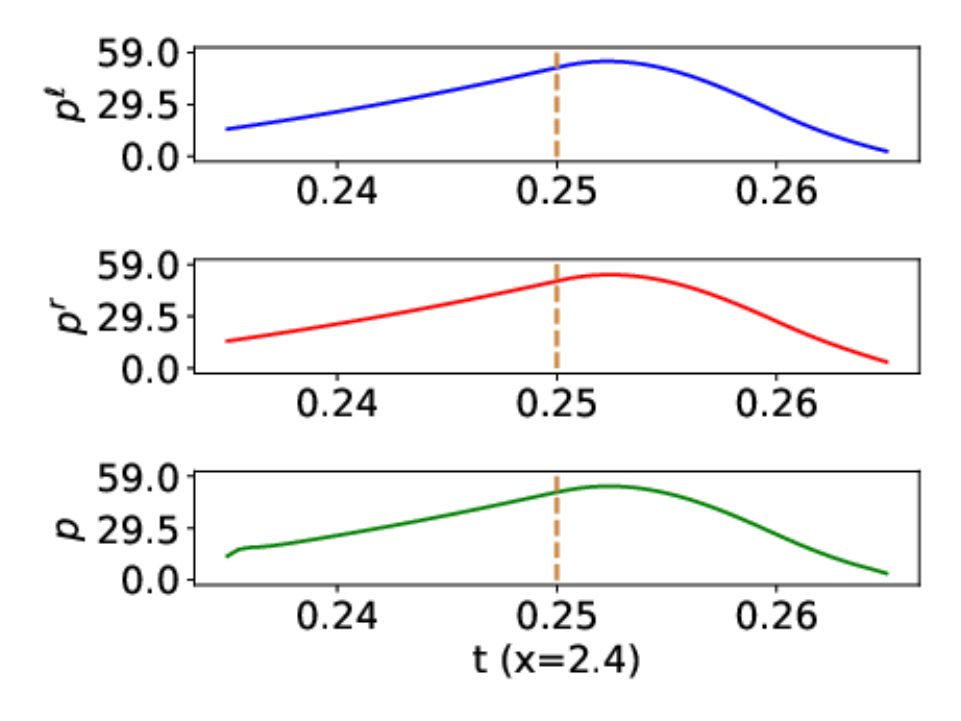}\\
    (a)&(b)&(c)&(d)
    \end{tabular}
    \caption{Comparison between single‑step and multi‑step evolution for the B$\to$T system~\eqref{eq_effect_phase_boundary_slope} with \(s = 0\). (a) Single‑step evolution errors at \(x = 2.4\); (b) $20$‑step evolution errors at \(x = 2.4\); (c) probability mass functions of the change point distribution and their Wasserstein barycenter obtained from single‑step evolution errors; (d) probability mass functions and their barycenter obtained from $20$‑step evolution errors.}
    \label{multievo}
\end{figure}

We compare the evolution errors~\eqref{eq_sequences} and the change point detection results obtained using single‑step versus $20$‑step evolution. The test case is the two‑phase system~\eqref{eq_effect_phase_boundary_slope} with slope \(s = 0\). To clearly highlight differences between the two evolution schemes, we use noise‑free data. Figure~\ref{multievo}(a) and (b) show the error sequences produced by single‑step and $20$‑step evolution, respectively; Figure~\ref{multievo}(c) and (d) display the corresponding change point probability mass functions (defined in~\eqref{pdelta}) and their respective Wasserstein barycenters~\eqref{ot}. We observe that single‑step evolution yields a more accurate change point location and produces probability mass functions with sharper, more distinguishable peaks.
\section{Proof of Proposition~\ref{proposition}}\label{proof of propo}
\begin{proof}

Denote $\kappa:=(\gamma(x)-\tau_m)/(\tau_{m+1}-\tau_m)$. By the one-sided Taylor expansion, we have
\begin{equation}
\begin{aligned}
u(x,\tau_{m+1}) & =u(x,\gamma(x))+(1-\kappa) \Delta t \partial_t^+u (x,\gamma(x))+\cO\left((1-\kappa)^2\Delta t^2\right)
\end{aligned}
\end{equation}
and
\begin{equation}
\begin{aligned}
u(x,\tau_{m}) & =u(x,\gamma(x))-\kappa\Delta t \partial_t^-u (x,\gamma(x))+\cO\left(\kappa^2\Delta t^2\right)
\end{aligned}
\end{equation}
Furthermore, by the assumption on the scheme's consistency~\eqref{assumption},
\begin{equation}
\begin{aligned}
\widehat{\mathcal{F}}_\ell(\bm U_m)&=u_t(x,\tau_{m})+\cO(\Delta x^q)\\
&= \partial^-_tu(x,\gamma(x)) + \cO(\Delta t)+\cO(\Delta x^q).
\end{aligned}
\end{equation}
Hence,
\begin{equation}
\begin{aligned}
e^\ell(m,x) & =u(x,\tau_{m+1})-u(x,\tau_{m})-\Delta t\widehat{\mathcal{F}}_\ell(\bm U_m)+\cO(\Delta x^p)\\
&=(1-\kappa)\Delta t\partial_t^+u(x,\gamma(x))+\kappa\Delta t\partial_t^-u(x,\gamma(x))-\Delta t\widehat{\mathcal{F}}_\ell(\bm U_m)+\cO(\Delta x^p)\\
&=(1-\kappa)\Delta t\partial_t^+u(x,\gamma(x))+\kappa\Delta t\partial_t^-u(x,\gamma(x))-\Delta t(\partial^-_tu(x,\gamma(x)) + \cO(\Delta t)+\cO(\Delta x^q))+\cO(\Delta x^p)\\
&=(1-\kappa)\Delta t \left(\partial_t^+u(x,\gamma(x))-\partial_t^-u(x,\gamma(x))\right)+\cO(\Delta t^2)+\cO(\Delta t\Delta x^q)+\cO(\Delta x^p)\\
&=\left(\partial_t^+u(x,\gamma(x))-\partial_t^-u(x,\gamma(x))\right)\cO(\Delta t)+\cO(\Delta x^p),
\end{aligned}
\end{equation}where the last equality holds since $\cO(\Delta t\Delta x^q)$ can be absorbed by $\cO(\Delta t)$.
\end{proof}

\bibliographystyle{siam}
\bibliography{main_new2}
\end{document}